\documentclass[a4paper,11pt]{article}

\usepackage{amsfonts,amssymb,amsmath,amsthm,latexsym}
\usepackage{mathrsfs}
\usepackage{graphicx,graphics}
\usepackage{hyperref}

\usepackage[english]{babel}
\makeatletter
\def\@captype{figure}
\makeatother

 \makeatletter
 \@addtoreset{equation}{section}
 \makeatother

 \makeatletter
 \@addtoreset{enunciato}{section}
 \makeatother

 \newcounter{enunciato}[section]

 \newtheorem{ittheorem}{Theorem}
 \newtheorem{itcorollary}{Corollary}
 \newtheorem{itlemma}{Lemma}
 \newtheorem{itproposition}{Proposition}
 \newtheorem{itdefinition}{Definition}
 \newtheorem{itremark}{Remark}
 \newtheorem{itclaim}{Claim}
 \newtheorem{itfact}{Fact}
 \newtheorem{itconjecture}{Conjecture}

 \newenvironment{theorem}{\addtocounter{enunciato}{1}
 \begin{ittheorem}}{\end{ittheorem}}
 
 \newenvironment{corollary}{\addtocounter{enunciato}{1}
 \begin{itcorollary}}{\end{itcorollary}}

 \newenvironment{lemma}{\addtocounter{enunciato}{1}
 \begin{itlemma}}{\end{itlemma}}

 \newenvironment{proposition}{\addtocounter{enunciato}{1}
 \begin{itproposition}}{\end{itproposition}}

 \newenvironment{definition}{\addtocounter{enunciato}{1}
 \begin{itdefinition}}{\end{itdefinition}}

 \newenvironment{remark}{\addtocounter{enunciato}{1}
 \begin{itremark}}{\end{itremark}}

 \newenvironment{claim}{\addtocounter{enunciato}{1}
 \begin{itclaim}}{\end{itclaim}}

 \newenvironment{fact}{\addtocounter{enunciato}{1}
 \begin{itfact}}{\end{itfact}}

 \newenvironment{conjecture}{\addtocounter{enunciato}{1}
 \begin{itconjecture}}{\end{itconjecture}}

 \newcommand{\be}[1]{\begin{equation}\label{#1}}
 \newcommand{\ee}{\end{equation}}

 \newcommand{\bl}[1]{\begin{lemma}\label{#1}}
 \newcommand{\el}{\end{lemma}}

 \newcommand{\br}[1]{\begin{remark}\label{#1}}
 \newcommand{\er}{\end{remark}}

 \newcommand{\bt}[1]{\begin{theorem}\label{#1}}
 \newcommand{\et}{\end{theorem}}

 \newcommand{\bd}[1]{\begin{definition}\label{#1}}
 \newcommand{\ed}{\end{definition}}

 \newcommand{\bcl}[1]{\begin{claim}\label{#1}}
 \newcommand{\ecl}{\end{claim}}

 \newcommand{\bfact}[1]{\begin{fact}\label{#1}}
 \newcommand{\efact}{\end{fact}}

 \newcommand{\bp}[1]{\begin{proposition}\label{#1}}
 \newcommand{\ep}{\end{proposition}}

 \newcommand{\bc}[1]{\begin{corollary}\label{#1}}
 \newcommand{\ec}{\end{corollary}}

 \newcommand{\bcj}[1]{\begin{conjecture}\label{#1}}
 \newcommand{\ecj}{\end{conjecture}}

 \newcommand{\bpr}{\begin{proof}}
 \newcommand{\epr}{\end{proof}}

 \newcommand{\bprt}[1]{\begin{proofoft}{\it\ref{#1}}.\,\,}
 \newcommand{\eprt}{\end{proofoft}}
 
 \newcommand{\bprc}[1]{\begin{proofofc}{\it\ref{#1}}.\,\,}
 \newcommand{\eprc}{\end{proofofc}}
 
 \newcommand{\bprp}[1]{\begin{proofofp}{\it\ref{#1}}.\,\,}
 \newcommand{\eprp}{\end{proofofp}}
 
 \newcommand{\bprl}[1]{\begin{proofofl}{\it\ref{#1}}.\,\,}
 \newcommand{\eprl}{\end{proofofl}}

 \newcommand{\bi}{\begin{itemize}}
 \newcommand{\ei}{\end{itemize}}

 \newcommand{\ben}{\begin{enumerate}}
 \newcommand{\een}{\end{enumerate}}
 
 \newcommand{\ess}{ess \,sup}

\newcommand{\argmax}{\operatornamewithlimits{argmax}}

\newcommand{\dist}{\operatornamewithlimits{dist}}

\newenvironment{proofoft}{\noindent 
{\bf Proof of Theorem\,}}{\hspace*{\fill}$\halmos$\medskip}
\newenvironment{proofofc}{\noindent 
{\bf Proof of Corollary\,}}{\hspace*{\fill}$\halmos$\medskip}
\newenvironment{proofofp}{\noinde
nt 
{\bf Proof of Proposition\,}}{\hspace*{\fill}$\halmos$\medskip}
\newenvironment{proofofl}{\noindent 
{\bf Proof of Lemma\,}}{\hspace*{\fill}$\halmos$\medskip}
\newcommand{\halmos}{\rule{1ex}{1.4ex}}

\def \da {\downarrow}

\def \PP {{\mathbb P}}
\def \EE {{\mathbb E}}

\def \cc {{\rm c}}

\def \E {{\mathbb E}}
\def \N {{\mathbb N}}
\def \P {{\mathbb P}}
\def \Q {{\mathbb Q}}
\def \R {{\mathbb R}}
\def \Z {{\mathbb Z}}

\def \cI {{\mathcal I}}

\def \cF {{\mathcal F}}

\newcommand{\one}{{\mathchoice {1\mskip-4mu\mathrm l}
         {1\mskip-4mu\mathrm l}
         {1\mskip-4.5mu\mathrm l}
         {1\mskip-5mu\mathrm l}}}

\begin{document}

\title{The parabolic Anderson model in a\\ 
dynamic random environment: basic properties of\\ 
the quenched Lyapunov exponent}

\author{\renewcommand{\thefootnote}{\arabic{footnote}}
D.\ Erhard
\footnotemark[1]
\\
\renewcommand{\thefootnote}{\arabic{footnote}}
F.\ den Hollander
\footnotemark[2]
\\
\renewcommand{\thefootnote}{\arabic{footnote}}
G.\ Maillard
\footnotemark[3]
}

\footnotetext[1]{
Mathematical Institute, Leiden University, P.O.\ Box 9512,
2300 RA Leiden, The Netherlands,\\
{\sl erhardd@math.leidenuniv.nl}
}
\footnotetext[2]{
Mathematical Institute, Leiden University, P.O.\ Box 9512,
2300 RA Leiden, The Netherlands,\\
{\sl denholla@math.leidenuniv.nl}
}
\footnotetext[3]{
CMI-LATP, Aix-Marseille Universit\'e,
39 rue F. Joliot-Curie, F-13453 Marseille Cedex 13, France,\\
{\sl maillard@cmi.univ-mrs.fr}
}

\date{\today}
\maketitle


\begin{abstract}

In this paper we study the parabolic Anderson equation $\partial u(x,t)/\partial t 
= \kappa\Delta u(x,t) + \xi(x,t)u(x,t)$, $x\in\Z^d$, $t\geq 0$, where the $u$-field 
and the $\xi$-field are $\R$-valued, $\kappa \in [0,\infty)$ is the diffusion 
constant, and $\Delta$ is the discrete Laplacian. The $\xi$-field plays the role of 
a \emph{dynamic random environment} that drives the equation. The initial condition 
$u(x,0)=u_0(x)$, $x\in\Z^d$, is taken to be non-negative and bounded. The solution 
of the parabolic Anderson equation describes the evolution of a field of particles 
performing independent simple random walks with binary branching: particles jump at 
rate $2d\kappa$, split into two at rate $\xi \vee 0$, and die at rate $(-\xi) \vee 0$. 
Our goal is to prove a number of \emph{basic properties} of the solution $u$ under 
assumptions on $\xi$ that are as weak as possible.
These properties will serve as a jump board for later refinements.
 
Throughout the paper we assume that $\xi$ is stationary and ergodic under 
translations in space and time, is not constant and satisfies $\E(|\xi(0,0)|)
<\infty$, where $\E$ denotes expectation w.r.t.\ $\xi$. Under a mild assumption 
on the tails of the distribution of $\xi$, we show that the solution to the 
parabolic Anderson equation exists and is unique for all $\kappa\in [0,\infty)$. 
Our main object of interest is the \emph{quenched Lyapunov exponent} $\lambda_0
(\kappa) = \lim_{t\to\infty} \frac{1}{t}\log u(0,t)$. It was shown in G\"artner, 
den Hollander and Maillard~\cite{GdHM11} that this exponent exists and is constant 
$\xi$-a.s., satisfies $\lambda_0(0) = \E(\xi(0,0))$ and $\lambda_0(\kappa) 
> \E(\xi(0,0))$ for $\kappa \in (0,\infty)$, and is such that $\kappa\mapsto
\lambda_0(\kappa)$ is globally Lipschitz on $(0,\infty)$ outside any neighborhood 
of $0$ where it is finite. Under certain weak space-time mixing assumptions on 
$\xi$, we show the following properties: (1) $\lambda_0(\kappa)$ does not depend 
on the initial condition $u_0$; (2) $\lambda_0(\kappa)<\infty$ for all $\kappa\in
[0,\infty)$; (3) $\kappa\mapsto\lambda_0(\kappa)$ is continuous on $[0,\infty)$ 
but not Lipschitz at $0$. We further conjecture: (4) $\lim_{\kappa\to\infty} 
[\lambda_p(\kappa) - \lambda_0(\kappa)]=0$ for all $p \in \N$, where $\lambda_p
(\kappa)=\lim_{t\to\infty} \frac{1}{pt}\log\E([u(0,t)]^p)$ is the $p$-th 
\emph{annealed Lyapunov exponent}. (In \cite{GdHM11} properties (1), (2) and (4) 
were not addressed, while property (3) was shown under much more restrictive 
assumptions on $\xi$.) Finally, we prove that our weak space-time mixing conditions 
on $\xi$ are satisfied for several classes of interacting particle systems. 

\newpage

{\it Acknowledgment.} DE and FdH were supported by ERC Advanced Grant 267356 VARIS. 
GM was supported by the CNRS while on sabbatical leave at EURANDOM, Eindhoven,
The Netherlands, during the academic year 2010--2011. 

\end{abstract}


\section{Introduction and main results}
\label{S1}

Section~\ref{S1.1} defines the parabolic Anderson model and provides motivation, 
Section~\ref{S1.2} describes our main targets and their relation to the literature, 
Section~\ref{S1.3} contains our main results, while Section~\ref{S1.4} 
discusses these results and state a conjecture.


\subsection{The parabolic Anderson model (PAM)}
\label{S1.1}

The parabolic Anderson model is the partial differential equation
\begin{equation}
\label{pA}
\frac{\partial}{\partial t}u(x,t) = \kappa\Delta u(x,t) + \xi(x,t)u(x,t),
\qquad x\in\Z^d,\,t\geq 0.
\end{equation}
Here, the $u$-field is $\R$-valued, $\kappa\in [0,\infty)$ is the diffusion
constant, $\Delta$ is the discrete Laplacian acting on $u$ as
\begin{equation}
\label{dL}
\Delta u(x,t) = \sum_{{y\in\Z^d} \atop {\|y-x\|=1}} [u(y,t)-u(x,t)]
\end{equation}
($\|\cdot\|$ is the $l_1$-norm), while
\begin{equation}
\label{rf}
\xi = (\xi_t)_{t \geq 0} \mbox{ with } \xi_t = \{\xi(x,t) \colon\,x\in\Z^d\}
\end{equation}
is an $\R$-valued random field playing the role a of \emph{dynamic random environment}
that drives the equation. As initial condition for (\ref{pA}) we take 
\begin{equation}
\label{ic}
\blacktriangleright\quad u(x,0) = u_0(x),\,x\in\Z^d, 
\mbox{ with } u_0  \mbox{ non-negative and bounded}.
\end{equation}

One interpretation of (\ref{pA}) and (\ref{ic}) comes from \emph{population dynamics}. 
Consider the special case where $\xi(x,t) = \gamma\bar\xi(x,t)-\delta$ with
$\delta,\gamma \in (0,\infty)$ and $\bar\xi$ an $\N_0$-valued random field. 
Consider a system of two types of particles, $A$ (catalyst) and $B$ (reactant), 
subject to:
\begin{itemize}
\item[--]
$A$-particles evolve autonomously according to a prescribed dynamics with 
$\bar\xi(x,t)$ denoting the number of $A$-particles at site $x$ at time $t$;
\item[--]
$B$-particles perform independent simple random walks at rate $2d\kappa$ 
and split into two at a rate that is equal to $\gamma$ times the number of 
$A$-particles present at the same location at the same time;
\item[--]
$B$-particles die at rate $\delta$;
\item[--]
the average number of $B$-particles at site $x$ at time $0$ is $u_{0}(x)$.
\end{itemize}
Then
\begin{equation}
\label{uint}
\begin{array}{lll}
u(x,t) &=& \hbox{the average number of $B$-particles at site $x$ at time $t$}\\
       && \hbox{conditioned on the evolution of the $A$-particles}.
\end{array}
\end{equation}

The $\xi$-field is defined on a probability space $(\Omega,\cF,\PP)$. Throughout the 
paper we assume that 
\begin{equation}
\label{staterg}
\begin{aligned}
&\blacktriangleright\quad \xi \mbox{ is \emph{stationary} and \emph{ergodic} under 
translations in space and time.}\\
&\blacktriangleright\quad \xi \mbox{ is \emph{not constant} and } \E(|\xi(0,0)|)<\infty.
\end{aligned}
\end{equation}
Without loss of generality we may assume that $\E(\xi(0,0)) = 0$. 


\subsection{Main targets and related literature}
\label{S1.2}

The goal of the present paper is to prove a number of \emph{basic properties} about the 
Cauchy problem in (\ref{pA}) with initial condition (\ref{ic}). In this section we 
describe these properties informally. Precise results will be stated in Section~\ref{S1.3}.

\medskip\noindent
$\bullet$ {\bf Existence and uniqueness of the solution.}
For \emph{static} $\xi$, i.e.,
\begin{equation}
\label{timeindeprf}
\xi = \{\xi(x)\colon\,x \in \Z^d\},
\end{equation}
existence and uniqueness of the solution to \eqref{pA} with initial condition \eqref{ic} 
were addressed by G\"artner and Molchanov~\cite{GM90}. Namely, for arbitrary, deterministic $q\colon\,
\Z^d \to \R$ and $u_0\colon\,\Z^d \to [0,\infty)$, they considered the equation
\begin{equation}
\label{timeindeppA}
\begin{cases}
\frac{\partial}{\partial t}u(x,t) = \kappa\Delta u(x,t) + q(x)u(x,t),\\
u(x,0) = u_0(x),
\end{cases}
\qquad x\in\Z^d,\,t\geq 0,
\end{equation}
with $u_0$ non-negative, and showed that there exists a non-negative solution if 
and only if the Feynman-Kac formula
\begin{equation}
\label{FKtimeindep}
v(x,t) = E_x\left(\exp\left\{\int_0^t 
q(X^\kappa(s))\,ds\right\}\,u_0(X^\kappa(t))\right)
\end{equation}
is finite for all $x$ and $t$. Here, $X^\kappa=(X^\kappa(t))_{t \geq 0}$ is the 
continuous-time simple random walk jumping at rate $2d\kappa$ (i.e., the Markov 
process with generator $\kappa\Delta$) starting in $x$ under the law $P_{x}$. 
Moreover, they showed that $v$ in \eqref{FKtimeindep} is the minimal non-negative 
solution to (\ref{timeindeppA}). From these considerations they deduced a criterion 
for the almost sure existence of a solution to equation (\ref{timeindeppA}) when $q=\xi$. This result was later extended to \emph{dynamic}
$\xi$ by Carmona and Molchanov~\cite{CM94}, who proved the following.

\bp{existtimedeppA}
{\rm (Carmona and Molchanov \cite{CM94})}
Suppose that $q\colon\,\Z^d \times [0,\infty) \to \R$ is such that $q(x,\cdot)$ is 
locally integrable for every $x$. Then, for every non-negative initial condition 
$u_0$, the deterministic equation
\begin{equation}
\label{timedepdetpA1}
\begin{cases}
\frac{\partial}{\partial t}u(x,t) = \kappa\Delta u(x,t) + q(x,t)u(x,t),\\
u(x,0) = u_{0}(x),
\end{cases}
\qquad x\in\Z^d,\,t\geq 0, 
\end{equation}
has a non-negative solution if and only if the Feynman-Kac formula
\begin{equation}
\label{FKtimedepdet}
v(x,t) = E_x\left(\exp\left\{\int_0^t q(X^\kappa(s),t-s)\,ds\right\}
u_0(X^\kappa(t))\right) 
\end{equation}
is finite for all $x$ and $t$. Moreover, $v$ in {\rm \eqref{FKtimedepdet}} is the 
minimal non-negative solution to {\rm (\ref{timedepdetpA1})}.
\ep

To complement Proposition~\ref{existtimedeppA}, we need to find a condition on $\xi$ 
that leads to uniqueness of \eqref{FKtimedepdet}. This will be the first of our targets.
To answer the question of uniqueness for \emph{static} $\xi$, G\"artner and 
Molchanov~\cite{GM90} introduced the following notion.

\bd{perbel}
A field $q=\{q(x)\colon\,x \in \Z^d\}$ is said to be percolating from below if for 
every $\alpha \in \R$ the level set $\{x\in\Z^d\colon\,q(x) \leq \alpha\}$ contains 
an infinite connected component. Otherwise $q$ is said to be non-percolating from below.
\ed

\noindent
It was shown in \cite{GM90} that if $q$ is non-percolating from below, then 
\eqref{timeindeppA} has at most one non-negative solution. We will show that 
a similar condition suffices for \emph{dynamic} $\xi$, namely, \eqref{timedepdetpA1}
has at most one non-negative solution when there is a $T>0$ such that
\begin{equation}
\label{nonpercT}
q^T = \{q^T(x)\colon\,x\in\Z^d\} \quad \mbox{ with } 
\quad q^T(x) = \sup_{0 \leq t\leq T} q(x,t)
\end{equation}
is non-percolating from below (Theorem~\ref{Unique} below).  This (surprisingly 
weak) condition is fulfilled $\xi$-a.s.\ for $q=\xi$ for most choices of $\xi$. Moreover, we show that this solution
is given by the Feynman-Kac formula (Theorem~\ref{Exist} below). 

\medskip\noindent
$\bullet$ {\bf Quenched Lyapunov exponent and initial condition.}
The \emph{quenched Lyapunov exponent} associated with (\ref{pA}) with initial condition 
$u_0$ is defined as
\begin{equation}
\label{Lyapexploc1}
\lambda_0^{u_0}(\kappa) = \lim_{t\to\infty} \frac{1}{t} \log u(0,t).
\end{equation}
G\"artner, den Hollander and Maillard~\cite{GdHM11} showed that if $u_0$ has finite support, 
then the limit exists $\xi$-a.s.\ and in $L^{1}(\PP)$, is $\xi$-a.s.\ constant, and does 
not depend on $u_0$. A natural question is whether the same is true for $u_0$ bounded with 
infinite support. This question was already addressed by Drewitz, G\"artner, Ramirez and 
Sun~\cite{DGRS12}. Define 
\begin{equation}
\label{Lyapexp2}
\overline{\lambda}_0^{u_0}(\kappa) 
= \lim_{t \to \infty}\frac{1}{t}
\log E_0\left(\exp\left\{\int_{0}^{t}\xi(X^{\kappa}(s),s)\,ds\right\}
u_0(X^{\kappa}(t))\right).
\end{equation}

\bp{shapeThm}
{\rm (Drewitz, G\"artner, Ramirez and Sun \cite{DGRS12})}\\
{\rm (I)} If $\xi$ satisfies the first line of \eqref{staterg} and is bounded, then 
$\overline{\lambda}_0^{\one}(\kappa)$ exists $\xi$-a.s.\ and in $L^{1}(\P)$, and 
is $\xi$-a.s.\ constant.\\ 
{\rm (II)} If, in addition, $\xi$ is reversible in time or symmetric in space, then, 
for all $u_0$ subject to {\rm (\ref{ic})}, $\overline{\lambda}_0^{u_0}(\kappa)$ exists 
$\xi$-a.s.\ and in $L^{1}(\P)$, and coincides with $\overline{\lambda}_0^{\one}
(\kappa)$.
\ep

\noindent
The time-reversal that distinguishes $\lambda_0^{\one}(\kappa)$ from $\overline{
\lambda}_0^{\one}(\kappa)$ is non-trivial. Under appropriate space-time mixing 
conditions on $\xi$, we show how Proposition~\ref{shapeThm} can be used to settle
the existence of $\lambda_0^{u_0}(\kappa)$ with the same limit for all $u_0$ 
subject to (\ref{ic}) (Theorem~\ref{nonlocLyp} below).

\medskip\noindent
$\bullet$ {\bf Finiteness of the quenched Lyapunov exponent.}
On the one hand it follows by an application of Jensen's inequality that $\lambda_0^{u_0}(\kappa) \geq \E(\xi(0,0))$ for all $\kappa$ (see Theorem 1.2(iii) in G\"artner, den Hollander and Maillard~\cite{GdHM11} for the details), while on the other hand if 
$\xi$ is bounded from above, then also $\lambda_0^{u_0}(\kappa)<\infty$ for all $\kappa$. 
For unbounded $\xi$ the same is expected to be true under a mild assumption on the 
positive tail of $\xi$. However, settling this issue seems far from easy. The only 
two choices of $\xi$ for which finiteness has been established in the literature 
are an i.i.d.\ field of Brownian motions (Carmona and Molchanov~\cite{CM94}) 
and a Poisson random field of independent simple random walks (Kesten and 
Sidoravicius~\cite{KS03}). We will show that finiteness holds under an appropriate 
mixing condition on $\xi$ (Theorem~\ref{expgrowthres} below).

\medskip\noindent
$\bullet$ {\bf Dependence on $\kappa$.}
In \cite{GdHM11} it was shown that $\lambda_0^{\delta_0}(0)
= \E(\xi(0,0))$, $\lambda_0^{\delta_0}(\kappa)> \E(\xi(0,0))$ for $\kappa\in (0,\infty)$, and 
$\kappa\mapsto\lambda_{0}^{\delta_0}(\kappa)$ is globally Lipschitz outside any neighborhood of 
zero where it is finite. Under certain strong ``noisiness'' assumptions on $\xi$, it was further 
shown that continuity extends to zero while the Lipschitz property does not. It remained unclear,
however, which characteristics of $\xi$ are really necessary for the latter two properties to hold. 
We will show that if $\xi$ is a Markov process, then in essence a weak condition on its Dirichlet form is 
enough to ensure continuity (Theorem~\ref{continuity} and Corollary~\ref{contex} below),
whereas the non Lipschitz property holds under a weak assumption on the fluctuations of $\xi$ (Theorem~\ref{nonLipsch}). Finally, by the ergodicity of $\xi$ 
in space, it is natural to expect (see Conjecture~\ref{largekappa} below) that $\lim_{\kappa\to\infty} 
[\lambda_p^{\delta_0}(\kappa)-\lambda_0^{\delta_0}(\kappa)] = 0$ for all $p\in\N$, where 
\begin{equation}
\label{annLya}
\lambda_p^{\delta_0}(\kappa) = \lim_{t\to\infty} \frac{1}{pt} \log \E([u(0,t)]^p) 
\end{equation}
is the $p$-th \emph{annealed Lyapunov exponent} (provided this exists). It was proved for three special choices of $\xi \colon$~(1) independent 
simple random walks; (2) the symmetric exclusion process; (3) the symmetric voter model,
(for references, see \cite{GdHM11}),
that, when $d$ is large enough,  $\lim_{\kappa\to\infty} \lambda^{\delta_0}_p(\kappa) =\E(\xi(0,0))$, $p\in\N_0$. It is known from Carmona and Molchanov~\cite{CM94} that $\lim_{\kappa\to\infty} 
\lambda^{\delta_0}_p(\kappa) =\frac{1}{2}\neq\E(\xi(0,0))$ for all $p\in\N$ when $\xi$ is an i.i.d.\ field of Brownian 
motions.
\begin{remark}
\label{annealedvsquencheed}
We expect that one can define the $p$-th annealed Lyapunov exponent even for non-integer values of $p$ and that in this case $\lim_{p\to 0}\lambda_p^{\delta_0}(\kappa)=\lambda_0^{\delta_0}(\kappa)$. 
This was indeed established by Cranston, Mountford and Shiga \cite{CMS02} when $\xi$ is an i.i.d. field of Brownian motions.

\end{remark}


\subsection{Main results}
\label{S1.3}

This section contains five definitions of space-time mixing assumptions on $\xi$, six 
theorems subject to these assumptions, as well as examples of $\xi$ for which these 
assumptions are satisfied. The material is organized as Sections~\ref{S1.3.1}--\ref{S1.3.4}. 
The first theorem refers to the deterministic PAM, the other four theorems to the random 
PAM. Recall that the initial condition $u_0$ is assumed to be non-negative and bounded. 
Further recall that $\xi$ satisfies \eqref{staterg}.


\subsubsection{Definitions: Space-time blocks, G\"artner-mixing, G\"artner-regularity and G\"artner-volatility}
\label{S1.3.1}

$\bullet$ {\bf Good and bad space-time blocks.}
For $A \geq 1$, $R\in\N$, $x\in\Z^d$ and $k,b,c\in\N_0$, define the space-time blocks
\begin{equation}
\label{Bblocks}
\tilde{B}_R^{A}(x,k;b,c) = \left(\prod_{j=1}^{d}\big[(x(j)-1-b)A^R,(x(j)+1+b)A^R\big)
\cap\Z^d\right)\times[(k-c)A^R,(k+1)A^R),
\end{equation}
abbreviate $B_R^{A}(x,k)=\tilde{B}_R^{A}(x,k;0,0)$, and define the space-blocks
\begin{equation}
\label{Qbox}
Q_R^{A}(x) = x + [0,A^R)^d\cap\Z^d.
\end{equation}

\begin{figure}[htbp]
\begin{center}
\setlength{\unitlength}{0.25cm}
\begin{picture}(15,10)(0,2.5)
\put(-10,-3){\vector(1,0){23}}
\put(-10,-3){\vector(0,1){13}}
\put(-3,-1){\framebox(7,8)[]{}}
\put(-2,0){\line(0,1){4}}
\put(13.5,-4){time}
\put(-14.5,9){space}
\end{picture}
\vspace{2cm}
\caption{\small The box represents $B_{R}^{A}(x,k)$. The line is a possible 
realization of $Q_{R}^{A}(x)$.}
\label{pictFiniteness1}
\end{center}
\end{figure}
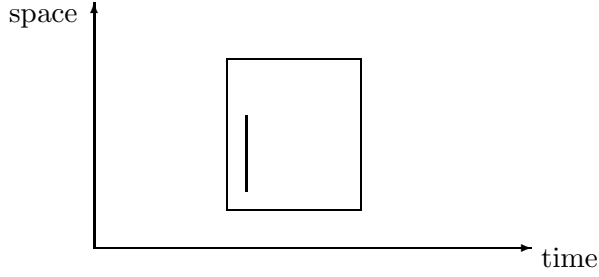

\noindent
It is convenient to extend the $\xi$-process to negative times, to obtain a two-sided
process $\overline{\xi}=(\overline{\xi}_t)_{t\in\R}$. Abbreviate $M=\mathrm{\ess}\,[\xi(0,0)]$.

\begin{figure}[htbp]
\vspace{1cm}
\begin{center}
\setlength{\unitlength}{0.25cm}
\begin{picture}(28,23)(0,2.5)
\put(-10,-3){\vector(1,0){43}}
\put(-10,-3){\vector(0,1){33}}
\put(-1.5,21.5){\dashbox{0.5}(3,4)[]{}}
\put(22.5,21.5){\dashbox{0.5}(3,4)[]{}}
\put(-1.5,1.5){\dashbox{0.5}(3,4)[]{}}
\put(22.5,1.5){\dashbox{0.5}(3,4)[]{}}
\put(-4.5,20.5){\dashbox{0.5}(6,7)[]{}}
\put(19.5,20.5){\dashbox{0,5}(6,7)[]{}}
\put(-4.5,0.5){\dashbox{0.5}(6,7)[]{}}
\put(19.5,0.5){\dashbox{0.5}(6,7)[]{}}
\put(-2,21){\framebox(6,6)[]{}}
\put(22,21){\framebox(6,6)[]{}}
\put(-2,1){\framebox(6,6)[]{}}
\put(22,1){\framebox(6,6)[]{}}
\put(-6,18){\framebox(10,12)[]{}}
\put(18,18){\framebox(10,12)[]{}}
\put(-6,-2){\framebox(10,12)[]{}}
\put(18,-2){\framebox(10,12)[]{}}
\put(-0.9,-3.5){\vector(1,0){6}}
\put(-0.9,-3,5){\vector(-1,0){6}}
\put(-4.5,-5.5){$(c+1)A^{R+1}$}
\put(-6.5,-2.5){$\circledast_1$}
\put(-6.5,9.5){$\circledast_2$}
\put(3.5,9.5){$\circledast_3$}
\put(3.5,-2.5){$\circledast_4$}
\put(-5,0){$\circledast_5$}
\put(-2.5,0.5){$\circledast_6$}
\put(34,-4){time}
\put(-14.5,29){space}
\put(11,24){\vector(1,0){6}}
\put(11,24){\vector(-1,0){6}}
\put(6,25){$(a_2-c-1)A^{R+1}$}
\put(-1,14){\vector(0,1){3}}
\put(-1,14){\vector(0,-1){3}}
\put(0,14){$(a_1-2b-2)A^{R+1}$}
\end{picture}
\vspace{2cm}
\caption{\small The dashed blocks are $R$-blocks, i.e., $B_R^{A}(x,k)$ (inner) and 
$\tilde{B}_R^{A}(x,k;b,c)$ (outer) for some choice of $A,x,k,b,c$. The solid 
blocks are $(R+1)$-blocks, i.e., $B_{R+1}^{A}(y,l)$ (inner) and $\tilde{B}_{R+1}^{A}(y,l;b,c)$ (outer) for a choice of $A,y,l,b,c$ such that they contain the corresponding $R$-blocks.
Furthermore, $\{\circledast_i\}_{i=1,2,3,4,5,6}$ represents the space-time coordinates $\circledast_1 = ((y-1-b)A^{R+1},(l-c)A^{R+1})$,$\circledast_2=((y+1+b)A^{R+1},(l-c)A^{R+1})$, $\circledast_3= ((y+1+b)A^{R+1},(l+1)A^{R+1})$, $\circledast_4= ((y-1-b)A^{R+1}, (l+1)A^{R+1})$, $\circledast_5=((x-1-b)A^{R}, (k-c)A^R)$ and
$\circledast_6= ((y-1)A^{R+1},lA^{R+1})$.}
\label{pictFiniteness2}
\end{center}
\end{figure}
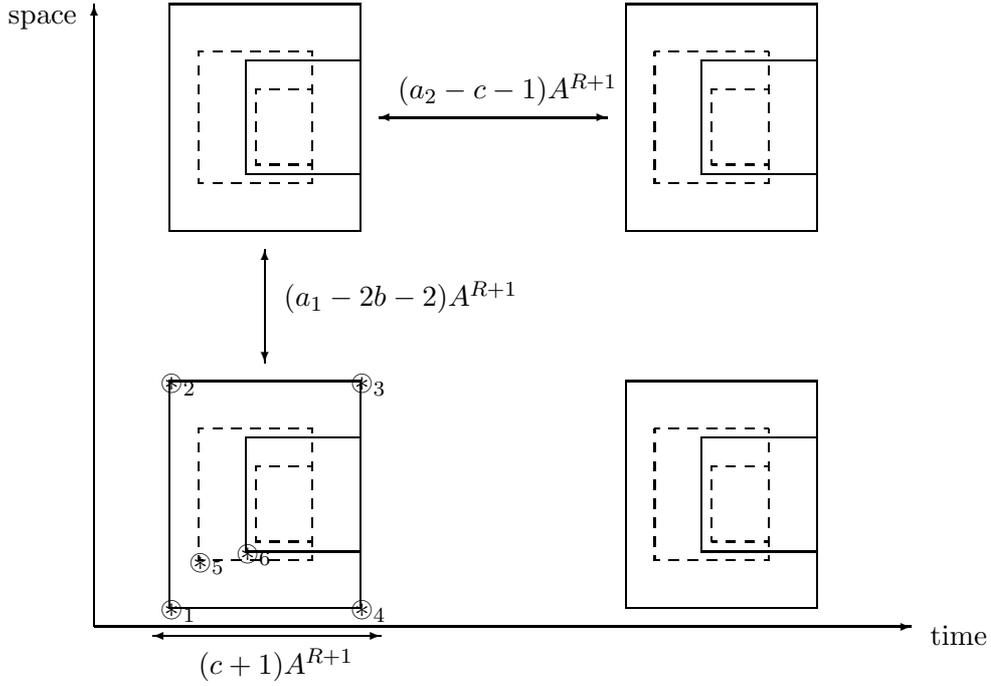

\bd{Cgood}
For $A \geq 1$, $R \in \N$, $x \in \Z^d$, $k\in\N$, $C\in [0,M]$ and $b,c\in\N_0$, the 
block $B_R^{A}(x,k)$ is called $(C,b,c)$-good when
\begin{align}
\label{BRgood}
\sum_{z \in Q_R^{A}(y)}\xi(z,s) \leq CA^{Rd}
\qquad \forall\,y\in\Z^d,\,s\geq 0\colon\,
Q_R^{A}(y)\times\{s\} \subseteq \tilde{B}_R^{A}(x,k;b,c). 
\end{align}
Otherwise it is called $(C,b,c)$-bad. 
\ed

\medskip\noindent
$\bullet$ {\bf G\"artner-mixing.}
For $A\geq 1$, $R\in\N$, $x\in\Z^d$, $k\in\N$, $C\in [0,M]$ and $b,c\in\N_0$, let 
\begin{equation}
\label{Aevent}
\begin{aligned}
&\mathcal{A}_R^{A,C}(x,k;b,c)\\ 
&= \big\{\mbox{$B_{R+1}^{A}(x,k)$ is $(C,b,c)$-good, 
but contains an $R$-block that is $(C,b,c)$-bad}\big\}. 
\end{aligned}
\end{equation}
In terms of these events we define the following \emph{space-time mixing conditions}
(see Fig.~\ref{pictFiniteness2}).
For $D \subset \Z^d\times\R$,
let $\sigma(D)$ be the $\sigma$-field generated by $\{\overline{\xi}(x,t)\colon\,(x,t)
\in D\}$.

\bd{Gartnermix}{\rm {\bf [G\"artner-mixing]}}\\
For $a_1,a_2 \in \N$, denote by $\Delta_n(a_1,a_2)$ the set of $\Z^d\times\N$-valued 
sequences $\{(x_i,k_i)\}_{i=0}^n$ that are increasing with respect to the lexicographic 
ordering of $\Z^d\times\N$ and are such that for all $0\leq i<j\leq n$ 
\begin{equation}
\label{mod} 
x_j\equiv x_i\mod a_1 \quad \mbox{ and } \quad k_j\equiv k_i\mod a_2.
\end{equation}
{\rm (a)} 
$\xi$ is called $(A,C,b,c)$-type-I G\"artner-mixing when there are $a_1, a_2 \in \N$ and a constant $K>0$ such 
that there is an $R_0 \in\N$ such that, for all $R\in\N$ with $R\geq R_0$ and all $n\in \N$,
\begin{equation}
\label{wGartner}
\begin{aligned}
\sup_{(x_i,k_i)_{i=0}^n\in\Delta_n(a_1,a_2)}
\P\bigg(\bigcap_{i=0}^{n}\mathcal{A}_R^{A,C}(x_i,k_i;b,c) \bigg)
\leq K\bigg (A^{(1+2d)}\bigg)^{-R(1+d)n}.
\end{aligned}
\end{equation} 
{\rm (b)} 
$\xi$ is called $(A,C,b,c)$-type-II G\"artner-mixing when for each family of events 
\begin{equation}
\mathcal{A}_{i}^{R} \in \sigma(B_{R+1}^{A}(x_i,k_i)),
\qquad (x_i,k_i)_{i=0}^{n} \in  \Delta_n(a_1,a_2),
\end{equation}
that are invariant under space-time shifts and satisfy
\begin{equation}
\label{Aeventtozero}
\lim_{R\to \infty}\P(\mathcal{A}_{i}^{R})=0,
\end{equation}
there are $a_1, a_2\in\N$ and a constant $K>0$ such that for each $\delta>0$ there is an $R_0\in\N$ such that, 
for all $n\in\N$,
\begin{equation}
\label{sGartner}
\begin{aligned}
\P\bigg(\bigcap_{i=0}^{n}
\left\{B_{R+1}^{A}(x_i,k_i) \mbox{ is $(C,b,c)$-good}, \mathcal{A}_{i}^{R}\right\}\bigg)
\leq K\delta^{n}\qquad R\geq R_0, R\in\N.
\end{aligned}
\end{equation}
{\rm (c)} 
$\xi$ is called type-I, respectively type-II, G\"artner-mixing, when there are $A\geq 1$, 
$C \in [0,M]$, $R\in\N$, $b,c\in \N$ such that $\xi$ is $(A,C,b,c)$-type-I, respectively, 
$(A,C,b,c)$-type-II, G\"artner-mixing.
\ed

\bd{Gartnerhypmix}{\rm {\bf [G\"artner-hyper-mixing]}}\\
{\rm (a)} 
$\xi$ is called G\"artner-positive-hyper-mixing when\\
{\rm (a1)} 
$\E[e^{q\sup_{s\in[0,1]} \xi(0,s)}] < \infty$ 
for all $q\geq 0$.\\
{\rm (a2)}
There are $b,c\in\N$ and a constant $C$ such that for each $A_0>1$ 
one can find $A\geq A_0$ such that $\xi$ is $(A,C,b,c)$-type-I 
G\"artner-mixing.\\
{\rm (a3)}
There are $R_0,C_0\geq 1$ such that 
\begin{equation}
\label{condexpgrowthB}
\P\left(\sup_{s \in [0,1]} \frac{1}{|B_{R}|}\sum_{y \in B_{R}} \xi(y,s) 
\geq C \right) \leq |B_{R}|^{-\alpha} \quad \forall\, R\geq R_0, C\geq C_0,
\end{equation}
for some $\alpha > (1+2d)(2+d)/d$, where $B_{R} = [-R,R]^d\cap\Z^d$.\\
{\rm (b)}
$\xi$ is called G\"artner-negative-hyper-mixing, when $-\xi$ is 
G\"artner-positive-hyper-mixing.
\ed

\br{remintgartmix}
{\rm If $\xi$ is bounded from above, then $\xi$ is G\"artner-positive-hyper-mixing. For those examples where $\xi(x,t)$ represents ``the number of particles at site $x$ 
at time $t$'', we may view G\"artner-mixing as a consequence of the fact that there 
are not enough particles in the blocks $\tilde{B}_R^{A}(x_i,k_i;b,c)$ that manage 
to travel to the blocks $\tilde{B}_R^{A}(x_j,k_j;b,c)$. Indeed, if there is a bad 
block on scale $R$ that is contained in a good block on scale $R+1$, then in some 
neighborhood of this bad block the particle density cannot be too large. This also 
explains why we must work with the extended blocks $\tilde{B}_R^{A}(x,k;b,c)$ 
instead of with the original blocks $B_{R}^{A}(x,k;0,0)$. Indeed, the surroundings 
of a bad block on scale $R$ can be bad when it is located near the boundary of a good 
block on scale $R+1$ (see Fig.~\ref{pictFiniteness2}).}
\er

\medskip\noindent
$\bullet$ {\bf G\"artner-regularity and G\"artner-volatility.}
Recall that $||\cdot||$ denotes the lattice-norm, see the line following (\ref{dL}). We say that 
$\Phi\colon\,[0,t]\to\Z^d$ is a path when
\begin{equation}
\label{path}
\|\Phi(s)-\Phi(s-)\| \leq 1 \qquad \forall\,s \in [0,t].
\end{equation}
We write $\Phi \in B_R$  when 
$\|\Phi(s)\| \leq R$ for all $s \in [0,t]$ and denote by $N(\Phi,t)$ the number of jumps of $\Phi$ up to time $t$.

\bd{Gartnerposreg}{\rm {\bf [G\"artner-regularity]}}\\
$\xi$ is called G\"artner-regular when\\
{\rm (a)} 
$\xi$ is G\"artner-negative-hyper-mixing and G\"artner-positive-hyper-mixing.\\
{\rm (b)} 
There are $t_0>0$ and $n_0\in\N$ such that for every $\delta_1 >0$ there is a 
$\delta_2 = \delta_2(\delta_1)>0$ such that
\begin{equation}
\label{contcond1}
\begin{aligned}
&\PP\left(\sum_{j=1}^n \int_{(j-1)t+1}^{jt} \xi\big(\Phi((j-1)t +1),s\big)\,ds 
\geq \delta_1 nt\right) \leq e^{-\delta_2 nt}\\
&\qquad\qquad\qquad\forall\,t \geq t_0,\,n \geq n_0,\,\Phi \in B_{tn}.
\end{aligned}
\end{equation}
\ed

\bd{Gartnernegreg}{\rm {\bf [G\"artner-volatility]}}\\
$\xi$ is called G\"artner-volatile when\\
{\rm (a)} $\xi$ is G\"artner-negative-hyper-mixing.\\
{\rm (b)}
\begin{equation}
\label{nonLipschcond1}
\lim_{t \to \infty} \frac{1}{\log t} 
\EE\left(\Big|\int_0^t [\xi(0,s)-\xi(e,s)]\,ds\Big| \right) 
= \infty \qquad \mbox{for some } e \in \Z^d \mbox{ with } \|e\|=1, 
\end{equation}
\ed

\br{remcondC2}
{\rm Corollary~\ref{contex} below will show that condition {\rm (b)} in Definition~\ref{Gartnerposreg} 
is satisfied as soon as the Dirichlet form of $\xi$ is non-degenerate, i.e., has a unique zero (see 
Section \ref{S7}).}
\er


\subsubsection{Theorems: Uniqueness, existence, finiteness and initial condition}
\label{S1.3.2}
Recall the definition of $q^{T}$ (see (\ref{nonpercT})), 
the condition on $u_0$ in (\ref{ic}) and the condition on $\xi$ in (\ref{staterg}).
\bt{Unique}{\rm {\bf [Uniqueness]}}
Consider a deterministic $q\colon\,\Z^d \times [0,\infty) \to \R$ such that:\\
{\rm (1)} There is a $T>0$ such that $q^T$ is non-percolating from below.\\
{\rm (2)} $q^T(x) < \infty$ for all $T>0$ and $x \in \Z^d$.\\
Then the Cauchy problem
\begin{equation}
\label{timedepdetpA2}
\begin{cases}
\frac{\partial}{\partial t}u(x,t) = \kappa\Delta u(x,t) + q(x,t)u(x,t),\\
u(x,0) = u_0(x),
\end{cases}
\qquad x\in\Z^d,\,t\geq 0, 
\end{equation}
has at most one non-negative solution.
\et

\bt{Exist}{\rm {\bf [Existence]}}
Suppose that:\\
{\rm (1)} $s \mapsto \xi(x,s)$ is locally integrable for every $x$, $\xi$-a.s.\\
{\rm (2)} $\EE(e^{q\xi(0,0)}) < \infty$ for all $q \geq 0$.\\
Then the function defined by the Feynman-Kac formula
\begin{equation}
\label{FK}
u(x,t) = E_x \left(\exp\left\{\int_0^t \xi(X^\kappa(s), t-s)\, 
ds\right\} u_0(X^\kappa(t))\right)
\end{equation}
solves {\rm (\ref{pA})} with initial condition $u_0$.
\et
From now on we assume that $\xi$ satisfies the conditions of Theorems~\ref{Unique}--\ref{Exist} (where $q$ is replaced by $\xi$ in Theorem~\ref{Unique}).
\bt{expgrowthres}{\rm {\bf [Finiteness]}} 
If $\xi$ is G\"artner-positive-hyper-mixing, then $\lambda_0^{\delta_0}(\kappa)<\infty$.
\et

\noindent
From now on we also assume that $\xi$ satisfies the conditions of Theorem~\ref{expgrowthres}.
The following result extends G\"artner, den Hollander and Maillard~\cite{GdHM11}, Theorem 1.1, 
in which it was shown that for the initial condition $u_0=\delta_0$ the quenched Lyapunov 
exponent exists and is constant $\xi$-a.s.

\bt{nonlocLyp}{\rm {\bf [Initial Condition]}}
If $\xi$ is reversible in time or symmetric in space, type-II
G\"artner-mixing and G\"artner-negative-hyper-mixing, then $\lambda_{0}^{u_0}(\kappa)
= \lim_{t\to\infty} \frac{1}{t}\log u(0,t)$ exists $\xi$-a.s.\ and in $L^{1}(\P)$, is 
constant $\xi$-a.s., and is independent of $u_0$.
\et


\subsubsection{Theorems: Dependence on $\kappa$}
\label{S1.3.3}

\bt{continuity}{\rm {\bf [Continuity at $\kappa=0$]}}
If $\xi$ is G\"artner-regular, 
then $\kappa \mapsto \lambda_0^{\delta_0}(\kappa)$ is  continuous at zero.
\et

\bt{nonLipsch}{\rm {\bf [Not Lipschitz at $\kappa=0$]}}
If $\xi$ is G\"artner-volatile,
then $\kappa \mapsto \lambda_0^{\delta_0}(\kappa)$ is not Lipschitz continuous in zero.
\et

\br{nonLipschrm}
{\rm Theorem~\ref{nonLipsch} was already shown in \cite{GdHM11},
under the additional assumption that $\xi$ is bounded from below.}
\er

\subsubsection{Examples}
\label{S1.3.4}

We state two corollaries in which we give examples of classes of $\xi$ for which
the conditions in Theorems~\ref{expgrowthres}--\ref{continuity} are satisfied.

\bc{mixingex}{\rm {\bf [Examples for Theorems \ref{expgrowthres}--\ref{nonlocLyp}]}}\\
{\rm (1)} 
Let $X=(X_t)_{t\geq 0}$ be a stationary and ergodic $\R$-valued Markov process.
Let $(X_{\cdot}(x))_{x\in\Z^d}$ be independent copies of $X$. Define $\xi$ by 
$\xi(x,t)=X_t(x)$. If
\begin{equation}
\label{Marksup}
\E\left[e^{q\sup_{s\in[0,1]} X_s}\right] < \infty \qquad \forall\, q\geq 0,
\end{equation}
then $\xi$ fulfills the conditions of Theorem~{\rm \ref{expgrowthres}}. If, moreover, 
the left-hand side of (\ref{Marksup}) is finite for all $q\leq 0$, then $\xi$ satisfies 
the conditions of Theorem~{\rm \ref{nonlocLyp}}.\\
{\rm (2)} 
Let $\xi$ be the zero-range process with rate function $g\colon\N_0\to (0,\infty)$, 
$g(k) = k^{\beta}$, $\beta \in (0,1]$, and transition probabilities given by a simple 
random walk on $\Z^d$. If $\xi$ starts from the product measure $\pi_{\rho}$, 
$\rho\in (0,\infty)$, with marginals
\begin{equation}
\label{ZRmeasure}
\pi_{\rho}\big\{\eta \in \N_0^{\Z^d}\colon\, \eta(x)=k\big\} =
\begin{cases}
\gamma\,\frac{\rho^{k}}{g(1) \times\cdots\times g(k)}, &\mbox{ if $k>0$}, \\
\gamma,  &\mbox{ if $k=0$},
\end{cases}
\end{equation}
where $\gamma \in (0,\infty)$ is a normalization constant, then $\xi$ satisfies the 
conditions of Theorems~{\rm \ref{expgrowthres}--\ref{nonlocLyp}}.
\ec

\bc{contex}{\rm {\bf [Examples for Theorem \ref{continuity}]}}
{\rm (1)} 
If $\xi$ is a bounded interacting particle system in the so-called $M<\varepsilon$ regime 
(see Liggett~{\rm \cite{L85}}), then the conditions of Theorem~{\rm \ref{continuity}} 
are satisfied.\\
{\rm (2)}
If $\xi$ is the exclusion process with an irreducible, symmetric and transient random walk
transition kernel, then the conditions of Theorem~{\rm \ref{continuity}} are satisfied.\\
{\rm (3)} 
If $\xi$ is the dynamics defined by
\begin{equation}
\label{ISRW}
\xi(x,t) = \sum_{y \in\Z^d} \sum_{j=1}^{N_y} \delta_{Y_{j}^{y}(t)}(x),
\end{equation}
where $\{Y_{j}^{y}\colon\,y\in\Z^d,1\leq j \leq N_y, Y_{j}^{y}(0) = y\}$ is a collection of independent 
continuous-time simple random walks jumping at rate one,
and $(N_y)_{y\in\Z^d}$ is a Poisson random field with intensity $\nu$ for some $\nu \in 
(0,\infty)$. If $d\geq 3$, then the conditions of Theorem~{\rm \ref{continuity}} are satisfied.
\ec

Corollaries~\ref{mixingex}--\ref{contex} list only a few examples that match the conditions. 
It is a separate problem to verify these conditions for as broad a class of interacting particle 
systems as possible.


\subsection{Discussion and a conjecture}
\label{S1.4}

\begin{figure}[htbp]
\vspace{1cm}
\begin{center}
\setlength{\unitlength}{0.25cm}
\begin{picture}(20,12)(17,2)
\put(0,2){\line(22,0){22}}
\put(0,2){\line(0,11){11}}
{\thicklines
\qbezier(0,2)(.4,8.8)(3,8)
\qbezier(3,8)(3.8,7.8)(5,6.6)
\qbezier(5,6.6)(9,2.6)(21,2.3)
}
\put(-1.3,1){$0$}
\put(0,2){\circle*{.35}}
\put(23,1.8){$\kappa$}
\put(-1,14){$\lambda^{u_0}_0(\kappa)$}
\put(29,2){\line(22,0){22}}
\put(29,2){\line(0,11){11}}
{\thicklines
\qbezier(29,2)(29.2,8)(52,12)
}
\put(27.8,1.1){$0$}
\put(29,2){\circle*{.35}}
\put(52,1.8){$\kappa$}
\put(28,14){$\lambda^{u_0}_0(\kappa)$}
\end{picture}
\caption{\small Qualitative picture of $\kappa\mapsto\lambda_0^{u_0}(\kappa)$
in the weakly, respectively, strongly catalytic regime.}
\label{fig-lambda0}
\end{center}
\end{figure}
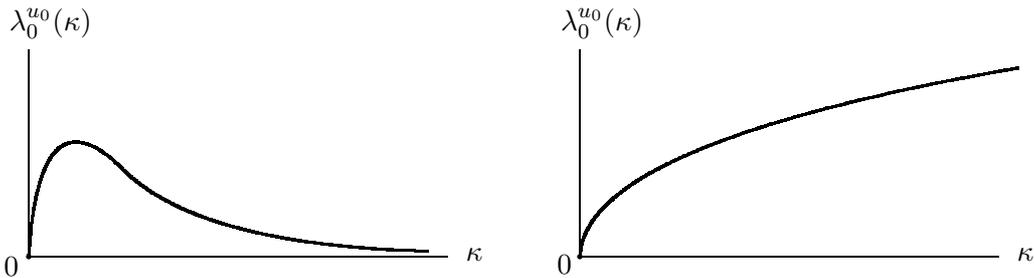

The proofs of Theorems~\ref{Unique}--\ref{nonLipsch} and Corollaries~\ref{mixingex}--\ref{contex} 
are given in Sections \ref{S2}--\ref{S7}. The content of Theorems~\ref{Unique}--\ref{nonLipsch} 
is summarized in Fig.~\ref{fig-lambda0}.

The importance of $\lambda^{u_0}_0(\kappa)$ within the \emph{population dynamics} 
interpretation of the parabolic Anderson model, as explained in Section~\ref{S1.1}, 
is the following. For $t>0$, randomly draw a $B$-particle from the population of 
$B$-particles at the origin. Let $L_t$ be the random time this $B$-particle and 
its ancestors have spent on top of $A$-particles. By appealing to the ergodic theorem, 
it may be shown that $\lim_{t\to\infty} L_t/t=\lambda^{u_0}_0(\kappa)$ a.s. Thus, 
$\lambda^{u_0}_0(\kappa)$ is the fraction of time the \emph{best} $B$-particles 
spend on top of $A$-particles, where best means that they come from the fastest growing 
family (``survival of the fittest''). Fig.~\ref{fig-lambda0} shows that for all 
$\kappa \in (0,\infty)$ \emph{clumping} occurs: the limiting fraction is strictly 
larger than the density of $A$-particles. In the limit as $\kappa\downarrow 0$ 
the clumping vanishes because the motion of the $A$-particles is ergodic in time. 
The clumping is hard to suppress for $\kappa\downarrow 0$: even a tiny bit of 
mobility allows the best $B$-particles and their ancestors to successfully
``hunt down'' the $A$-particles.

In the limit as $\kappa\to\infty$ we expect the quenched Lyapunov exponent to 
merge with the \emph{annealed Lyapunov exponents} defined in \eqref{annLya},
with $\delta_0$ replaced by $u_0$.

\begin{conjecture}
\label{largekappa}
$\lim_{\kappa\to\infty} [\lambda_p^{u_0}(\kappa)-\lambda_0^{u_0}(\kappa)]=0$ 
for all $p\in\N$.
\end{conjecture}

\noindent
The reason is that for large $\kappa$ the $B$-particles can easily find the largest 
clumps of $A$-particles and spend most of their time there, so that it does not 
matter much whether the largest clumps are close to the origin or not.

It remains to identify the scaling behaviour of $\lambda^{u_0}_0(\kappa)$ for 
$\kappa\downarrow 0$ and $\kappa\to\infty$. Under strong noisiness conditions on 
$\xi$, it was shown in G\"artner, den Hollander and Maillard~\cite{GdHM11} that 
$\lambda^{u_0}(\kappa)$ tends to zero like $1/\log(1/\kappa)$ (in a rough 
sense), while it tends to $\E(\xi(0,0))$ as $\kappa\to\infty$. For the \emph{annealed} 
Lyapunov exponents $\lambda^{u_0}_p(\kappa)$, $p\in\N$, there is no singular 
behavior as $\kappa\downarrow 0$, in particular, they are Lipschitz continuous at
$\kappa=0$ with $\lambda^{u_0}_p(0)>\E(\xi(0,0))$. For three specific choices 
of $\xi$ it was shown that $\lambda^{u_0}_p(\kappa)$ with $u_0\equiv 1$ decays like 
$1/\kappa$ as $\kappa\to\infty$ (see \cite{GdHM11} and references therein). A 
distinction is needed between the \emph{strongly catalytic regime} for which 
$\lambda^{u_0}_p(\kappa) = \infty$ for all $\kappa \in [0,\infty)$, and the 
\emph{weakly catalytic regime} for which $\lambda^{u_0}_p(\kappa) < \infty$ for 
all $\kappa\in [0,\infty)$. (These regimes were introduced by G\"artner and 
den Hollander~\cite{GdH06} for independent simple random walks.) We expect 
Conjecture~\ref{largekappa} to be valid in both regimes. 


\section{Existence and uniqueness of the solution}
\label{S2}

In this section we prove Theorem~\ref{Unique} (uniqueness; Section~\ref{S2.1}) 
and Theorem~\ref{Exist} (existence; Section~\ref{S2.2}).


\subsection{Uniqueness}
\label{S2.1}

The proof of Theorem \ref{Unique} is based on the following lemma.

\bl{monotonicity} 
Let $q_i\colon \Z^d \times [0,\infty)\to\R$, $i\in\{1,2\}$, satisfy conditions 
{\rm (1)--(2)} in Theorem~{\rm \ref{Unique}} and be such that, for a given initial 
condition $u_0$, the two corresponding Cauchy problems
\begin{equation}
\label{Cauchy12}
\begin{cases}
\frac{\partial}{\partial t}u_i(x,t) = \kappa\Delta u_i(x,t) + q_i(x,t)u_i(x,t),\\
u_i(x,0) = u_0(x),
\end{cases}
\qquad x\in\Z^d,\,t\geq 0,\,i\in\{1,2\},
\end{equation}
have a solution. If there exists a $T>0$ such that $q_1(x,t) \geq q_2(x,t)$ for all 
$x \in \Z^d$ and $t \in [0,T]$, then $u_1(x,t) \geq u_2(x,t)$ for all $x \in \Z^d$ 
and $t\in [0,T]$, where $u_1$ and $u_2$ are any two solutions of \eqref{Cauchy12}.
\el

We first prove Theorem~\ref{Unique} subject to Lemma~\ref{monotonicity}.

\bpr
Note from Definition~\ref{perbel} that whenever $q^T$ is non-percolating from below 
for $T=T_0$ for some $T_0>0$, then the same is true for all $T \geq T_0$. Fix $T 
\geq T_0$, and let $u$ be a non-negative solution of (\ref{timedepdetpA2}) with zero 
initial condition, i.e., $u_0(x) = 0$ for all $x \in \Z^d$. It is sufficient to prove 
that $u(x,t) = 0$ for all $x\in\Z^d$ and $t\in [0,T]$. 

Let $v$ be the solution of the Cauchy problem
\begin{equation}
\label{qTCauchy}
\begin{cases}
\frac{\partial}{\partial t}v(x,t) = \kappa \Delta v(x,t) + q^T(x) v(x,t),\\
v(x,0) = v_0(x) = 0,
\end{cases}
\qquad x\in\Z^d, t\in[0,T],
\end{equation}
which exists because the corresponding Feynman-Kac representation is zero by G\"artner 
and Molchanov~\cite{GM90}, Lemma 2.2. By Lemma \ref{monotonicity} it follows that 
$0\leq u\leq v$ on $\Z^d \times [0,T]$. Using that $q^T$ is non-percolating from below, 
we may apply \cite{GM90}, Lemma 2.3, to conclude that (\ref{qTCauchy}) has at most one 
solution. Hence $u = v = 0$ on $\Z^d \times [0,T]$, which gives the claim.
\epr

We next prove Lemma~\ref{monotonicity}. 

\bpr
Fix $R\in\N$. Let $B_R = [-R,R]^d \cap \Z^d$, $\mathrm{int}(B_R) = (-R,R)^d \cap \Z^d$, and 
$\partial B_R = B_R \backslash \mathrm{int}(B_R)$. If $u_1$ and $u_2$ are solutions of 
(\ref{Cauchy12}) on $\Z^d \times [0,\infty)$, then they are also solutions on 
$B_R \times [0,T]$. More precisely, for $i\in\{1,2\}$, $u_i$ is a solution of 
the Cauchy problem
\begin{equation}
\label{Cauchyblock}
\left\{\begin{array}{ll}
\frac{\partial}{\partial t}v(x,t) = \kappa \Delta v(x,t) +  q_i(x,t) v(x,t), 
&(x,t) \in \mathrm{int}(B_R) \times [0,T],\\
v(x,0) = u_0(x),
&x \in B_R,\\
v(x,t) = u_i(x,t),
&(x,t) \in \partial B_R \times [0,T].
\end{array}
\right.
\end{equation}
Recall that $q_1 \geq q_2$ on $\Z^d \times [0,T]$. Choose $c_R^T$ such that
\begin{equation}
\label{cNT}
c_R^T > \max_{x \in B_R,\,t \in [0,T]} q_1(x,t)
\geq \max_{x \in B_R,\,t \in [0,T]} q_2(x,t),
\end{equation}
and abbreviate
\begin{equation}
\label{vdef}
\left\{
\begin{array}{ll}
v(x,t) = e^{-c_R^Tt}\,\big[u_1(x,t)-u_2(x,t)\big], 
&(x,t) \in B_R \times [0,T],\\
\bar{Q}_i = q_i - c_R^T, 
&i \in \{1,2\}.
\end{array}
\right.
\end{equation} 
Then, by \eqref{Cauchyblock}, $v$ satisfies
\begin{equation}
\label{transformedeq}
\left\{
\begin{array}{lll}
\frac{\partial}{\partial t} v(x,t) = \kappa \Delta v(x,t) 
&\\ 
\qquad\qquad\qquad + e^{-c_R^Tt}\,\bar{Q}_1(x,t)u_1(x,t) 
- e^{-c_R^Tt}\,\bar{Q}_2(x,t)u_2(x,t), 
&(x,t) \in \mathrm{int}(B_R) \times [0,T], \\
v(x,0) = 0,
&x \in B_R, \\
v(x,t) = e^{-c_R^Tt}\,\big[u_1(x,t)- u_2(x,t)\big], 
&(x,t) \in \partial B_R \times [0,T].
\end{array}
\right.
\end{equation}
Now, suppose that there exists a $(x_\ast,t_\ast) \in \mathrm{int}(B_R) \times [0,T]$ such that 
\begin{equation}
\label{minvdef}
v(x_{\ast}, t_{\ast}) = \min_{x \in \mathrm{int}(B_R),\,t \in [0,T]} v(x,t) < 0.
\end{equation}
Then
\begin{equation}
\label{contra1}
\frac{\partial}{\partial t}v(x_\ast,t_\ast) \leq 0
\end{equation}
and
\begin{equation}
\label{contra2}
\Delta v(x_\ast,t_\ast) = \sum_{ {y \in \Z^d} \atop {\|y-x_\ast\|=1} } 
\big[v(y,t_{\ast}) - v(x_{\ast},t_{\ast})\big] \geq 0.
\end{equation}
Moreover, by (\ref{cNT}--\ref{vdef}) and (\ref{minvdef}),
\begin{equation}
\label{contra3}
\begin{aligned}
&e^{-c_R^Tt_\ast}\,\bar{Q}_1(x_\ast,t_\ast)u_1(x_\ast,t_\ast) 
- e^{-c_R^Tt_\ast}\,\bar{Q}_2(x_\ast,t_\ast)u_2(x_\ast,t_\ast)\\
&\qquad = \big[q_1(x_\ast,t_\ast) - c_R^T\big]\,v(x_\ast,t_\ast)
+ [q_1(x_\ast,t_\ast) - q_2(x_\ast,t_\ast)]\,e^{-c_R^Tt_\ast}\,u_2(x_\ast,t_\ast)
> 0.
\end{aligned}
\end{equation}
But (\ref{contra1}--\ref{contra3}) contradict the first line of (\ref{transformedeq})
at $(x,t)=(x_\ast,t_\ast)$. Hence \eqref{minvdef} fails, and so it follows from 
\eqref{vdef} that $u_1(x,t) \geq u_2(x,t)$ for all $x \in \mathrm{int}(B_R)$ and $t \in [0,T]$. 
Since $R$ can be chosen arbitrarily, the claim follows.
\epr


\subsection{Existence}
\label{S2.2}

In the sequel we use the abbreviations
\begin{eqnarray}
\label{timerevIntegral}
\cI^{\kappa}(a,b,c) &= \int_a^b \xi(X^\kappa(s),c-s)ds, \qquad 0 \leq a \leq b \leq c,\\
\label{timeIntegral}
\overline{\cI}^{\kappa}(a,b,c) &= \int_a^b \xi(X^\kappa(s),c+s)ds, 
\qquad 0 \leq a \leq b \leq c.
\end{eqnarray}

\bpr
To prove Theorem~\ref{Exist}, by Proposition~\ref{existtimedeppA} it is enough to show 
that 
\begin{equation}
\label{FK1}
E_x\left(e^{\cI^{\kappa}(0,t,t)}
u_0(X^{\kappa}(t))\right)<\infty
\qquad \forall x\in\Z^d,\,t\geq 0.
\end{equation}
Since $u_0$ is assumed to be non-negative and bounded (recall (\ref{ic})), without 
loss of generality we may take $u_0 \equiv 1$. We give the proof for $x=0$, the 
extension to $x \in \Z^d$ being straightforward.
\noindent
Fix $q \in \Q \cap [0,\infty)$. Using Jensen's inequality and the stationarity of 
$\xi$, we have (recall (\ref{staterg}))
\begin{equation}
\begin{aligned}
\EE\left(E_0\left(e^{\cI^{\kappa}(0,q,q)}\right)\right) 
&=E_0\left(\EE\left(e^{\cI^{\kappa}(0,q,q)}\right)\right)\\
&\leq E_0\left(\EE\left(\frac{1}{q} \int_0^q 
\exp\Big\{q \xi(X^\kappa(s),q-s)\Big\}\,ds\right)\right)\\ 
&= E_0\left(\frac{1}{q} \int_0^q \EE\Big(
\exp\Big\{q \xi(0,0)\Big\}\Big)\,ds\right)\\
&=\EE\big(e^{q \xi(0,0)}\big)<\infty,
\end{aligned}
\end{equation}
where the finiteness follows by condition {\rm (2)}. Hence, for every $q \in \Q \cap [0,\infty)$ 
there exists a set $A_q$ with $\P(A_q)=1$ such that 
\begin{equation}
\label{+FK}
E_0\left(e^{\cI^{\kappa}(0,q,q)}\right)<\infty
\qquad \forall\,\xi \in A_q.
\end{equation}
To extend (\ref{+FK}) to $t\in [0,\infty)$, note that, by the Markov property of $X^\kappa$ 
applied at time $q-t$, $q>t$, we have
\begin{equation}
\label{endexist}
\begin{aligned}
E_0\left(e^{\cI^{\kappa}(0,q,q)}\right)
&\geq E_0\left(e^{\cI^{\kappa}(0,q,q)}
\one\Big\{X^\kappa(r)=0\,\,\forall\,r \in [0,q-t]\Big\}\right)\\
&= e^{\int_{0}^{q-t}\xi(0,q-s)\,ds}
P_0\Big(X^\kappa(r)=0\,\,\forall\,r \in [0,q-t]\Big) 
E_0\left(e^{\cI^{\kappa}(0,t,t)} \right).
\end{aligned}
\end{equation}
Because $s \mapsto \xi(0,s)$ is locally integrable $\xi$-a.s.\ by condition {\rm (1)}, we 
have $\int_{0}^{q-t}\xi(0,q-s)\,ds>- \infty$ $\xi$-a.s. The claim now follows from 
(\ref{+FK}--\ref{endexist}) by picking $q \in  \Q \cap [0,\infty)$ and $t \in [0,\infty)$.
\epr


\section{Finiteness of the quenched Lyapunov exponent}
\label{S3}

In this section we prove Theorem~\ref{expgrowthres}. In Section~\ref{S3.1} we sketch the 
strategy of the proof. In Sections~\ref{S3.2}--\ref{S3.6} the details are worked out.


\subsection{Strategy of the proof}
\label{S3.1}

The proof uses ideas from Kesten and Sidoravicius~\cite{KS03}. Fix $C,b,c$ according to 
our assumptions on $\xi$. For $j\in\N$ and $t>0$, define the set of random walk paths
\begin{align}
\label{PiL}
\Pi(j,t)= \Bigg\{&\Phi\colon\, [0,t] \to \Z^d\colon\, \Phi\mbox{ makes } j \mbox{ jumps, }
\Phi([0,t]) \subseteq [-C_1t\log t,C_1t\log t]^{d}\cap\Z^d\Bigg\},
\end{align}
where $C_1$ will be determined later on. Abbreviate $[C_1]_t=[-C_1t\log t,C_1t\log t]^{d}
\cap\Z^d$. For $A\geq 1$, $R\in\N$ and $\Phi \in \Pi(j,t)$, define
\begin{eqnarray}
\label{psiPhi}
\Psi_{R}^{A}(\Phi) 
&=& \mbox{number of good $(R+1)$-blocks crossed by $\Phi$ containing a bad $R$-block},\\
\label{psiL}
\Psi_{R}^{A,j} &=& \sup_{\Phi \in \Pi(j,t)} \Psi_{R}^{A}(\Phi),\\
\label{XiPhi}
\Xi_{R}^{A}(\Phi) &=& \mbox{number of bad $R$-blocks crossed by $\Phi$},\\
\label{PhiL}
\Xi_{R}^{A,j} &=& \sup_{\Phi \in \Pi(j,t)} \Xi_{R}^{A}(\Phi).
\end{eqnarray}
The proof comes in 5 steps, organized as Sections~\ref{S3.2}--\ref{S3.6}: (1) the 
Feynman-Kac formula may be restricted to paths contained in $[C_1]_t$; (2) there 
are no bad $R$-blocks for sufficiently large $R$; (3) the Feynman-Kac formula can 
be estimated in terms of bad $R$-blocks; (4) bounds can be derived on the number of 
bad $R$-blocks; (5) completion of the proof.


\subsection{Step 1: Restriction to $[C_1]_t$}
\label{S3.2}

\bl{supremumlem}
Fix $C_1>0$. Suppose that $\E\left(e^{q\sup_{s \in [0,1]} \xi(0,s)}\right)<\infty$ for
all $q>0$. Then:\\
{\rm (a)}
$\xi$-a.s. 
\begin{equation}
\label{boxsupremum}
\limsup_{t \to \infty} \left[\frac{1}{\log t}\sup\left\{\xi(x,s)\colon\,x \in [C_1]_t, 
0\leq s\leq t\right\}\right] \leq 1.
\end{equation}
{\rm (b)} 
$\xi$-a.s.\ there exists a $t_0\geq 0$ such that, for all $t\geq t_0$ and $x \notin [C_1]_t$,
\begin{equation}
\label{outersupremum}
\sup_{s \in [0,t]} \xi(x,s) \leq \log\|x\|.
\end{equation}
\el

\bpr
(a) For any $\theta>0$ and $t \geq 1$, we may estimate
\begin{equation}
\label{boxsupest1}
\begin{aligned}
&\P\left(\exists\,x \in [C_1]_t\colon\, \sup_{s\in[0,t]}\xi(x,s) \geq \log t\right)\\
&\leq \sum_{x \in [C_1]_t}\sum_{k=0}^{\lfloor t \rfloor} 
\P\left(\sup_{s\in[k,k+1]} \xi(x,s) \geq \log t\right)\\
&\leq (2C_1t\log t +1)^{d} (\lfloor t  \rfloor +1) \exp\{-\theta \log t\}
\E\left(\exp\left\{\theta \sup_{s \in [0,1]} \xi(0,s)\right\}\right).
\end{aligned}
\end{equation}
Choosing $\theta >2(d+1)+1$, we get that the right-hand side is summable over $t \in \N$. 
Hence, by the Borel-Cantelli Lemma, we get the claim.

\medskip\noindent
(b) The proof is similar and is omitted.
\epr

\noindent
The main result of this section reads:

\bl{outerlem}
There exists a $C_0>0$ such that $\xi$-a.s.\ there exists a $t_0 >0$ such that
\begin{equation}
\label{outerFK}
E_0\left(e^{\cI^{\kappa}(0,t,t)}
\one{\{X^{\kappa}([0,t]) \not\subseteq [C_1]_t}\}\right) \leq e^t
\qquad \forall\,t\geq t_0, C_1\geq C_0.
\end{equation}
\el

\bpr
See Kesten and Sidoravicius~\cite{KS03}, Eq.\ (2.38). We only sketch the main idea.
Take a realization $\Phi\colon\,[0,t]\to\Z^d$ of a random walk path that leaves the 
box $[C_1]_t$. Then $\|\Phi\| = \max\{\|x\|\colon\,x\in\Phi([0,t])\} > C_1 t\log t$. By 
Lemma \ref{supremumlem}, 
\begin{equation}
\sup_{s \in [0,t]} \sup_{\|x\|\leq \|\Phi\|} \xi(x,s) \leq \log \|\Phi\|,
\end{equation}
and so we can estimate
\begin{equation}
\label{outerFKest3}
\begin{aligned}
&E_0\left(e^{\cI^{\kappa}(0,t,t)}
\one{\{X^{\kappa}([0,t]) \not\subseteq [C_1]_t}\}\right)\\
& \leq E_0\left(\exp\left\{t \sup_{s \in [0,t]} \log \|X^{\kappa}(s)\|\right\}
\one{\{X^{\kappa}([0,t]) \not\subseteq [C_1]_t}\}\right).
\end{aligned}
\end{equation}
The rest of the proof consists of balancing the exponential growth of the term with
the supremum against the superexponential decay of $P_0(X^{\kappa}([0,t])\not\subseteq 
[C_1]_t)$. See \cite{KS03} for details.
\epr


\subsection{Step 2: No bad $R$-blocks for large $R$}
\label{S3.3}

\bl{exclusionlem} 
Fix $C_0 >0$ according to Lemma~{\rm \ref{outerlem}}, and suppose that $\xi$ satisfies 
condition {\rm (a3)} in the G\"artner-positive-hyper-mixing definition. Then for every $C_1 \geq 
C_0$ and $\varepsilon>0$ there exists an $A=A(\varepsilon)>2$ such that 
\begin{equation}
\label{exclusionest}
\begin{aligned}
\P\Big(\Xi_{R}^{A,j} >0
\mbox{ for some }R\geq \varepsilon \log t \mbox{ and some } j\in\N_0\Big)
\end{aligned}
\end{equation}
is summable over $t\in \N$. (It suffices to choose $A=\lfloor e^{1/a(1+2d)\varepsilon}\rfloor$ 
for some $a>1$.) 
\el

\bpr
Fix $C_1\geq C_0$, $A>2$ and assume that $\Xi_{R}^{A,j}>0$ for some $j\in\N_0$. Then there is 
a bad $R$-block $B_R^{A}(x,k)$ that intersects $[C_1]_t\times[0,t]$. Hence there is a pair 
$(y,s) \in \Z^d\times[0,\infty)$ such that $Q_R^{A}(y)\times\{s\} \subseteq \tilde{B}_R^{A}
(x,k;b,c)$ and 
\begin{equation}
\label{firstestlem3}
\sum_{z\in Q_R^{A}(y)}\xi(z,s)> CA^{Rd}.
\end{equation}
In particular, $x$ and $s$ satisfy $\dist(y,[C_1]_t)\leq (b+2)A^{R}$ and $s\in [0,t+A^{R}]$. 
Hence, for $\varepsilon >0$,
\begin{equation}
\label{secondestlem3}
\begin{aligned}
&\P\Big(\Xi_{R}^{A,j} >0\mbox{ for some } R \geq \varepsilon\log t, j\in\N_0\Big)\\
&\leq \sum_{R \geq \varepsilon\log t}
\P\Bigg(\sum_{z\in Q_{R}^{A}(y)}\xi(z,s)> CA^{Rd}\mbox{ for some }(y,s)\colon\\ 
&\qquad\qquad\qquad\dist(y,[C_1]_t)\leq (b+2)A^{R}, s\in [0,t+A^{R}]\Bigg)\\
&\leq \sum_{R\geq \varepsilon\log t}\,\sum_{y:\,\dist(y,[C_1]_t)\leq (b+2)A^{R}}
\sum_{k=0}^{\lfloor t+A^{R} \rfloor} \P\Big(\exists s\in[k,k+1)\colon\,
\sum_{z\in Q_{R}^{A}(y)}\xi(z,s)>CA^{Rd}\Big).
\end{aligned}
\end{equation}
By assumption (\ref{condexpgrowthB}), we may bound the two inner sums by
\begin{equation}
\label{thirdestlem3}
(2C_1t\log t +1+(b+2)A^{R})^{d}\times(\lfloor t+A^{R} \rfloor +1)\times (2A^{R}+1)^{-d\alpha}
\stackrel{\mathrm{def}} = G(R,t).
\end{equation}
Recall the definition of $\alpha$ (see below \ref{condexpgrowthB})), to see that one can
choose $A$ as described in the formulation of Lemma \ref{exclusionlem} to get that
\begin{equation}
\label{finalestlem3}
\sum_{R\geq \varepsilon \log t} G(R,t)
\end{equation}
is summable over $t\in \N$.
\epr


\subsection{Step 3: Estimate of the Feynman-Kac formula in terms of bad blocks}
\label{S3.4}

\bl{FKbadblock}
Fix $\varepsilon >0$ and $A>2$. For all $C_1\geq C_{0}$ (where $C_0$ is determined by 
Lemma~{\rm \ref{outerlem}}),
\begin{equation}
\label{FKbadblockest}
\begin{aligned}
&E_{0}\left(e^{\cI^{\kappa}(0,t,t)}
\one{\{X^{\kappa}([0,t])\subseteq[C_1]_t\}}\right)\\
&\leq \sum_{j\in\N_0} \frac{(2dt\kappa)^{j}}{j!}
\exp\left\{t(CA^{d}-2d\kappa) + \sum_{R=1}^{\infty}CA^{(R+1)d}A^{R}\Xi_{R}^{A,j}\right\}.
\end{aligned}
\end{equation} 
\el

\bpr
See \cite{KS03}, Lemma 9. We sketch the proof. Note that
\begin{equation}
\label{FKsum}
\begin{aligned}
&E_0\left(e^{\cI^{\kappa}(0,t,t)}\one{\{X^{\kappa}([0,t])\subseteq[C_1]_t\}}\right)
= \sum_{j\in \N_0}e^{-2dt\kappa}\frac{(2dt\kappa)^j}{j!}\\
&\qquad\times \sum_{x_1,x_2,\dots ,x_j \in \Z^d} \frac{1}{(2d)^j}
E_0\left(\exp\left\{\sum_{i=1}^{j}\int_{S_{i-1}}^{S_i}\xi(x_{i-1},t-u)\, du 
+ \int_{S_j}^{t}\xi(x_j,t-u)\, du\right\}\right),
\end{aligned}
\end{equation}
where $j$ is the number of jumps, $0=x_0,x_1,\dots,x_j$, $x_i \in[C_1]_t, i\in\{0,1,\cdots,
j\}$, are the nearest-neighbor sites visited, and $0=S_0<S_1<\dots< S_j<t$ are the jump 
times. To analyze (\ref{FKsum}), fix $A>2$, $R\in \N$ as well as $0=s_0<s_1<\dots<s_j$ 
and a path $\Phi$ with these jump times, and define
\begin{equation}
\label{LambdaR}
\begin{aligned}
\Lambda_{R}(\Phi) = &\bigcup_{i=1}^{j} \left\{u \in [s_{i-1},s_{i})\colon\, 
CA^{Rd} < \xi(x_{i-1},t-u)\leq CA^{(R+1)d}\right\} \\
&\qquad \bigcup \left\{u\in [s_j,t)\colon\, CA^{Rd}< \xi(x_j,t-u)\leq CA^{(R+1)d}\right\}.
\end{aligned}
\end{equation}
The contribution of $\Phi$ to the exponential in (\ref{FKsum}) may be bounded from 
above by
\begin{equation}
\label{contribution}
tCA^d + \sum_{R=1}^{\infty}CA^{(R+1)d}|\Lambda_{R}(\Phi)|,
\end{equation}
where the first term comes from the space-time points $(x_{i-1},t-u)$ with $\xi(x_{i-1},t-u)
\leq CA^d$. If $CA^{Rd} < \xi(x_{i-1},t-u) \leq CA^{(R+1)d}$, then $(x_{i-1},t-u)$ belongs 
to a bad $R$-block. There are at most $\Xi_{R}^{A,j}$ such blocks, and any path spends at 
most a time $A^{R}$ in each $R$-block. Hence
\begin{equation}
\label{LambdaRest}
|\Lambda_{R}(\Phi)|\leq A^{R}\Xi_{R}^{A,j}.
\end{equation}
The claim now follows from (\ref{FKsum}), (\ref{contribution}--\ref{LambdaRest}) and the fact 
that there are at most $(2d)^j$ nearest-neighbor paths $(0=x_0,x_1,x_2,\dots, x_j)$ that are 
contained in $[C_1]_t$.
\epr


\subsection{Step 4: Bound on the number of bad blocks}
\label{S3.5}

The goal of this section is to provide a bound on the number of bad blocks on all scales
simultaneously (Lemma \ref{multiscalelem} below). In Section~\ref{S3.6} we will combine
Lemmas~\ref{outerlem} and \ref{FKbadblock}--\ref{multiscalepsi} to prove 
Theorem~\ref{expgrowthres}.

\bl{multiscalelem}
Fix $\varepsilon>0$, pick $A$ according to Lemma~{\rm \ref{exclusionlem}}
and assume that $\Xi_{R}^{A,j}= 0$ for all $R\geq \lceil \varepsilon \log t\rceil$. 
Then, for some $C_2>0$,
\begin{eqnarray}
\label{multiscalepsi}
&\P\Big(\Psi_{R}^{A,j} \geq (t+j)(A^{(1+2d)})^{-R}
\mbox{ for some } R\in\N \mbox{ and some } j\in\N_0\Big),\\
\label{multiscalephi}
&\P\Big(\Xi_{R}^{A,j} \geq C_2 (t+j)(A^{(1+2d)})^{-R} 
\mbox{ for some } R\in\N \mbox{ and some } j\in\N_0\Big),
\end{eqnarray}
are summable on $t\in\N$. 
\el

\noindent
The proof of Lemma~\ref{multiscalelem} is based on Lemmas \ref{mainlem}--\ref{badblocklem}
below. The first estimates for fixed $R$ the probability that there is a large number of 
good $(R+1)$-blocks containing a bad $R$-block, the second gives a recursion bound on the 
number of bad blocks in terms of $\Psi_{R}^{A,j}$.

\bl{mainlem}
Suppose that $\xi$ satisfies condition {\rm (a2)} in Definition~{\rm \ref{Gartnerhypmix}}. 
Then, for $R$ large enough, $j \in \N_0$ and $A$ chosen according to Lemma~\ref{exclusionlem}, for some constant $C_3>0$
\begin{equation}
\label{mainest}
\P\left(\Psi_{R}^{A,j} \geq (t+j)(A^{(1+2d)})^{-R}\right) 
\leq \exp\left\{-C_3(t+j)(A^{(1+2d)})^{-R}\right\}.
\end{equation}
\el

\bl{badblocklem}
Fix $\varepsilon >0$, and pick $A$ according to Lemma~{\rm \ref{exclusionlem}}. 
Assume that $\Xi_{R}^{A,j}= 0$ for all $R\geq \lceil \varepsilon \log t\rceil$.
Then, with $N=\lceil \varepsilon\log t\rceil$,
\begin{equation}
\label{badblockest}
\Xi_{R}^{A,j}\leq 2^dA^{(1+d)}\sum_{i=0}^{N-R-1}2^{id}A^{i(1+d)}\Psi_{R+i}^{A,j}.
\end{equation}
\el

\noindent
The proofs of Lemmas~\ref{mainlem}, \ref{badblocklem} and \ref{multiscalelem}
are given in Sections~\ref{S3.5.1}, \ref{S3.5.2} and \ref{S3.5.3}, respectively.


\subsubsection{Proof of Lemma~\ref{mainlem}}
\label{S3.5.1}

\bpr
Throughout the proof, $x,x'\in\Z^d$ and $k,k' \in \N$. The idea of the proof is to divide 
space-time blocks into equivalence classes, such that in each equivalence class blocks 
are far enough away from each other so that they can be treated as being independent (see 
Fig.~\ref{pictFiniteness2}). The proof comes in four steps.

\medskip\noindent
{\bf 1.} Fix $A>2$ according to Lemma \ref{exclusionlem}, fix $R\in \N$ and take $a_1,a_2,b,c
\in \N_0$ according to condition {\rm (a2)} in Definition~\ref{Gartnerhypmix}. We say that 
$(x,k)$ and $(x',k')$ are equivalent if and only if
\begin{equation}
\label{equiv}
x\equiv x' \mod a_1 \quad \mbox{ and } \quad k\equiv k' \mod a_2.
\end{equation}
This equivalence relation divides $\Z^d\times\N$ into $a_1^da_2$ equivalence classes.
We write $\sum_{(x^*,k^*)}$ to denote the sum over all equivalence classes. Furthermore, 
we define
\begin{equation}
\label{defchi}
\begin{aligned}
\chi^{A}(x,k) = &\one{\left\{B_{R+1}^{A}(x,k)\mbox{ is good, but contains a 
bad $R$-block}\right\}}.
\end{aligned}
\end{equation}
We tacitly assume that all blocks under consideration intersect $[C_1]_t\times[0,t]$. 
Then
\begin{equation}
\label{firstestpsilem}
\begin{aligned}
&\P\left(\Psi_{R}^{A,j} \geq (t+j)(A^{(1+2d)})^{-R}\right)\\
&\leq \sum_{(x^*,k^*)}\P\Bigg(\begin{array}{ll}&\exists\mbox{ a path with $j$ jumps that intersects
at least }\\
&(t+j)(A^{(1+2d)})^{-R}/a_1^{d}a_2
\mbox{ blocks $B_{R+1}^{A}(x,k)$}\\
&\mbox{ with $\chi^{A}(x,k)=1$, $(x,k)\equiv (x^*,k^*)$}\end{array}\Bigg).
\end{aligned}
\end{equation}

\medskip\noindent
{\bf 2.}
For the rest of the proof we fix an equivalence class and
define $\rho_{R}^{1/(1+d)} = (A^{(1+2d)})^{-R}$ (recall \ref{wGartner}).
To control the number of different ways to cross a prescribed number
of $R$-blocks we consider enlarged blocks.
To that end, as in \cite{KS03}, take $\nu
= \lceil \rho_{R}^{-1/(1+d)} \rceil$ and define space-time blocks
\begin{equation}
\label{Btildeblock}
\begin{aligned}
\tilde{B}_{R}^{A}(x,k) = \left(\prod_{j=1}^{d}
[\nu(x(j)-1)A^{R},\nu(x(j)+1)A^{R})\cap\Z^d\right)
\times[\nu kA^{R},\nu(k+1)A^{R}).
\end{aligned}
\end{equation}
By the same reasoning as in \cite{KS03}, we see that at most 
\begin{equation}
\label{muj}
\mu(j)\stackrel{\mathrm{def}}=3^{d}\left(\frac{t+j}{\nu A^{R+1}}+2\right)
\end{equation}
blocks $\tilde{B}_{R+1}^{A}(x,k)$ can be crossed by a path $\Phi$ with $j$ jumps. 
We write 
\begin{equation}
\label{eq:allseq}
\bigcup_{(x_i,k_i)} \tilde{B}_{R+1}^{A}(x_i,k_i)
\quad \mbox{and} \quad \sum_{\tilde{B}_{R+1}^{A}(x_i,k_i)}
\end{equation}
to denote the union over at most $\mu(j)$ blocks $\tilde{B}_{R+1}^{A}(x_i,k_i)$, $0\leq i\leq \mu(j)-1$, and
to denote the sum over all possible sequences of blocks $\tilde{B}_{R+1}^{A}(x_i,k_i)$,
that may be crossed by a path $\Phi$ with $j$ jumps, respectively.
As each block $B_{R+1}^{A}(x,k)$ that may be crossed by such a path is contained in 
the union in (\ref{eq:allseq}), we may estimate the probability in (\ref{firstestpsilem}) from above by
\begin{equation}
\label{upbound}
\begin{aligned}
 \sum_{\tilde{B}_{R+1}^{A}(x_i,k_i)}
\P\bigg(\begin{array}{ll}&\mbox{the union in (\ref{eq:allseq})
contains at least $(t+j)\rho_{R}^{1/(1+d)}/a_1^{d}a_2$}\\
&\mbox{blocks $B_{R+1}^{A}(x,k)$
with $\chi^{A}(x,k)=1$, $(x,k)\equiv (x^*,k^*)$ }\end{array}\bigg).
\end{aligned}
\end{equation}
To estimate the probability in (\ref{upbound})
write
\begin{equation}
\label{abbrev1}
\begin{aligned}
\mathcal{A}^n((x_0,k_0),\ldots,(x_{\mu(j)-1},k_{\mu(j)-1})) 
= \bigg\{\hspace{-.4cm}\begin{array}{ll}&\mbox{the union in (\ref{eq:allseq}) contains $n$ blocks
$B_{R+1}^{A}(x,k)$}\\
&\mbox{with $\chi^{A}(x,k)=1$, $(x,k)\equiv (x^*,k^*)$}\end{array}\hspace{-.2cm}\bigg\}.\\
\end{aligned}
\end{equation}
Since the union (\ref{eq:allseq}) contains at most $L= \nu^{(1+d)}\mu(j)$ blocks
$B_{R+1}^{A}(x,k)$, the probability in (\ref{upbound}) are bounded from above by
\begin{equation}
\label{eq:partition}
\sum_{n=\frac{(t+j)\rho_{R}^{1/(1+d)}}{a_1^da_2}}^{L}
\P\Big(\mathcal{A}^n((x_0,k_0),\ldots,(x_{\mu(j)-1},k_{\mu(j)-1})) \Big).
\end{equation}
Note that there are $\binom{L}{n}$ ways of choosing $n$ blocks
$B_{R+1}^{A}(x,k)$ with $\chi^{A}(x,k)=1$ out of $L$ $(R+1)$-blocks.
Hence, by condition {\rm (a2)} in Definition~\ref{Gartnerhypmix}, 
(\ref{eq:partition}) is at most
\begin{equation}
\label{eq:binomial}
\sum_{n=\frac{(t+j)\rho_{R}^{1/(1+d)}}{a_1^da_2}}^{L}
\binom{L}{n} \rho_{R}^{n}
\leq \Big(1-\rho_{R}\Big)^{-L} K 
\P\Bigg( T_{L}\geq \frac{(t+j)\rho_{R}^{1/(1+d)}}{a_1^da_2}\Bigg),
\end{equation}
where $T_L= \mathrm{BIN}(L,\rho_{R})$.

\medskip\noindent
{\bf 3.}
To estimate the binomial random variable, note that by Bernstein's inequality 
(compare with \cite{KS03}, Lemma 11), there is a constant $C'$ such that, for all 
$\lambda \geq 2\E(T)$,
\begin{equation}
\label{Bernstein}
\P(T_L\geq \lambda)\leq \exp\{-C'\lambda\}.
\end{equation}
We may assume that $\rho_{R}^{-1/(1+d)} \in \N$, so that
\begin{equation}
\label{Texpectation}
\E(T_L) = \nu^{(1+d)}\mu(j)\rho_{R}
= 3^{d}\left(\frac{t+j}{A^{R+1}}\,\rho_{R}^{1/(1+d)}+2\right),
\end{equation}
and hence, by Lemma \ref{exclusionlem} and the fact that $R\leq \varepsilon \log t$,
\begin{equation}
\label{Tcompare}
\frac{(t+j)\rho_{R}^{1/(1+d)}}{a_1^{d}a_2}\geq 2\E(T_L).
\end{equation}
Since $a_1,a_2$ are independent of $R$, we may estimate, using (\ref{Bernstein}),
\begin{equation}
\label{appofBernstein}
\begin{aligned}
\P\left(T_L\geq \frac{(t+j)\rho_{R}^{1/(1+d)}}{a_1^{d}a_2}\right)
&\leq \exp\left\{-C'(t+j)\rho_{R}^{1/(1+d)}\right\} \\
&= \exp\left \{-C'(t+j)(A^{(1+2d)})^{-R}\right\}.
\end{aligned}
\end{equation}
It rests to show that the first term on the right hand side in (\ref{eq:binomial}) does not contribute.
Note that $1/(1-\rho_{R}) = 1+ \rho_{R}/(1-\rho_{R})$, so that
\begin{equation}
\label{logest}
\log\bigg(\frac{1}{1-\rho_R}\bigg)\leq \frac{\rho_{R}}{1-\rho_R}.
\end{equation} 
Thus, if we assume $\rho_r^{-1/(1+d)}\in\N$,
\begin{equation}
\label{restest}
L\log\bigg(\frac{1}{1-\rho_R}\bigg) \leq \frac{\mu(j)}{1-\rho_R}.
\end{equation}
Inserting (\ref{restest}) into the second term on the right hand side
of (\ref{eq:binomial}), comparing it with the right hand side of (\ref{appofBernstein}), 
and recalling the definition of $\mu(j)$ in (\ref{muj}),
we see that the asymptotic of (\ref{eq:binomial}) is determined by the 
probability term in (\ref{appofBernstein}).

\medskip\noindent
{\bf 4.} 
Finally, we estimate (\ref{firstestpsilem}). 
\begin{claim}
\label{cl:cardinality}
There is $C>0$ such that the number of summands in (\ref{upbound})
is bounded from above by $e^{C\mu(j)}$.
\end{claim}
Before we proof the claim, we show how one deduces Lemma \ref{mainlem} from it. Insert (\ref{appofBernstein})
into (\ref{upbound}) and use Claim \ref{cl:cardinality} to obtain that
\begin{equation}
\label{fourthestpsilem}
\begin{aligned}
\P\left(\Psi_{R}^{A,j} \geq (t+j)(A^{(1+2d)})^{-R}\right)
\leq Ka_1^da_2 \exp\left \{-C'(t+j)(A^{(1+2d)})^{-R}\right\}.
\end{aligned}
\end{equation}
We now prove Claim \ref{cl:cardinality}.
\bpr
Assume that $d=1$. Divide time 
into intervals of length $\nu A^{R+1}$ and fix an integer-valued sequence\\ $(l_1,l_2,\ldots,l_{t/\nu A^{R+1}})$ 
such that
\begin{equation}
\label{eq:lsum}
\sum_{i=1}^{t/\nu A^{R+1}}l_i\leq \mu(j).
\end{equation}
We first estimate the number of ways in which $\Phi$ can visit $l_i$ $R$-blocks $\tilde{B}_{R+1}^{A}(x,k)$ 
in the $i$-th time interval for $i=1,\ldots,t/\nu A^{R+1}$. Let $\tilde{B}_{R+1}^{A}(x_0,k_0),\tilde{B}_{R+1}^{A}(x_1,k_1),
\ldots,\tilde{B}_{R+1}^{A}(x_{l_i-1},k_{l_i-1})$ be a possible sequence. Then $x_0=0$, $x_1=\pm 1$, 
$x_2 \in (-1,x_1+1)$ if $x_1=1$ or $x_2 \in (x_1-1,1)$ if $x_1=-1$, etc. For each $j$ such 
that $l_i+1< j \leq l_{i+1}$, $x_j$ can take two values depending on $x_{j-1}$. Note 
that if $j= l_i +1$ for some $j$ and $i$, then $x_j$ is the space label of a block at the 
beginning of a new time interval. Consequently, $x_j\in\{x_{j-1}-1,x_{j-1},x_{j-1}+1\}$.
Hence, for a fixed choice of $(l_1,l_2,\ldots,l_{t/\nu A^{R+1}})$, there are at most $3^{l_1}
\times 3^{l_2}\times\cdots\times 3^{l_{t/\nu A^{R+1}}}\leq 3^{\mu(j)}$ possibilities to choose a sequence 
of blocks with a prescribed sequence $(l_1,l_2,\ldots,l_{t/\nu A^{R+1}})$. The above arguments, 
together with the fact that for some $a,b \in (0,\infty)$ there are no more than $(a/l) 
e^{b\sqrt{l}}$ such sequences $(l_1,l_2,\ldots, l_{t/\nu A^{R+1}})$ (see Hardy and Ramanujan \cite{HR18} and Erd\"os
\cite{E42}), yield the claim. 
The extension to $d \geq 2$ is straightforward.
\epr
This finishes the proof of Lemma \ref{mainlem}.
\epr


\subsubsection{Proof of Lemma \ref{badblocklem}}
\label{S3.5.2}

\bpr
We first show that
\begin{equation}
\label{recursive1}
\Xi_{R}^{A,j} \leq 2^dA^{(1+d)}\Xi_{R+1}^{A,j} + 2^dA^{(1+d)}\Psi_{R}^{A,j}.
\end{equation}
In order to see why (\ref{recursive1}) is true, take a bad $R$-block $B_{R}^{A}(x,k)$ that 
is crossed by a path with $j$ jumps. Then there are two possibilities. Either $B_{R}^{A}(x,k)$ 
is contained in a bad $(R+1)$-block, or all $(R+1)$-blocks that contain $B_{R}^{A}(x,k)$ 
are good. Since an $(R+1)$-block contains $A^{(1+d)}$ $R$-blocks, and there are at most 
$2^d$ $(R+1)$-blocks, which may contain a given $R$-block, the first term in the above sum 
bounds the number of bad $R$-blocks contained in a bad $(R+1)$-block. In contrast, the 
second term bounds the number of bad $R$-blocks contained in a good $(R+1)$-block. Hence 
we obtain (\ref{recursive1}).

We can now prove the claim. Apply (\ref{recursive1}) iteratively to the terms in the sum, 
i.e., replace $\Xi_{R+i}^{A,j}$ by
\begin{equation}
\label{recursive2}
2^dA^{(1+d)}\Xi_{R+i+1}^{A,j} +2^d A^{(1+d)}\Psi_{R+i}^{A,j}.
\end{equation}
This yields
\begin{equation}
\label{appofrecursion}
\Xi_{R}^{A,j} \leq 2^dA^{(1+d)}\sum_{i=0}^{N-R-1}2^{id}A^{i(1+d)}\Psi_{R+i}^{A,j},
\end{equation}
from which the claim follows.
\epr


\subsubsection{Proof of Lemma \ref{multiscalelem}}
\label{S3.5.3}

\bpr
Fix $\varepsilon>0$ 
and $0<R\leq \varepsilon\log t$. Then, by Lemma~\ref{mainlem},
\begin{equation}
\label{multiscalepsiest1}
\begin{aligned}
&\P\left(\Psi_{R}^{A,j} \geq 
(t+j)(A^{(1+2d)})^{-R} \mbox{ for some $j\in\N_0$}\right)\\
&\leq \sum_{j\in\N_0} \P\left(\Psi_{R}^{A,j} 
\geq (t+j)(A^{(1+2d)})^{-R}\right)\\
&\leq \sum_{j\in\N_0} \exp\left\{-tC_3(t+j)(A^{(1+2d)})^{-R}\right\}\\
&\leq \exp\left\{-C_3t(A^{(1+2d)})^{-\varepsilon \log t}\right\} \sum_{j\in\N_0}
\exp\left\{-C_3j(A^{(1+2d)})^{-\varepsilon \log t}\right\}.
\end{aligned}
\end{equation}
Recall Lemma \ref{exclusionlem}, which implies that $\delta\stackrel{\mathrm{def}}
=\log (A)\varepsilon(1+2d)<1$. Consequently, the right-hand side of (\ref{multiscalepsiest1}) 
is at most
\begin{equation}
\label{multiscalepsiest2}
\begin{aligned}
&\exp\left\{-C_3 t^{(1-\delta)}\right\} 
\,\frac{1}{1-\exp\left\{-C_3t^{-\delta}\right\}}\\
& \leq \frac{1}{C_3}\exp\left\{-C_3t^{(1-\delta)}\right\}
\,t^{\delta}\,\exp\left\{C_3t^{-\delta}\right\}.
\end{aligned}
\end{equation}
It therefore follows that
\begin{equation}
\label{multiscalepsiestfinal}
\begin{aligned}
&\P\left(\Psi_{R}^{A,j} \geq (t+j)(A^{(1+2d)})^{-R}
\mbox{ for some $j\in\N_0$, $R\in\N$}\right)\\
&\leq \frac{1}{C_3}t^{\delta}\,\varepsilon\log t\,
\exp\left\{-C_3t^{(1-\delta)} + C_3t^{-\delta}\right\},
\end{aligned}
\end{equation}
which is summable over $t\in \N$. In order to prove the second statement, suppose that none 
of the events in (\ref{multiscalepsi}) occurs. With Lemma~\ref{badblocklem} we may estimate
\begin{equation}
\label{multiscalephiest1}
\begin{aligned}
\Xi_{R}^{A,j}&\leq 2^dA^{(1+d)}\sum_{i=0}^{N-R-1}2^{id}A^{i(1+d)}\Psi_{R+i}^{A,j}\\
&\leq 2^dA^{(1+d)}\sum_{i=0}^{N-R-1}(t+j)2^{id}A^{i(1+d)}(A^{(1+2d)})^{-i-R}\\
&\leq 2^dA^{(1+d)}(t+j)A^{-R(1+2d)}\sum_{i\in\N_0} 2^{id}A^{-id}\\
&\stackrel{\mathrm{def}}= (t+j)A^{-R(1+2d)}C_2,
\end{aligned}
\end{equation}
where we use that $A>2$ (see Lemma~\ref{exclusionlem}). 
\epr


\subsection{Step 5: Proof of Theorem \ref{expgrowthres}}
\label{S3.6}

Fix $\varepsilon>0$ and $A$ such that Lemma \ref{multiscalelem} applies. It follows from 
Lemma~\ref{exclusionlem} that
\begin{equation}
\label{appofexclusionlem}
\begin{aligned}
\P\Big(\Xi_{R}^{A,j} >0
\mbox{ for some }R\geq \varepsilon \log t, j\in\N_0\Big)
\end{aligned}
\end{equation}
is summable over $t\in \N$. Hence, by the Borel-Cantelli Lemma, there is an $t_0\in\N$ 
such that none of the events in the above probability occurs for integer $t\geq t_0$.
Thus, by Lemma~\ref{FKbadblock}, for all integer $t\geq t_{0}$ we have, with 
$N=\lfloor \varepsilon \log t\rfloor$,
\begin{equation}
\label{appFKbadblock}
\begin{aligned}
&E_{0}\left(e^{\cI^{\kappa}(0,t,t)}
\one{\big\{X^{\kappa}([0,t])\subseteq[C_1]_t\big\}}\right)\\
&\leq \sum_{j\in\N_0} \frac{(2dt\kappa)^{j}}{j!}\exp\left\{t(CA^{d}-2d\kappa) + 
\sum_{R=1}^{N}CA^{(R+1)d}A^{R}\Xi_{R}^{A,j}\right\}.
\end{aligned}
\end{equation}
Using the bound of Lemma~\ref{multiscalelem}, we have
\begin{equation}
\label{appomultiscalelem}
\begin{aligned}
\sum_{R=1}^{N}CA^{(R+1)d}A^{R}\Xi_{R}^{A,j}
&\leq \sum_{R=1}^{N}C(t+j)A^{(R+1)d}A^{R}A^{-R(1+2d)}C_2\\
&\leq (t+j)A^{d}C'\sum_{R\in\N} A^{-Rd} \leq C_4(t+j).
\end{aligned}
\end{equation}
We can therefore estimate the last line of (\ref{appFKbadblock}) by
\begin{equation}
\label{finalestthm}
\begin{aligned}
&\sum_{j\in\N_0} \frac{(2dt\kappa)^{j}}{j!}\exp\left\{t(CA^{d}-2d\kappa)+ C_4(t+j)\right\}\\
& = \exp\left\{t(CA_2^{d}-2d\kappa + C_4 + 2d\kappa e^{C_4})\right\}.
\end{aligned}
\end{equation}
From (\ref{finalestthm}) and Lemma~\ref{outerlem}, we obtain 
\begin{equation}
\label{discretelimsup}
\limsup_{\substack{t\to\infty\\ t\in\N}}\frac{1}{t}\log 
E_{0}\left(e^{\cI^{\kappa}(0,t,t)}\right) < \infty.
\end{equation}
To extend this to sequences along $\R$ instead of $\N$, note that
\begin{equation}
\label{exttoR}
\begin{aligned}
u(0,t) \leq u(0,n+1)\,e^{-\int_{t}^{n+1}\xi(0,s)\,ds}e^{2d\kappa (n+1-t)},
\qquad t\in[n,n+1].
\qquad 
\end{aligned}
\end{equation}
Since $\xi$ is ergodic in time, we have 
\begin{equation}
\label{apptimeergodic}
\lim_{t \to \infty}\frac{1}{t}\int_{t}^{\lceil t \rceil}\xi(0,s)\,ds = 0.
\end{equation}
Theorem~\ref{expgrowthres} follows from (\ref{exttoR}--\ref{apptimeergodic}).


\section{Initial condition}
\label{S4}

In this section we prove Theorem~\ref{nonlocLyp}. Section~\ref{S4.1} contains some preparations.
Section~\ref{S4.2} states three lemmas (Lemmas~\ref{nonloclem1}--\ref{nonloclem3} below) that 
are needed for the proof of Theorem~\ref{nonlocLyp}, which is given in Section~\ref{S4.3}. 
Section~\ref{S4.4} provides the proof of these three lemmas. 


\subsection{Preparations}
\label{S4.1}

In this section we first state and prove a lemma (Lemma~\ref{summablelem} below) that will be 
needed for the proof of Theorem~\ref{nonlocLyp}. After that we introduce some further notation
(Definitions~\ref{blockdef}--\ref{sufficientblock} below).

Fix $R_0 \in\N$ and take $A,C$ according to our assumption (type-II G\"artner-mixing). Set $N=CA^{R_0d}$ and abbreviate 
$\xi_N=(\xi\wedge N)\vee (-N)$. Let $u_N$ be the solution of (\ref{pA}) with $\xi$ replaced 
by $\xi_N$. Abbreviate (recall (\ref{timerevIntegral}--\ref{timeIntegral}))
\begin{eqnarray}
\label{timerevIntegralN}
\cI_N^{\kappa}(a,b,c) &= \int_a^b \xi_N(X^\kappa(s),c-s)\,ds, \qquad 0 \leq a \leq b \leq c,\\
\label{timeIntegralN}
\overline{\cI}_N^{\kappa}(a,b,c) &= \int_a^b \xi_N(X^\kappa(s),c+s)\,ds, 
\qquad 0 \leq a \leq b \leq c.
\end{eqnarray}

\bl{summablelem}
If, for all $N$ of the form $N=CA^{R_0d}$ and for all $\varepsilon >0$ and some sequence 
$(t_{r})_{r\in\N}$ of the form $t_r=rL$ with $L>0$,
\begin{equation}
\label{summablecond}
\P\left(E_0\left(e^{\overline{\cI}_N^{\kappa}(0,t_r,0)} \right) 
> e^{(\overline{\lambda}_0^{\one}(\kappa) + \varepsilon)t_r}\right)
\end{equation}
is summable on $r$, then Theorem~{\rm \ref{nonlocLyp}} holds.
\el

\bpr

Fix $\varepsilon >0$. Note that $u_N(0,t)$ has the same distribution as $E_0(e^{\overline
{\cI}_N^{\kappa}(0,t,0)})$, so that we can replace the latter by $u_N(0,t)$ in 
(\ref{summablecond}) without violating the summability condition. Thus, by the 
Borel-Cantelli Lemma, we have
\begin{equation}
\label{limsupiniThm}
\limsup_{r\to\infty}\frac{1}{t_r} \log E_0\left(e^{{\cI}_N^{\kappa}(0,t_r,t_r)}\right) 
\leq \overline{\lambda}_0^{\one}(\kappa) +\varepsilon \qquad \xi\mbox{-a.s.}
\end{equation}
The extension to sequences along $\R$ may be done as in the proof of 
Theorem~\ref{expgrowthres} (recall (\ref{exttoR})). Standard arguments yield 
\begin{align}
\label{limsupfin}
\limsup_{t \to \infty} \frac{1}{t}\log u_N(0,t) \leq \overline{\lambda}_0^{\one}(\kappa)
\qquad \xi\mbox{-a.s.}
\end{align}
To extend this to the solution of (\ref{pA}) with initial condition $u_0\equiv 1$, 
we estimate
\begin{equation}
\label{xidivision}
\begin{aligned}
&\int_{0}^{t}\xi(X^{\kappa}(s),t-s)ds \\
&\leq \int_0^t\xi(X^{\kappa}(s),t-s)\one{\{\xi(X^{\kappa}(s),t-s)\geq N\}}\, ds + 
\int_0^t\xi_N(X^{\kappa}(s),t-s)\, ds.
\end{aligned}
\end{equation}
Note that, by (\ref{appomultiscalelem}) and the arguments given in Lemma~\ref{FKbadblock}, 
we have, for $t\in \N$ sufficiently large,
\begin{equation}
\label{xiNest}
\sup_{\Phi \in \Pi(j,t)} \int_0^t\xi(\Phi(s),t-s)\one{\{\xi(\Phi(s),t-s)\geq N\}}\,ds
\leq (t+j)A^dC'\sum_{R=R_0}^{\infty} A^{-Rd}.
\end{equation}
Next, choose $M>1$ such that
\begin{equation}
\label{standardest}
\limsup_{t\to \infty}\frac{1}{t}\log 
E_0\left(e^{\cI^{\kappa}(0,t,t)}\one{\{N(X^{\kappa},t)>Mt\}}\right) < \overline{\lambda}_0^{\one}(\kappa).
\end{equation}
Then, by (\ref{xidivision}--\ref{xiNest}), for $t\in\N$ sufficiently large,
\begin{equation}
\label{summablelemmainest}
\begin{aligned}
&E_0\left(e^{\cI^{\kappa}(0,t,t)}\one{\{N(X^{\kappa},t)\leq Mt\}}\right)\\
&\leq \exp\left\{(M+1)tA^dC'\sum_{R=R_0}^{\infty}A^{-Rd}\right\}
E_0\left(e^{\cI_N^{\kappa}(0,t,t)}\one{\{N(X^{\kappa},t)\leq Mt\}}\right).
\end{aligned}
\end{equation}
We infer from (\ref{limsupfin}) and (\ref{standardest}--\ref{summablelemmainest}) that
\begin{equation}
\label{summablelemconcl}
\limsup_{\substack{t \to \infty\\t\in\N}} \frac{1}{t}u(0,t) 
\leq (M+1)A^dC'\sum_{R=R_0}^{\infty}A^{-Rd} + \overline{\lambda}_0^{\one}(\kappa).
\end{equation}
Taking the limit $R_0\to \infty$, $R_0\in\N$, we obtain
\begin{equation}
\label{limsuptinN}
\limsup_{\substack{t\to \infty\\ t\in\N}}\frac{1}{t}\log u(0,t) 
\leq \overline{\lambda}_0^{\one}(\kappa).
\end{equation}
The extension to sequences along $\R$ may again be done as in the proof of 
Theorem~\ref{expgrowthres} (recall (\ref{exttoR})). Furthermore, Proposition~\ref{shapeThm} 
gives $\lambda_0^{\delta_0}(\kappa) = \overline{\lambda}_0^{\delta_0}(\kappa) 
= \overline{\lambda}_0^{\one}(\kappa)$, so that
\begin{equation}
\label{nonlocloccomp}
\limsup_{t \to \infty} \frac1t \log u(0,t)\leq \lambda_0^{\delta_0}(\kappa).
\end{equation}
By monotonicity, the reverse inequality holds with the limsup replaced by the liminf. It 
follows that $\lambda_0^{\one}(\kappa)$ exists and equals $\lambda_0^{\delta_0}(\kappa)$. 
A further monotonicity argument shows that the same is true for $\lambda_0^{u_0}(\kappa)$ 
for any initial condition $u_0$ subject to \eqref{ic}.
\epr

In view of Lemma \ref{summablelem}, our target is to prove \eqref{summablecond}. We fix $M$ 
subject to (\ref{standardest}), $N$ of the form $N=CA^{R_0d}$, $\varepsilon>0$ small, and 
write $t$ as $t=rA^R$, $r,R\in \N$, $A>2$. Note that the choice of $M$ implies that it is 
enough to concentrate on path with at most $Mt$ jumps.

We proceed by introducing space-time blocks and dividing them into good blocks and bad 
blocks, respectively, into $N$-sufficient blocks and $N$-insufficient blocks (compare 
with (\ref{BRgood} and Fig.~\ref{pictFiniteness2}).

\begin{figure}[htbp]
\begin{center}
\setlength{\unitlength}{0.35cm}
\begin{picture}(22,15)(0,2)
\put(-1,5){\vector(1,0){22}}
\put(-1,5){\vector(0,1){11}}
\put(7,15){\line(1,0){7}}
\put(14,15){\line(0,-1){9}}
\put(14,6){\line(-1,0){7}}
\put(7,6){\line(0,1){2.25}}
\put(7,12.75){\line(0,1){2.25}}
\linethickness{1mm}
\put(7,8.25){\line(0,1){4.5}}
\put(22,4.8){time}
\put(-2,17){space}
\put(8,9){$B_R^{A}(x,k)$}
\put(6.3,3.5){$kA^R$}
\put(12.5,3.5){$(k+1)A^R$}
\put(-9.5,5.7){$(x(j)-1)4MA^R$}
\put(-9.5,14.5){$(x(j)+1)4MA^R$}
\end{picture}
\caption{\small The thick line is the subpedestal.}
\label{Defsufficient}
\end{center}
\end{figure}
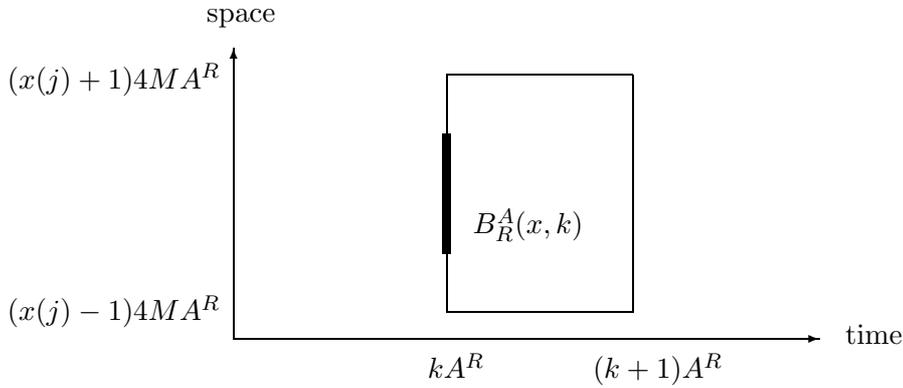

\bd{blockdef}
For $x \in \Z^d$, $k\in\N$ and $b,c \in \N_0$, define (see Fig.~\ref{Defsufficient})
\begin{equation}
\label{spacetimeblock}
\begin{aligned}
&\hat{B}_R^{A}(x,k;b,c)\\ 
&\qquad = \left(\prod_{j=1}^{d}[(x(j)-1-b)4MA^R,(x(j)+1+b)4MA^R)\cap\Z^d\right)
\times[(k-c)A^R,(k+1)A^R)
\end{aligned}
\end{equation}
and abbreviate $\hat{B}_R^{A}(x,k)= \hat{B}_R^{A}(x,k;0,0)$.
\ed

\noindent
For $S \subset \Z^d$, let $\partial S$ denote the inner boundary of $S$. For $S \times S' 
\subset \Z^d \times\R$, let $\Pi_1(S \times S')$ denote the projection of $S \times S'$ 
onto the first $d$ coordinates (the spatial coordinates).

\bd{pedestaldef}
The subpedestal of $B_R^{A}(x,k)$ is defined as
\begin{equation}
\label{subpedestal}
\begin{aligned}
\hat{B}_{R}^{A,\mathrm{sub}}(x,k)=\Big\{&y \in \Pi_1(\hat{B}_R^{A}(x,k))\colon\,\\
& |y(j)-z(j)|\geq 2MA^R,\,j\in\{1,2,\dots,d\}\,\forall\,z\in
\partial\Pi_1(\hat{B}_R^{A}(x,k))\Big\}\times \{kA^R\}.
\end{aligned}
\end{equation}
\ed

\bd{sufficientblock}
A block $\hat{B}_R^{A}(x,k)$ is called $N$-sufficient when, for every $y \in \Pi_1
(\hat{B}_{R}^{A,\mathrm{sub}}(x,k))$ (see Fig.~{\rm \ref{Defsufficient}}),
\begin{equation}
\label{good}
E_y\left(e^{\overline{\cI}_N^{\kappa}(0,A^R,kA^R)}
\one{\big\{N(X^{\kappa},A^R)\leq MA^R\big\}}\right) 
\leq e^{(\overline{\lambda}_0^{\one}(\kappa)+\varepsilon)A^R}.
\end{equation}
Otherwise $\hat{B}_R^{A}(x,k)$ is called $N$-insufficient. A subpedestal is called 
$N$-sufficient/$N$-insufficient when its corresponding block is $N$-sufficient/$N$-insufficient.
\ed

\noindent
The notion of good/bad is similar as in Definition~\ref{Cgood} with the only difference 
that $B_{R}^{A}(x,k)$ is replaced by $\hat{B}_{R}^{A}(x,k)$ and $B_{R}^{A}(x,k;b,c)$ by 
$\hat{B}_{R}^{A}(x,k;b,c)$. Similarly as in (\ref{PhiL}), define $\hat{\Xi}_{R}^{A,j}$ 
to be the maximal number of bad $R$-blocks a path with $j$ jumps can cross.

\subsection{Three lemmas}
\label{S4.2}

For the proof of Theorem~\ref{nonlocLyp} we need Lemmas~\ref{nonloclem1}--\ref{nonloclem3}  
below. The first says that each block is $N$-sufficient with a large probability (and is 
comparable with \cite{CMS02}, Lemma 4.3), the second controls the number of bad blocks
(and is comparable with Lemma \ref{multiscalelem}), the third estimates the number of $N$-insufficient 
blocks that are good and are visited by a typical random walk path (see \cite{CMS02}, 
Lemma 4.4 and \cite{KS03}, Lemma 11).

\bl{nonloclem1}
Fix $A\in\N$. For every $\tilde{\delta} >0$ there is an $R_0=R_0(A,\tilde{\delta})\in\N$ 
such that
\begin{equation}
\label{nonlocgood1}
\P(\hat{B}_R^{A}(x,k) \mbox{ is $N$-sufficient }) \geq 1-\tilde{\delta} 
\qquad \forall\,\,R\geq R_0,\,x \in \Z^d,\,k \in \N.
\end{equation}
\el

\bl{nonloclem2}
For every $A_0>2$ there is an $A\in\N$ with $A\geq A_0$ such that, for some $C>0$
independent of $A$,
\begin{equation}
\label{badblocknonloc}
\P\bigg(\hat{\Xi}_{R}^{A,j} \geq C(t+j)(A^{(1+2d)})^{-R} 
\mbox{ for some $j\in\N_0$ and $R\in\N$}\bigg)
\end{equation}
is summable on $t\in\N$.
\el

\bl{nonloclem3}
Let $C(Mt,\eta t/A^R)$ be the event that there is a path $\Phi$ with $\Phi(0)=0$ and
$N(\Phi,t) \leq Mt$ that up to time $t$ crosses more than $\eta t/A^R$ $N$-insufficient 
subpedestals of a good $R$-block. Then, under the G\"artner-mixing type-II condition, 
for every $\eta>0$ there is an $A$ (which can be chosen as in Lemma~{\rm \ref{nonloclem2}}) 
and $R_0\in\N$ such that, for some $c_1>0$,
\begin{align}
\label{nonlocgood2}
\P\left(C(Mt,\eta t/A^R)\right) \leq e^{-c_1\eta t/A^R}, \qquad \forall R\geq R_0.
\end{align}
\el


\subsection{Proof of Theorem \ref{nonlocLyp}}
\label{S4.3}

\bpr
The proof comes in three steps. Fix $0<\eta < \varepsilon$, and choose $A, R\geq R_0$ 
according to Lemmas~\ref{nonloclem1}--\ref{nonloclem3}. 

\medskip\noindent
{\bf 1.} 
Consider all random walk path that start in zero, make $0 \leq j \leq Mt$ jumps, and 
attain values $\{x_1,x_2,\dots,x_{t/N_0-1}\}$ at times $kA^R$, $k \in \{1,2,\dots,t/A^R-1\}$. 
By the Markov property,
\begin{equation}
\label{fixedpathest}
\begin{aligned}
&E_0\left(e^{\overline{\cI}_N^{\kappa}(0,t,0)}
\prod_{k=1}^{t/A^R-1}\one{\{X^{\kappa}(kA^R)=x_k\}}\right)
\leq \prod_{k=0}^{t/A^R-1}E_{x_k}\left(e^{\overline{\cI}_N^{\kappa}(0,A^R,kA^R)}\right),
\end{aligned}
\end{equation}
where $x_0=0$. Let $I$ and $S$ be the sets of indices $k$ such that $(x_k,kA^R)$ is in 
an $N$-insufficient, respectively, $N$-sufficient subpedestal. Then the right-hand side of 
\eqref{fixedpathest} can be rewritten as
\begin{equation}
\label{insufsuf}
\prod_{k \in I} E_{x_k}\left(e^{\overline{\cI}_N^{\kappa}(0,A^R,kA^R)}\right)
\prod_{k \in S} E_{x_k}\left(e^{\overline{\cI}_N^{\kappa}(0,A^R,kA^R)}\right).
\end{equation}
Because of Lemmas~\ref{nonloclem2} and \ref{nonloclem3}, there is a measurable set, 
independent of $j$ and of $\xi$-probability at least $1-e^{-c_1\eta t/A^R}$, such 
that
\begin{equation}
\label{Iest}
|I|\leq \eta t/A^{R} + C(t+j)(A^{(1+2d)})^{-R}.
\end{equation}
Since $\xi_N\leq N$, on a set of that probability the first term in \eqref{insufsuf} 
can be estimated from above by $e^{N\eta t}\exp\{NC(t+j)/A^{2Rd}\}$.

\medskip\noindent
{\bf 2.} 
Pick a realization of $\xi$ which satisfies (\ref{Iest}). To bound the second 
term in (\ref{insufsuf}), we split this term up as
\begin{equation}
\label{splittingpart1}
\begin{aligned}
\prod_{k \in S} & \bigg[ E_{x_k}\left(e^{\overline{\cI}_N^{\kappa}(0,A^R,kA^R)}
\one{\big\{N(X^{\kappa},A^R)\leq MA^R\big\}}\right)\\
& + E_{x_k}\left(e^{\overline{\cI}_N^{\kappa}(0,A^R,kA^R)}
\one{\big\{N(X^{\kappa},A^R)>MA^R\big\}}\right)\bigg],
\end{aligned}
\end{equation}
which can be written as
\begin{equation}
\label{splittingpart2}
\begin{aligned}
\sum_{J\subset S}&\bigg[\prod_{k \in J}  
E_{x_k}\left(e^{\overline{\cI}_N^{\kappa}(0,A^R,kA^R)}
\one{\big\{N(X^{\kappa},A^R)\leq MA^R\big\}}\right)\\
&\times \prod_{k \notin J} E_{x_k}
\left(e^{\overline{\cI}_N^{\kappa}(0,A^R,kA^R)}
\one{\big\{N(X^{\kappa},A^R)> MA^R\big\}}\right)\bigg].
\end{aligned}
\end{equation}
Take $c\gg1$. Then, for $M$ large enough, $P_{x_k}\left(N(X,A^R)>MA^R\right) \leq e^{-cA^R}$. 
Hence the second term in \eqref{splittingpart2} can be bounded from above by 
$e^{A^R(-c+N)(t/A^R-|J|)}$. Recall the definition of a $N$-sufficient block, to bound the 
sum in \eqref{splittingpart2} by
\begin{align}
\label{upperbound}
e^{t(-c+N)}\left(1+ e^{A^R(\overline{\lambda}_0^{\one}
(\kappa)+\varepsilon +c-N)}\right)^{t/A^R}.
\end{align}
Summing over all possible values $(x_1,x_2,\dots,x_{t/A^R-1})$ compatible with a path $\Phi$ 
such that $\Phi(0)=0$ and $N(\Phi,t)= j$, and fixing $\eta \leq \varepsilon$, we obtain
\begin{equation}
\label{prefinalest}
\begin{aligned}
&E_0\left(e^{\overline{\cI}_N^{\kappa}(0,t,0)}
\one{\{N(X^{\kappa},t)=j\}}\right)\\
&= \sum_{x_1,x_2,\dots,x_{t/A^R-1}} 
E_0\left(e^{\overline{\cI}_N^{\kappa}(0,t,0)}
\prod_{k=1}^{t/A^R-1}\one{\big\{X^{\kappa}(kA^R)=x_k\big\}}\right)\\
&\qquad \qquad \qquad\times p_{A^R}(0,x_1)\times \cdots \times 
p_{A^R}(x_{t/A^R-2},x_{t/A^R-1})\\
&\leq \exp\left\{\frac{jNC}{A^{2Rd}}\right\}
\exp\left\{t\left(N\eta + \frac{NC}{A^{2Rd}}-c+N\right)\right\} 
\left(1+ e^{A^R(\overline{\lambda}_0^{\one}(\kappa)+\varepsilon +c-N)}\right)^{t/A^R}\\
&\stackrel{\mathrm{def}}= \exp\left\{\frac{jNC}{A^{2Rd}}\right\} C_1^{N}(t,A^R,\varepsilon),
\end{aligned}
\end{equation}
where $(p_s(x,y))_{s\geq 0,x,y \in \Z^d}$ denote the transition probabilities of a 
continuous-time simple random walk jumping at rate $\kappa$. 

\medskip\noindent
{\bf 3.} 
We proceed by summing over the number of jumps, to obtain
\begin{equation}
\label{finalest}
\begin{aligned}
&E_0\left(e^{\overline{\cI}_N^{\kappa}(0,t,0)}
\one{\{N(X^{\kappa},t)\leq Mt\}}\right)\\
&\leq \sum_{j=0}^{Mt} e^{-2d\kappa t}\frac{(2d\kappa t)^j}{j!} 
\exp\left\{\frac{jNC}{A^{2Rd}}\right\} C_1^{N}(t,A^R,\varepsilon)\\
&\leq e^{-2d\kappa t} \exp\left\{2d\kappa t e^{NC/A^{2Rd}}\right\}
C_1^{N}(t,A^R,\varepsilon)
\stackrel{\mathrm{def}}=C_2^{N}(t,A^R,\varepsilon).
\end{aligned}
\end{equation}
Thus, we have shown that there is an $A >2$ such that 
for each $R\geq R_0$
\begin{align}
\label{summable}
\P\left[E_0\left(e^{\overline{\cI}_N^{\kappa}(0,rA^R,0)}
\one{\big\{N(X^{\kappa},rA^R)\leq MrA^R\big\}}\right) 
> C_2^{N}(rA^R,A^R,\varepsilon)\right] \leq e^{-r\eta c_1},
\end{align}
which is summable on $r\in \N$. By the boundedness of $\xi_N$, the same is true
without the indicator in the expectation (after a possible enlargement of $C_2^N$ by $\varepsilon$). Further note that 
\begin{equation}
\label{Cestimate}
\begin{aligned}
&\frac{1}{rA^R}\log C_2^{N}(rA^R,A^R,\varepsilon)\\
&= -2d\kappa + 2d\kappa e^{NC/A^{2Rd}}
+N\eta+ \frac{NC}{A^{2Rd}}
+ \overline{\lambda}_0^{\one}(\kappa)+\varepsilon 
+ \frac{1}{A^{R}}\log\left(e^{-A^{R}(\overline{\lambda}_0^{\one}(\kappa)
+\varepsilon +c-N)}+1\right),
\end{aligned}
\end{equation}
so that $C_2^{N}(rA^{R},A^R,\varepsilon)$ is indeed of the form $\overline{\lambda}_0^{\one}
(\kappa)+\varepsilon$.  Thus, we have proved Lemma~\ref{summablelem} and hence 
Theorem~\ref{nonlocLyp}.
\epr


\subsection{Proof of Lemmas \ref{nonloclem1}--\ref{nonloclem3}}
\label{S4.4}


\subsubsection{Proof of Lemma \ref{nonloclem1}}
\label{S4.4.1}

\bpr
The proof comes in three steps and uses ideas from Cranston, Mountford and Shiga~\cite{CMS02}, 
Lemma 4.3, and Drewitz, G\"artner, Ramirez and Sun~\cite{DGRS12}, Lemma 4.3.

\medskip\noindent
{\bf 1.}
Suppose that we already showed that, $\xi$-a.s.\ and independently of the realization 
of $\xi$, for all $\varepsilon>0$ there is an $\eta'>0$ such that, for all $0<\eta<\eta'$,
\begin{equation}
\label{DGRSest1}
\begin{aligned}
\limsup_{R \to \infty} 
\sup_{\substack{ x,y \in [-2MA^R, 2MA^R]^{d}\cap\Z^d \\ \|x-y\|\leq \eta A^R} }
\frac{1}{A^R} 
&\Big|\log E_x\left(e^{\overline{\cI}_N^{\kappa}(0,A^R,0)}
\one{\Big\{N(X^{\kappa},A^R)\leq MA^R\Big\}}\right)\\
& - \log E_y\left(e^{\overline{\cI}_N^{\kappa}(0,A^R,0)}
\one{\Big\{N(X^{\kappa},A^R)\leq MA^R\Big\}}\right)\Big| \leq \varepsilon.
\end{aligned}
\end{equation}
We show how (\ref{DGRSest1}) can be used to obtain the claim.

\medskip\noindent
{\bf 2.} By Proposition~\ref{shapeThm} for fixed $\tilde{\delta} >0$ there is an $R_0
=R_0(A,\tilde{\delta})\in\N$ such that, for all $R\geq R_0$,
\begin{equation}
\label{goodest1}
\P\Bigg(E_0\left(e^{\overline{\cI}^{\kappa}(0,A^R,0)}
\one{\Big\{N(X^{\kappa},A^R)\leq MA^R\Big\}}\right) 
\leq e^{(\overline{\lambda}_0^{\one}(\kappa)+\varepsilon)A^R}\Bigg) \geq 1-\tilde{\delta}.
\end{equation}
To extend this to $\overline{\cI}_N$, note that for each $t\in\N$
\begin{equation}
\label{truncxi}
\int_0^t\xi_N(X^{\kappa}(s),s)\,ds 
\leq \int_0^t\xi(X^{\kappa}(s),s)\,ds 
+ \int_0^t-\xi(X^{\kappa}(s),s)\one\{\xi(X^{\kappa}(s),s)<-N\}\,ds.
\end{equation}
Thus, given a realization of $X^{\kappa}$ with no more than $Mt$ jumps, by the fact that 
$\xi$ is G\"artner-negative-hyper-mixing, (\ref{appomultiscalelem}) and the arguments 
given in the proof of Lemma~\ref{FKbadblock}, the second term in the right-hand side is 
at most
\begin{equation}
(M+1)tA^dC'\sum_{R=R_0}^{\infty}A^{-Rd}.
\end{equation}
Hence (\ref{goodest1}) remains true when we replace $\overline{\cI}$ by $\overline{\cI}_N$.
According to (\ref{DGRSest1}), this estimate also holds when we replace $0$ by any $x$ with 
$\|x\| \leq \eta A^R$ for $\eta$ small enough, independently of the realization of $\xi$. 
Consequently, for any $\tilde{\delta}>0$ there is an $R_0\in\N$ such that 
\begin{align}
\label{goodest2}
\P\left(\sup_{x\colon\,\|x\|\leq \eta A^R}
E_x\left(e^{\overline{\cI}_N^{\kappa}(0,A^R,0)}
\one{\Big\{N(X^{\kappa},A^R)\leq MA^R\Big\}}\right) 
\leq e^{(\overline{\lambda}_0^{\one}(\kappa)+\varepsilon)A^R}\right) 
\geq 1-\tilde{\delta} \qquad R\geq R_0.
\end{align}
Next, note that $[-2MA^R,2MA^R]^d$ can be divided into $K$ boxes, with $K \sim 4^dM^d/\eta^d$, 
of the form $x_i(A^R)+ B_{\eta A^R}$, where the $x_i(A^R)$'s are separated by $\eta A^R$. By 
the stationarity of $\xi$ in space, we have
\begin{equation}
\label{shiftgood}
\begin{aligned}
&\P\left(E_0\left(e^{\overline{\cI}_N^{\kappa}(0,A^R,0)}
\one{\Big\{N(X^{\kappa},A^R)\leq MA^R\Big\}}\right)
\leq e^{(\overline{\lambda}_0^{\one}(\kappa)+\varepsilon)A^R}\right)\\
&= \P\left(E_{x_i(A^R)}\left(e^{\overline{\cI}_N^{\kappa}(0,A^R,0)}
\one{\Big\{N(X^{\kappa},A^R)\leq MA^R\Big\}}\right)
\leq e^{(\overline{\lambda}_0^{\one}(\kappa)+\varepsilon)A^R}\right).
\end{aligned}
\end{equation}
Thus, for the same choice of $R_0$ as in (\ref{goodest2}),
\begin{equation}
\label{goodest3}
\begin{aligned}
&\P\left(\sup_{y\colon\,y \in x_i(A^R)+B_{\eta A^R}}
E_y\left(e^{\overline{\cI}_N^{\kappa}(0,A^R,0)}
\one{\Big\{N(X^{\kappa},A^R)\leq MA^R\Big\}}\right) 
\leq e^{(\overline{\lambda}_0^{\one}(\kappa)+\varepsilon)A^R}\right)\\ 
&\qquad \geq 1-\tilde{\delta} \qquad R\geq R_0.
\end{aligned}
\end{equation}
Since $K$ is independent of $R_0$, we may conclude that
\begin{equation}
\label{goodfinal}
\P\left(B_{R}^{A}(0,0) \mbox{ is $N$-sufficient}\right)\geq 1-\tilde{\delta}.
\end{equation}
By the stationarity of $\xi$ in space and time, the same statement holds for any block 
$B_R^{A}(x,k)$, which proves the claim. It therefore remains to prove (\ref{DGRSest1}).

\medskip\noindent
{\bf 3.} 
Since for $M$ large the event $\Big\{N(X^{\kappa},A^R)>MA^R\Big\}$ does not contribute 
on an exponential scale, in order to prove \eqref{DGRSest1} it suffices to show that
\begin{equation}
\label{DGRSest2}
\begin{aligned}
\limsup_{R \to \infty} 
\sup_{\substack{ x,y \in [-2MA^R,2MA^R]^{d}\cap\Z^d \\ \|x-y\|\leq \eta A^R} }
\frac{1}{A^R} \Bigg|\log\frac{ E_x\left(e^{\overline{\cI}_N^{\kappa}(0,A^R,0)}\right)} 
 { E_y\left(e^{\overline{\cI}_N^{\kappa}(0,A^R,0)}\right)}\Bigg| \leq \varepsilon.
\end{aligned}
\end{equation}
To that end, we show that we can restrict ourself to contributions coming from random walk 
paths that stay within a certain distance of $[-2M,2M]^d\cap\Z^d$. More precisely, 
$\xi$-a.s.\ there is a box $B_L=[-L,L]^{d}\cap \Z^d$, independent of $0<\eta <1$ and 
containing $[-2M,2M]^d\cap \Z^d$, such that 
\begin{equation}
\label{outsideBR}
\begin{aligned}
\sup_{x\in [-2MA^R,2MA^R]^d\cap\Z^d} 
\sum_{\substack{w \in \Z^d \\ w \notin A^RB_L}} 
&\Bigg[E_x\left(e^{\overline{\cI}_N^{\kappa}(0,\eta A^R,0)} 
\delta_{w}(X^{\kappa}(\eta A^R))\right)
E_w\left(e^{\overline{\cI}_N^{\kappa}(0,(1-\eta)A^R,\eta A^R)}\right)\Bigg]\\
&\leq e^{A^R} \sup_{x \in [-2MA^R,2MA^R]^d\cap \Z^d} 
P_x\left(X^{\kappa}(\eta A^R) \notin A^RB_L\right)\\
&= e^{A^R} P_0\Big(X^{\kappa}(\eta A^R) \notin A^R B_{L-2M}\Big)\\
& \leq e^{A^R} \exp\Bigg\{-A^R(L-2M)\left(
\log\left(\frac{L-2M}{\kappa d}\right)-1\right)\Bigg\},
\end{aligned}
\end{equation}
where the last inequality follows from G\"artner and Molchanov~\cite{GM90}, Lemma 4.3. 
Consequently, we may concentrate in (\ref{DGRSest2}) on the contribution coming from 
paths that stay inside $A^RB_L \cap\Z^d$. Next, note that, $\xi$-a.s.\ and uniformly 
in $x,y \in [-2MA^R,2MA^R]^d\cap\Z^d$,
\begin{equation}
\label{ratio}
\begin{aligned}
&\frac{\sum_{w\in A^RB_L}E_x\left(e^{\overline{\cI}_N^{\kappa}(0,A^R,0)}
\delta_{w}(X^{\kappa}(\eta A^R))\right)}
{\sum_{w\in A^RB_L}E_y\left(e^{\overline{\cI}_N^{\kappa}(0,A^R,0)}
\delta_{w}(X^{\kappa}(\eta A^R))\right)}\\
& \leq \sup_{w \in A^RB_L} 
\frac{E_x\left(e^{\overline{\cI}_N^{\kappa}(0,\eta A^R,0)}
\delta_{w}(X^{\kappa}(\eta A^R))\right)}
{E_y\left(e^{\overline{\cI}_N^{\kappa}(0,\eta A^R,0)}
\delta_{w}(X^{\kappa}(\eta A^R))\right)}\\
& \leq e^{2N\eta A^R} \sup_{w \in A^RB_L} 
\frac{P_0\left(X^{\kappa}(\eta A^R)=w-x\right)}
{P_0\left(X^{\kappa}(\eta A^R)=w-y\right)}, \qquad 0<\eta<1.
\end{aligned}
\end{equation}
To obtain (\ref{DGRSest1}), it remains to estimate the probabilities in the last line. 
This can be done by applying bounds on probabilities for simple random walks (see 
\cite{DGRS12}, Lemma 4.3 for details).
\epr


\subsubsection{Proof of Lemma \ref{nonloclem2}}
\label{S4.4.2}

The only difference with the situation in the proof of Lemma~\ref{multiscalelem} is that 
we replaced the $R$-blocks $B_{R}^{A}(x,k)$ by the $R$-blocks $\hat{B}_{R}^{A}(x,k)$. 
However, this does not affect the proof. Thus, the proof of Lemma \ref{multiscalelem}
yields the claim.


\subsubsection{Proof of Lemma \ref{nonloclem3}}
\label{S4.4.3}

\bpr
The proof comes in two steps and is essentially a copy of the proof of Lemma~\ref{mainlem}. 
Throughout the proof, $x,x'\in\Z^d$, $k,k'\in\N$ and $\delta>0$ is fixed.

\medskip\noindent
{\bf 1.} 
Pick $a_1,a_2\in\N$ according to our main assumption. We say that $(x,k)$ and $(x',k')$ 
are equivalent if and only if
\begin{equation}
\label{equivnonloc}
x\equiv x' \mod a_1 \quad \mbox{ and } \quad k\equiv k' \mod a_2.
\end{equation}
This equivalence relation divides $\Z^d\times\N$ into $a_1^da_2$ equivalence classes. We 
write $\sum_{(x^*,k^*)}$ 
to denote the sum over all equivalence classes. Furthermore, we define
\begin{equation}
\label{defchinonloc}
\begin{aligned}
\hat{\chi}^{A}(x,k) 
= &\one\Big\{\hat{B}_R^{A}(x,k) \mbox{ is good, but has an $N$-insufficient subpedestal}\Big\}.
\end{aligned}
\end{equation}
Henceforth we assume that all blocks under consideration intersect $[-Mt,Mt]\times[0,t]$. 
Then we have
\begin{equation}
\label{firstestpsilemnonloc}
\begin{aligned}
&\P(C(Mt,\eta t/A^R))\\
&\leq \sum_{(x^*,k^*)}
\P\Bigg(\hspace{-.2cm}\begin{array}{ll}&\exists\mbox{a path with no more than $Mt$ jumps that intersects at least}\\
&\eta t/A^Ra_1^{d}a_2
\mbox{ blocks $\hat{B}_{R}^{A}(x,k)$ with $\hat{\chi}^{A}(x,k)=1$, $(x,k)\equiv (x^*,k^*)$}\end{array}\Bigg).
\end{aligned}
\end{equation} 

\medskip\noindent
{\bf 2.}
To proceed we fix $j\leq Mt$ and an equivalence class.
Let $\nu =\lceil \delta^{-1/(1+d)}\rceil$ and define the space block
\begin{equation}
\label{eq:starblock}
\overline{B}_{R}^{A}(x,k)= 
\Bigg(\prod_{j=1}^{d}[\nu(x(j)-1)4MA^{R}, \nu(x(j)+1)4MA^{R})\cap\Z^d\Bigg)
\times [\nu kA^{R},\nu(k+1)A^{R}).
\end{equation}
As in the proof of Lemma 3.6 we see that a path $\Phi$ with $j$ jumps crosses
at most
\begin{equation}
\label{eq:muj}
\mu(j) = 3^d\bigg(\frac{t+j}{\nu A^{R}} +2\bigg)
\end{equation}
blocks $\overline{B}_{R}^{A}(x,k)$. We write
\begin{equation}
\label{eq:union}
\bigcup_{(x_i,k_i)}\overline{B}_{R}^{A}(x_i,k_i)
\quad \mbox{and} \quad \sum_{\overline{B}_{R}^{A}(x_i,k_i)}
\end{equation}
to denote the union over at most $\mu(j)$ blocks $\overline{B}_{R}^{A}(x_i,k_i)$,
$0\leq i\leq \mu(j)-1$, and to denote the sum over all possible sequences of blocks $\overline{B}_{R}^{A}(x_i,k_i)$, that may be crossed by a path $\Phi$ with $j$ jumps, respectively.
To further estimate, define for any sequence $\overline{B}_{R}^{A}(x_0,k_0),
\ldots \overline{B}_{R}^{A}(x_{\mu(j)-1},k_{\mu(j)-1})$ of blocks 
(\ref{eq:union}) and any $n\in\N$ the event
\begin{equation}
\label{eq:partitionevent}
\begin{aligned}
&\hat{\mathcal{A}}^{n}((x_0,k_0),\ldots, (x_{\mu(j)-1},k_{\mu(j)-1}))
= \bigg\{\hspace{-.4cm}\begin{array}{ll}&\mbox{the union in (\ref{eq:union})
contains $n$ blocks $\hat{B}_{R}^{A}(x,k)$}\\
&\mbox{with $\hat{\chi}^{A}(x,k) = 1$,
$(x,k)\equiv (x^*,k^*)$}\end{array}\hspace{-.1cm}\bigg\}.
\end{aligned}
\end{equation}
By our assumption, there is $R_0\in\N$ such that the probability of the event in (\ref{eq:partitionevent}) may be bounded from above by $K\delta^{n}$.
We conclude in the same way as in Lemma 3.6, see also Claim \ref{cl:cardinality}, that for some $C>0$
\begin{equation}
\label{eq:binomialsum}
\begin{aligned}
\sum_{\overline{B}_{R}^{A}(x_i,k_i)}
&\sum_{n=\frac{\eta t}{A^{R}a_1^da_2}}^{L}
\P\Big(\hat{\mathcal{A}}^{n}((x_0,k_0),\ldots, (x_{\mu(j)-1},k_{\mu(j)-1}))\Big)\\
&\qquad\qquad\leq e^{C\mu(j)}(1-\delta)^{-L}K
\P\bigg(T_L\geq \frac{\eta t}{A^{R}a_1^da_2}\bigg),
\end{aligned}
\end{equation}
where $T_L=\mathrm{BIN}(L,\delta)$, $L=\nu^{(1+d)}\mu(j)$.
The same arguments as in Lemma~\ref{mainlem} 
yield that the binomially distributed random variable may be bounded from above by
\begin{equation}
\label{appofBernsteinnonloc}
\begin{aligned}
\exp\left\{-C'\eta t/A^Ra_1^da_2\right\}
\end{aligned}
\end{equation}
for some $C'>0$. 
Similarly as in Lemma \ref{mainlem}, if $\delta^{-1/1+d}\in\N$, 
the second term on the right hand side of (\ref{eq:binomialsum}) may be bounded from above by
\begin{equation}
\label{deltaest}
\frac{3^d(M+1)t\delta^{1/1+d}}{(1-\delta)A^R} + \frac{3^d2}{1-\delta}.
\end{equation}
Since $\delta$ tends to zero, if $R$ tends to infinity,
for $t$ large enough the second term on the right hand side of (\ref{eq:binomialsum}) does not contribute.
The same can be seen to be true for the first term on the right hand side of (\ref{eq:binomialsum}).
Finally, to estimate (\ref{firstestpsilemnonloc}), insert (\ref{appofBernsteinnonloc}) 
into (\ref{firstestpsilemnonloc}), to obtain
\begin{equation}
\label{fourthestpsilemnonloc}
\begin{aligned}
&\P\left(C(Mt,\eta t/A^{R})\right)
\leq Ka_1^da_2Mt \exp\left \{-C'\eta t/A^{R}\right\},
\end{aligned}
\end{equation}
which yields the claim.
\epr


\section{Continuity at $\kappa=0$}
\label{S5}

The proof of Theorem~\ref{continuity} is given in Section~\ref{S5.3}. 
It is based on Lemmas~\ref{contLem1}--\ref{contLem3} below, which are stated in 
Section~\ref{S5.1} and are proved in Section~\ref{S5.2}.


\subsection{Three lemmas}
\label{S5.1}

Fix $b \in (0,1)$, and define the set of paths
\begin{equation}
\label{Akappant}
\begin{aligned}
A^\kappa_{nt} = 
&\Bigg\{\Phi \colon\, [0,nt] \to \Z^d\colon\, N(\Phi,nt) \leq 
\frac{1}{\log(1/\kappa)^b}\,nt,\\ 
&\qquad\forall\,1 \leq j \leq n\,\, \exists\, x_j\in\Z^d\colon\, \|x_j\| 
\leq \frac{1}{\log(1/\kappa)^b}\,nt,\,
\Phi(s) = x_j \,\,\forall\, s \in[(j-1)t +1,jt)\Bigg\},
\end{aligned}
\end{equation}
i.e., paths of length $nt$ that do not jump in time intervals of length $t-1$ 
and whose number of jumps is bounded by $\frac{1}{\log(1/\kappa)^b}\,nt$. 
Note that $\kappa \mapsto A^\kappa_{nt}$ is 
non-decreasing.

\bl{contLem1} 
Suppose that $\xi$ satisfies condition {\rm (b)} in Definition~{\rm \ref{Gartnerposreg}}.
Then, $\xi$-a.s., for any sequence of positive numbers $(a_m)_{m\in\N}$ tending to 
zero there exists a strictly positive and non-increasing sequence $(\kappa_m)_{m\in\N}$ 
such that, for all $m \in \N$ and $0< \kappa \leq \kappa_m$, there exists a $t_m 
= t_m(\kappa_m)$ such that, for all $t \in \Q \cap [t_m,\infty)$, there exists an 
$n_m = n_m(\kappa_m,t)$ such that
\begin{equation}
\label{contLem1res}
\sup_{\Phi \in A^\kappa_{nt}} \sum_{j=1}^n \int_{(j-1)t+1}^{jt} 
\xi(\Phi(j-1)t+1),s)\,ds \leq a_mnt \qquad \forall\,n \geq n_m. 
\end{equation}
\el
\noindent
We say that two paths $\Phi_1$ and $\Phi_2$ on $[0,nt]$ are equivalent, written
$\Phi_1 \sim \Phi_2$, if and only if 
\begin{equation}
\label{Phisim}
\Phi_1|_{[(j-1)t+1,jt)} = \Phi_2|_{[(j-1)t+1,jt)} 
\qquad \forall\, 1 \leq j \leq n.
\end{equation} 
This defines an equivalence relation $\sim$, and we denote by $A^{\kappa,\sim}_{nt}$
the set of corresponding equivalent classes. The following lemma provides an 
estimation of the cardinality of $A^{\kappa,\sim}_{nt}$.

\bl{contLem2} 
$|A^{\kappa,\sim}_{nt}| \leq (nt/\log(1/\kappa)^b)2^n(2d)^{nt/\log(1/\kappa)^b}+1$.
\el

\bl{contLem3}
Suppose that $\xi$ is G\"artner-positive-hyper-mixing.
Then there are $A,C>0$ such that $\xi$-a.s. for $t\in\Q$ large enough
and any choice of disjoint subintervals $\mathcal{I}_1, \mathcal{I}_2,\cdots, \mathcal{I}_k$,
$k\in\N$, of $[0,t]$ such that $|\mathcal{I}_i|=|\mathcal{I}_l|$, $i,l\in\{1,2,\cdots,k\}$,
and each $R_0\in\N$ and each path $\Phi\in B_{t}$
\begin{equation}
\label{contlem3est}
\sum_{i=1}^{k}\int_{\mathcal{I}_i}\xi(\Phi(s),s)\,ds 
\leq k|\mathcal{I}_1|CA^{R_0d}+ (t+N(\Phi,t))C'A^d\sum_{R=R_0}^{\infty}A^{-Rd},
\end{equation}
for some constant $C'>0$.
\el


\subsection{Proof of Lemmas~\ref{contLem1}--\ref{contLem3}}
\label{S5.2}


\subsubsection{Proof of Lemma~\ref{contLem2}}
\label{S5.2.1}

\bpr
Fix an integer $k \leq nt$. We start by estimating the number of possible arrangements 
of the jumps in paths with $k$ jumps. Since we do not distinguish between two paths
that coincide on the intervals $[(j-1)t+1,jt)$, $1\leq j\leq n$, only the last jumps 
before the times $(j-1)t+1$, $1\leq j\leq n$, need to be considered. 

First, the number of arrangements with jumps in $0\leq l \leq k$ different intervals is 
$\binom{n}{l}$. Since the number of different intervals cannot exceed $n$, the number 
of different arrangements is bounded from above by $\sum_{l=1}^n \binom{n}{l} \leq 
2^n$. Next, there are $(2d)^k$ different points in $\Z^d$ that can be visited by a 
path with $k$ jumps. Therefore
\begin{equation}
|A^{\kappa,\sim}_{nt}| 
\leq \sum_{k=1}^{nt/\log(1/\kappa)^b} 2^n (2d)^k + 1 
\leq (nt/\log(1/\kappa)^b) 2^n (2d)^{nt/\log(1/\kappa)^b} + 1,
\end{equation}
which proves the claim.
\epr


\subsubsection{Proof of Lemma~\ref{contLem1}}
\label{S5.2.2}

\bpr
Choose $\kappa_1$ such that 
\begin{equation}
\frac{\log(2d)}{\log(1/\kappa_1)^b} < \delta_2(a_1),
\end{equation} 
and $t_1 = t_1(\kappa_1)$ such that 
\begin{equation}
t_1\left(-\delta_2 + \frac{\log(2d)}{\log(1/\kappa_1)^b}\right) < -\log 2.
\end{equation}
Then, by condition {\rm (b)} in Definition \ref{Gartnerposreg} and Lemma~\ref{contLem2}, 
for all $t \geq t_0 \vee t_1$ we have
\begin{equation}
\begin{aligned}
&\P\left(\sup_{\Phi \in A^{\kappa_1}_{nt}} \sum_{j=1}^n \int_{(j-1)t+1}^{jt}
\xi(\Phi((j-1)t+1),s)\,ds \geq a_1 nt\right)\\
&\qquad \leq \sum_{\Phi \in A^{\kappa_1,\sim}_{nt}} 
\P\left(\sum_{j=1}^n \int_{(j-1)t+1}^{jt} \xi(\Phi((j-1)t+1),s)\,ds 
\geq a_1 nt\right)\\
&\qquad \leq\left[(nt/\log(1/\kappa)^b) 2^n (2d)^{nt/\log(1/\kappa)^b} + 1\right] 
e^{-\delta_2 nt},
\end{aligned}
\end{equation}
which is summable on $n$. Hence, by the Borel-Cantelli lemma, there exists a set 
$B_{\kappa_1,t}$ with $\P(B_{\kappa_1,t})=1$ for which there exists an $n_0= n_0
(\xi,\kappa_1,t)$ such that
\begin{equation}
\label{Akappa1bd}
\sup_{\Phi \in A^{\kappa_1}_{nt}} \sum_{j=1}^n 
\int_{(j-1)t+1}^{jt} \xi(\Phi((j-1)t +1),s)\,ds 
\leq a_1 nt \qquad \forall\, n \geq n_0.
\end{equation}
Since $\kappa \mapsto A^\kappa_{nt}$ is non-decreasing, (\ref{Akappa1bd}) is true 
for all $0<\kappa \leq \kappa_1$. Define 
\begin{equation}
B_1 = \bigcap_{t \geq t_1, t \in \Q} B_{\kappa_1,t},
\end{equation}
for which still $\P(B_1)=1$. Similarly, we can construct sets $B_m$, $m\in\N\setminus\{1\}$, 
with $\P(B_m)=1$ such that on $B_m$ there exist $\kappa_m$ and $t_m = t_m(\kappa_m)$ such 
that for all $0<\kappa \leq \kappa_m$ and $t \geq t_m$ with $t \in \Q$ there exists 
an $n_0 = n_0(\xi,\kappa_m,t)$ such that
\begin{equation} 
\sup_{\Phi \in A^\kappa_{nt}} \sum_{j=1}^n \int_{(j-1)t+1}^{jt}
\xi(\Phi((j-1)t+1),s)\,ds 
\leq a_m nt \qquad \forall\, n \geq n_0.
\end{equation}
Hence $B= \cap_{m\in\N} B_{m}$ is the desired set. Note that we can control the value of 
$t_m$ by choosing $\kappa_m$ small enough. Indeed, with the right choice of $\kappa_m$, 
it follows that $t_{m-1}(\kappa_{m-1}) = t_m(\kappa_m)$ for all $m\in\N$.
\epr


\subsubsection{Proof of Lemma~\ref{contLem3}}
\label{S5.2.3}

\bpr
Fix $A,C$ as in Section \ref{S3}, $R_0\in\N$, a path $\Phi\in B_{t}$ and disjoint subintervals  $\mathcal{I}_1, \mathcal{I}_2,\cdots, \mathcal{I}_k$, $k\in\N$, of $[0,t]$ with equal length.
Note that
\begin{equation}
\label{disjint}
\begin{aligned}
\sum_{j=1}^{k}\int_{\mathcal{I}_j}\xi(\Phi(s),s)\, ds
&\leq \sum_{j=1}^{k}\int_{\mathcal{I}_j} \xi(\Phi(s),s)\one{\{\xi(\Phi(s),s) \leq CA^{R_0d}\}}\,ds\\
&\qquad + \int_{0}^{t} \xi(\Phi(s),s)\one{\{\xi(\Phi(s),s) > CA^{R_0d}\}}\,ds.
\end{aligned}
\end{equation}
By (\ref{appomultiscalelem}) and Lemma \ref{FKbadblock},
for $t\in \Q$ sufficiently large,
the second term on the right hand side
in (\ref{disjint}) may be bounded from above by
\begin{equation}
\label{disjintup}
(t+N(\Phi,t))C'A^d\sum_{R=R_0}^{\infty}A^{-Rd}.
\end{equation}
Inserting (\ref{disjintup}) into (\ref{disjint}) yields the claim.
\epr



\subsection{Proof of Theorem~\ref{continuity}}
\label{S5.3}

In this section we prove Theorem~\ref{continuity} with the help of
Lemmas~\ref{contLem1}--\ref{contLem3}. The proof comes in three steps, 
organized as Sections~\ref{S5.3.1}--\ref{S5.3.3}.


\subsubsection{Estimation of the Feynman-Kac representation on $A_{nT}^\kappa$}
\label{S5.3.1}

Consider the case $u_0(x)=\delta_0(x)$, $x\in\Z^d$. Recall (\ref{Lyapexploc1}) and 
(\ref{FK}), and estimate
\begin{equation}
\label{Lyapupbd}
\lambda_0^{\delta_0}(\kappa) \leq \lim_{n \to \infty} \frac{1}{nT} \log 
E_0\left(e^{\overline{\cI}{\kappa}(0,nT,0)}\right)<\infty,
\qquad T>0,
\end{equation}
where we reverse time, use that $X^\kappa$ is a reversible dynamics, and remove 
the constraint $X^\kappa(nT)=0$. Recalling (\ref{Akappant}) and (\ref{Phisim}),
we have
\begin{equation}
\begin{aligned}
E_0\left(e^{\overline{\cI}{\kappa}(0,nT,0)} 
\one_{A^\kappa_{nT}}\right)
&= \sum_{\Phi \in A^{\kappa,\sim}_{nT}} 
E_0\left(e^{\overline{\cI}{\kappa}(0,nT,0)}
\one_{A^\kappa_{nT}} \one_{\big\{X|_{[0,nT]} \sim \Phi\big\}}\right)\\
& = \sum_{\Phi \in A^{\kappa,\sim}_{nT}} 
E_0\left(\exp\left\{\sum_{j=1}^n \overline{\cI}^{\kappa}((j-1)T+1,jT,0)\right\}\right.\\
& \qquad \left. \times
\exp\left\{\sum_{j=0}^{n-1} \overline{\cI}^{\kappa}(jT,jT+1,0)\right\} 
\one_{A^\kappa_{nT}} \one_{\big\{X|_{[0,nT]} \sim \Phi\big\}}\right). 
\end{aligned}
\end{equation}
By the H\"older inequality with $p,q> 1$ and $1/p+1/q=1$, we have
\begin{equation}
\label{cont3}
\begin{aligned}
&E_0\left(e^{\overline{\cI}{\kappa}(0,nT,0)}
\one_{\{A^\kappa_{nT}\}}\right)\\
&\qquad\leq \sum_{\Phi \in A^{\kappa,\sim}_{nT}} 
E_0\left(\exp\left\{p\sum_{j=1}^n\overline{\cI}^{\kappa}((j-1)T+1,jT,0)\,ds\right\}
\one_{A^\kappa_{nT}} \one_{\big\{X|_{[0,nT]} \sim \Phi\big\}}\right)^{1/p}\\
&\qquad\qquad\times 
E_0\left(\exp\left\{q\sum_{j=0}^{n-1} \overline{\cI}^{\kappa}(jT,jT+1,0)
\right\}\one_{A^\kappa_{nT}}\right)^{1/q}.
\end{aligned}
\end{equation}
Next, fix $(a_m)_{m \in \N}$, $(\kappa_m)_{m \in \N}$, $t_m$ as in 
Lemma~\ref{contLem1}, choose $T>0$ such that
\begin{equation}
\label{Tkappa}
t_m \leq T = T(\kappa_m) = K \lfloor{\log(1/\kappa_{m})} \rfloor,
\qquad m \gg 1,
\end{equation}
where $K$ is a constant to be chosen later. 
For all $0<\kappa\leq\kappa_m$ and $n \geq n_m(\kappa_m,T(\kappa))$, 
by Lemma~\ref{contLem1} we have
\begin{equation}
\label{cont5}
\begin{aligned}
&\sum_{\Phi \in A^{\kappa,\sim}_{nT}} 
E_0\left(\exp\left\{p\sum_{j=1}^n \overline{\cI}^{\kappa}((j-1)T+1,jT,0)\right\}
\one_{A^\kappa_{nT}} \one_{\big\{X|_{[0,nT]} \sim \Phi\big\}}\right)^{1/p}\\
&\qquad\leq \sum_{\Phi \in A^{\kappa,\sim}_{nT}} e^{a_m nT}
P_0\left(A_{nT}^\kappa,X|_{[0,nT]} \sim \Phi\right)^{1/p}
\leq e^{a_m nT} |A^{\kappa, \sim}_{nT}|,
\end{aligned}
\end{equation} 
while by Lemma~\ref{contLem3} we have
\begin{equation}
\label{cont7}
\begin{aligned}
&E_0\left(\exp\left\{q\sum_{j=0}^{n-1}\overline{\cI}^{\kappa}(jT,jT+1,0)
\right\}\one_{A^\kappa_{nT}}\right)^{1/q}\\
&\qquad\leq 
\exp\left\{nCA^{R_0d}+ nT\left(1+\frac{1}{\log(1/\kappa)^b}\right)C'A^d\sum_{R=R_0}^{\infty}A^{-Rd}\right\}.
\end{aligned}
\end{equation}
From Lemma~\ref{contLem2} we know that $\frac{1}{nT}\log|A^{\kappa,\sim}_{nT}|$ 
tends to zero if we let first $n\to\infty$ and then $\kappa \da 0$. Therefore, 
combining (\ref{cont3}--\ref{cont7}) and using that $\lim_{\kappa\da 0} T
= \lim_{\kappa\da 0} T(\kappa) = 0$, we get
\begin{equation}
\label{endAkappa}
\limsup_{\kappa\da 0} \limsup_{n\to\infty} \frac{1}{nT} \log 
E_0\left(e^{\overline{\cI}{\kappa}(0,nT,0)}
\one_{A^\kappa_{nT}}\right) \leq \max\left\{a_m,C'A^d\sum_{R=R_0}^{\infty}A^{-Rd}\right\} .
\end{equation}


\subsubsection{Estimation of the Feynman-Kac representation on $[A_{nT}^\kappa]^\cc$} 
\label{S5.3.2}

The proof comes in three steps.

\medskip\noindent
{\bf 1.} 
We start by estimating the corresponding Feynman-Kac term on $[A^\kappa_{nT}]^\cc$. 
Split
\begin{equation}
[A^\kappa_{nT}]^\cc = B^\kappa_{nT} \cup C^\kappa_{nT}
\end{equation}
with
\begin{equation}
B^\kappa_{nT} = \left\{N(X^\kappa,nT) > \frac{1}{\log(1/\kappa)^b}\,nT\right\} 
\end{equation}
and
\begin{equation}
\label{CkappanTdef}
\begin{aligned}
C^\kappa_{nT} &= \Big\{\exists\,1 \leq j\leq n\colon\,  
\mbox{ for all } x\in\Z^d \mbox{ with } \|x\| \leq \frac{1}{\log(1/\kappa)^b}\,nT\\
&\qquad\qquad \mbox{ there exists a } s_j \in [(j-1)T+1,jT) \mbox{ such that } 
X^\kappa(s_j) \neq x \Big\}.
\end{aligned}
\end{equation}
Then
\begin{equation}
P_0\big(B^\kappa_{nT}\big) 
\leq \exp\left\{\big[-J_\kappa(1/\log(1/\kappa)^b) + o_n(1)\big]\,nT\right\},
\end{equation}
where 
\begin{equation}
J_\kappa(x) = x \log(x/2d\kappa) -x + 2d\kappa
\end{equation}
is the large deviation rate function of the rate-$2d\kappa$ Poisson process. Thus, 
by the H\"older inequality with $p,q>1$ and $1/p+1/q=1$, we have
\begin{equation}
\begin{aligned}
&E_0\left(e^{\overline{\cI}{\kappa}(0,nT,0)}
\one_{B^\kappa_{nT}}\right)\\
&\qquad \leq E_0\left(\exp\left\{p\overline{\cI}{\kappa}(0,nT,0)\right\}\right)^{1/p} 
\exp\left\{1/q\big[-J_\kappa(1/\log(1/\kappa)^b) + o_n(1)\big]\,nT\right\}.
\end{aligned}
\end{equation}
Recalling Theorem \ref{expgrowthres} and using that $\lim_{\kappa \da 0}J_\kappa(1/\log(1/\kappa)^b) 
= \infty$, we get
\begin{equation}
\lim_{\kappa\da 0}\lim_{n\to\infty} \frac{1}{nT} \log
E_0\left(e^{\overline{\cI}{\kappa}(0,nT,0)}
\one_{B^\kappa_{nT}}\right) = -\infty.
\end{equation}

\medskip\noindent
{\bf 2.} 
Note that
\begin{equation}
C^\kappa_{nT} \subseteq B^\kappa_{nT} \cup D^\kappa_{nT}
\quad \mbox{with} \quad
D^\kappa_{nT} = \Big(C^\kappa_{nT} \cap [B^\kappa_{nT}]^\cc\Big).
\end{equation}
Since we have just proved that the Feynman-Kac representation on $B^\kappa_{nT}$ 
is not contributing, we only have to look at the contribution coming from 
$D^\kappa_{nT}$, namely,
\begin{equation}
E_0\left(e^{\overline{\cI}{\kappa}(0,nT,0)}
\one_{D^\kappa_{nT}}\right).
\end{equation}
On the event $D_{nT}^{\kappa}$, the random walk $X^\kappa$ stays inside 
the box of radius $nT/\log(1/\kappa)^b$, and jumps during the time intervals 
$[(j-1)T+1,jT)$, $1\leq j\leq n$, defined in \eqref{CkappanTdef}. By the H\"older
inequality with $p,q>1$ and $1/p+1/q=1$, we have
\begin{equation}
E_0\left(e^{\overline{\cI}{\kappa}(0,nT,0)}
\one_{D^\kappa_{nT}}\right) \leq I\times II,
\end{equation}
where
\begin{equation}
\begin{aligned}
I &=
\left[E_0\left(\exp\left\{p\sum_{j=1}^n 
\overline{\cI}{\kappa}((j-1)T+1,jT,0)\right\}
\one_{D^\kappa_{nT}}\right)\right]^{1/p},\\
II &=
\left[E_0\left(\exp\left\{q\sum_{j=0}^{n-1} 
\overline{\cI}^{\kappa}(jT,jT+1,0)\right\}\one_{D^\kappa_{nT}}\right)\right]^{1/q}.
\end{aligned}
\end{equation}
Define
\begin{equation}
\label{Jdef}
J = \Big\{1\leq j\leq n\colon\,N(X^\kappa,[(j-1)T+1,jT))\geq 1\Big\},
\end{equation}
where $N(X^\kappa,\cI)$ is the number of jumps of the random walk $X^\kappa$
during the time interval $\cI$. Using that $X^\kappa$ is not jumping in the time 
intervals $[(j-1)T+1,jT))$, $j\in J^\cc$, we may write
\begin{equation}
\label{xisplit}
\begin{aligned}
&\sum_{j=1}^n \overline{\cI}^{\kappa}((j-1)T+1,jT,0)\\
&\qquad = \sum_{j\in J} \overline{\cI}^{\kappa}((j-1)T+1,jT,0) 
+ \sum_{j\in J^\cc} \int_{(j-1)T+1}^{jT}\xi(X^{\kappa}((j-1)T+1),s)\,ds.
\end{aligned}
\end{equation}

\medskip\noindent
{\bf 3.} 
To estimate the second term in the right-hand side of (\ref{xisplit}), pick any 
$\Phi^{X^\kappa}\in A_{nT}^\kappa$ such that $\Phi^{X^\kappa}= X^\kappa$
on $\cup_{j\in J^\cc} [(j-1)T+1,jT)$ and apply Lemma~\ref{contLem1}, to get
\begin{equation}
\label{Jcpart}
\begin{aligned}
&\sum_{j\in J^\cc} \int_{(j-1)T+1}^{jT}\xi(X^{\kappa}((j-1)T+1),s)\,ds\\
&\qquad \leq a_m nT 
- \sum_{j\in J}\int_{(j-1)T+1}^{jT}\xi(\Phi^{X^\kappa}(s),s)\, ds \quad \xi\mbox{-a.s.}
\end{aligned}
\end{equation}
Note that $\xi$ is G\"artner-negative-hyper-mixing, 
so that by Lemma \ref{contLem3} we
may estimate the second term on the right hand side of (\ref{Jcpart})
by 
\begin{equation}
\label{contlem3appl}
|J|(T-1)CA^{R_0d} + nT\left(1+\frac{1}{\log(1/\kappa)^b}\right)C'A^{d}\sum_{R=R_0}^{\infty}A^{-Rd}.
\end{equation}

To estimate the first term in the right-hand side of (\ref{xisplit}), apply
Lemma~\ref{contLem3}, to get
\begin{equation}
\label{pfoufff}
\begin{aligned}
&\sum_{k=1}^n\left[ E_0\left(\exp\left\{p\sum_{j\in J}
\overline{\cI}^{\kappa}((j-1)T+1,jT,0)\right\}
\one_{D^\kappa_{nT}}\one_{\{|J|=k\}}\right) e^{pk(T-1)CA^{R_0d}}\right]\\
&\qquad \leq \exp\left\{pnT\left(1+\frac{1}{\log(1/\kappa)^b}\right)C'A^{d}\sum_{R=R_0}^{\infty}A^{-Rd}\right\}\\
&\qquad\qquad\times\sum_{k=1}^n \exp\left\{2pk(T-1)CA^{R_0d}\right\}
P_0\left(D_{nT}^{\kappa},|J|=k\right).
\end{aligned}
\end{equation}
The distribution of $|J|$ is BIN$(n,1- e^{-2d\kappa (T-1)})$. Hence
the sum on the right hand side of (\ref{pfoufff}) is bounded from above by
\begin{equation}
\label{BinJpart}
\begin{aligned}
& \sum_{k=1}^n \binom{n}{k} \left(1-e^{-2d\kappa(T-1)}\right)^k 
e^{-2d\kappa(T-1)(n-k)} e^{2pk(T-1)CA^{R_0d}}\\
&\leq \left((1-e^{-2d\kappa(T-1)}) e^{2p(T-1)CA^{R_0d}} 
+ e^{-2d\kappa (T-1)}\right)^n.
\end{aligned}
\end{equation}
Combining (\ref{Jcpart}--(\ref{BinJpart})), we arrive at
\begin{equation}
\label{Ibd}
\begin{aligned}
I\leq
&e^{a_mnT}\exp\left\{2nT\left(1+\frac{1}{\log(1/\kappa)^b}\right)
C'A^{d}\sum_{R=R_0}A^{-Rd}\right\}\\
&\qquad\qquad\times \left((1-e^{-2d\kappa(T-1)}) e^{2p(T-1)CA^{R_0d}} 
+ e^{-2d\kappa (T-1)}\right)^{n/p}.
\end{aligned}
\end{equation}
On the other hand, by Lemma \ref{contLem3}, we have
\begin{equation}
\label{IIbd}
II \leq e^{nCA^{R_0d}}\exp\left\{nT\left(1+\frac{1}{\log(1/\kappa)^b}\right)
C'A^{d}\sum_{R=R_0}A^{-Rd}\right\}.
\end{equation}
Therefore, combining (\ref{Ibd}--\ref{IIbd}), we finally obtain
\begin{equation}
\label{DnTkappabd}
\begin{aligned}
&E_0\left(e^{\overline{\cI}^{\kappa}(0,nT,0)}
\one_{D^\kappa_{nT}}\right)\\
&\qquad \leq e^{a_m nT}e^{nCA^{R_0d}}\exp\left\{3nT\left(1+\frac{1}{\log(1/\kappa)^b}\right)
C'A^{d}\sum_{R=R_0}A^{-Rd}\right\}\\
&\qquad\qquad\times\left((1-e^{-2d\kappa(T-1)}) e^{2p(T-1)CA^{R_0d}} + e^{-2d\kappa(T-1)}\right)^{n/p}.
\end{aligned}
\end{equation}


\subsubsection{Final estimation}
\label{S5.3.3}

By (\ref{DnTkappabd}), we have
\begin{equation}
\label{DnTkappabd2}
\begin{aligned}
&\frac{1}{nT} \log E_0\left(e^{\overline{\cI}^{\kappa}(0,nT,0)}
\one_{D^\kappa_{nT}}\right)\\
&\qquad \leq a_m  + \frac{2d\kappa(T-1)}{pT}\,e^{2p(T-1)CA^{R_0d}} 
 + \frac{CA^{R_0d}}{T}+ 3\left(1+\frac{1}{\log(1/\kappa)^b}\right)C'A^d\sum_{R=R_0}A^{-Rd}.
\end{aligned}
\end{equation}
Abbreviate $M_2 = 2pCA^{R_0d}$ and recall (\ref{Tkappa}). Then the right-hand
side of (\ref{DnTkappabd2}) is asymptotically equivalent to
\begin{equation}
a_m + \frac{2d\kappa}{p} (1/\kappa)^{M_2K} + 3\left(1+\frac{1}{\log(1/\kappa)^b}\right)C'A^d\sum_{R=R_0}A^{-Rd}, \qquad \kappa\da 0.
\end{equation}
Choosing $K \leq 1/2M_2$, $K \in \Q$, and recalling (\ref{endAkappa}) 
we finally arrive at
\begin{equation}
\begin{aligned}
\lambda_{0}^{\delta_0}(\kappa)
\leq a_m+  \frac{2d\sqrt{\kappa}}{p}+ 3\left(1+\frac{1}{\log(1/\kappa)^b}\right)C'A^d\sum_{R=R_0}A^{-Rd},
\end{aligned}
\end{equation}
which tends to zero as $\kappa\da 0$, $R_0\to\infty$ and $m\to\infty$.


\section{No Lipschitz continuity at $\kappa=0$}
\label{S6}

In this section we prove Theorem~\ref{nonLipsch}. The proof is very close to that of 
G\"artner, den Hollander and Maillard~\cite{GdHM11}, Theorem 1.2(iii), where it is 
assumed that $\xi$ is bounded from below. For completeness we will repeat the main 
steps in that proof. 

\bpr
Fix $C_1>0$, write
(see \ref{PiL}),
\begin{equation}
\lambda_{0}^{\delta_0}(\kappa)
= \lim_{n \to \infty} \frac{1}{nT+1} \log 
E_0\left(e^{\overline{\cI}^{\kappa}(0,nT+1,0)}
\delta_0(X^{\kappa}(nT+1))\one{\Big\{X^{\kappa}([0,nT+1])\subseteq[C_1]_{nT+1}\Big\}}\right)
\end{equation}
and abbreviate
\begin{equation}
I_j^\xi(x) = \int_{(j-1)T+1}^{jT} \xi(x,s)\,ds, \qquad
Z_j^\xi = \argmax_{x \in \{0,e\}} I_j^\xi(x), \qquad 1\leq j\leq n.
\end{equation}
Consider the event
\begin{equation}
A^\xi 
= \left[\bigcap_{j=1}^n \{X^\kappa(t) = Z_j^\xi\,\,\forall\,t \in[(j-1)T+1,jT)\}\right]
\cap \left\{X^\kappa(nT+1)=0\right\}.
\end{equation} 
We have
\begin{equation}
\label{nonlip3}
\begin{aligned}
&E_0\left(\exp\left\{\overline{\cI}^{\kappa}(0,nT+1,0)\right\} 
\delta_0(X^\kappa(nT+1))\right)\\
&\qquad\geq E_0\left(\exp\left\{\sum_{j=1}^{n+1}\overline{\cI}^{\kappa}((j-1)T,(j-1)T+1,0)
+ \sum_{j=1}^n \overline{\cI}^{\kappa}((j-1)T+1,jT,0)\right\}\,\one_{A^\xi}\right).
\end{aligned}
\end{equation}
Using the reverse H\"older inequality with $q<0<p<1$ and $1/q+1/p=1$, we have
\begin{equation}
\label{genlowerbound}
\begin{aligned}
&E_0\left(\exp\left\{\sum_{j=1}^{n+1} \overline{\cI}^{\kappa}((j-1)T,(j-1)T+1,0) 
+ \sum_{j=1}^n \overline{\cI}^{\kappa}((j-1)T+1,jT,0)\right\}\,\one_{A^\xi} \right)\\
&\qquad\geq\left[E_0
\left(\exp\left\{q\sum_{j=1}^{n+1}
 \overline{\cI}^{\kappa}((j-1)T,(j-1)T+1,0)\right\}
 \one{\Big\{X^{\kappa}([0,nT+1])\subseteq[C_1]_{nT+1}\Big\}}\right)\right]^{1/q}\\
&\qquad\qquad\times
\left[E_0\left(\exp\left\{p\sum_{j=1}^n \overline{\cI}^{\kappa}((j-1)T+1,jT,0)\right\}\,
\one_{A^\xi}
\one{\Big\{X^{\kappa}([0,nT+1])\subseteq[C_1]_{nT+1}\Big\}}\right)\right]^{1/p}.
\end{aligned}
\end{equation}
To estimate the first term on the right hand side of (\ref{genlowerbound}),
fix $R_0\in\N$ and choose $A,C>0$ such that all results of Section
\ref{S3} are satisfied for $q\xi$.
Moreover, note that by a refinement of the arguments given in the proof of Lemma
\ref{FKbadblock} and (\ref{appomultiscalelem}),
one has $\xi$-a.s. for $nT+1\in\N$ sufficiently large  
\begin{equation}
\label{appmultinonLip}
\begin{aligned}
&\sum_{j=1}^{n+1}\int_{(j-1)T}^{(j-1)T+1}q\xi(X^{\kappa}(s),s)
\one{\left\{q\xi(X^{\kappa}(s),s)>-qCA^{R_0d}\right\}}\, ds\\
&\qquad\leq -\frac{nT+1+N(X^{\kappa},nT+1)}{T}qC'A^{d}\sum_{R=R_0}^{\infty}A^{-Rd}.
\end{aligned}
\end{equation}
Consequently,
\begin{equation}
\label{nonLipshort1}
\begin{aligned}
&E_0
\left(\exp\left\{q\sum_{j=1}^{n+1}
 \overline{\cI}^{\kappa}((j-1)T,(j-1)T+1,0)\right\}
 \one{\Big\{X^{\kappa}([0,nT+1])\subseteq[C_1]_{nT+1}\Big\}}\right)\\
 &\qquad\leq e^{-q(n+1)CA^{R_0d}}
 E_0\left(\exp\left\{-\frac{nT+1+N(X^{\kappa},nT+1)}{T}qC'A^{d}\sum_{R=R_0}^{\infty}A^{-Rd}\right\}\right),
\end{aligned}
\end{equation} 
which equals
\begin{equation}
\label{nonLipshort2}
\begin{aligned}
&e^{-q(n+1)CA^{R_0d}}\exp\left\{-\frac{nT+1}{T}qC'A^{d}\sum_{R=R_0}^{\infty}A^{-Rd}\right\}\\
&\qquad\times\exp\left\{2d\kappa(nT+1)\Big(e^{-\frac{q}{T}C'A^d\sum_{R=R_0}^{\infty}A^{-Rd}}-1\Big)\right\}.
\end{aligned}
\end{equation}
As in the proof of \cite{GdHM11}, Theorem 1.2(iii), we have
\begin{equation}
\label{nonlip5}
\begin{aligned}
&\left[E_0\left(\exp\left\{p\sum_{j=1}^n \overline{\cI}^{\kappa}((j-1)T+1,jT,0)\right\}\,
\one_{A^\xi}\right)\right]^{1/p}\\
&\qquad \geq \left[\exp\left\{[1+o_{n}(1)]\,np\,
\E\left(\max\{I_1^\xi(0),I_1^\xi(e)\}\right)\right\}
[p_1^\kappa(e)]^{n+1}\,e^{-2d\kappa n(T-1)}\right]^{1/p}.
\end{aligned}
\end{equation}
Combining (\ref{nonlip3}--\ref{nonlip5}), we arrive at
\begin{equation}
\begin{aligned}
&\frac{1}{nT+1} \log 
E_0\left(e^{\overline{\cI}^{\kappa}(0,nT+1,0)}\, 
\delta_0(X^{\kappa}(nT+1))\right)\\
&\qquad \geq -\frac{1}{(nT+1)}\,CA^{R_0d}(n+1)-\frac{1}{T}C'A^d\sum_{R=R_0}^{\infty}A^{-Rd}
+\frac{2d\kappa}{q}\Big(e^{-\frac{q}{T}C'A^d\sum_{R=R_0}^{\infty}A^{-Rd}}-1\Big)\\
&\qquad\qquad + \frac{1}{p(nT+1)}
\left[[1+o_n(1)]\,pn\,\E\left(\max\left\{I_1^\xi(0),I_1^\xi(e)\right\}\right)\right]
 -\frac{2d\kappa n(T-1)}{p(nT+1)}\\
 &\qquad\qquad + \frac{n+1}{p(nT+1)}\,
\log p_{1}^{\kappa}(e).
\end{aligned}
\end{equation}
Using that $p_1^\kappa(e)=\kappa[1+ o_\kappa(1)]$ as $\kappa \da 0$ and letting 
$n\to\infty$, we get that 
\begin{equation}
\begin{aligned}
&\lambda_0^{\delta_0}(\kappa) \geq -\frac{1}{(nT+1)}\,CA^{R_0d}(n+1)-\frac{1}{T}C'A^d\sum_{R=R_0}^{\infty}A^{-Rd}
+\frac{2d\kappa}{q}\Big(e^{-\frac{q}{T}C'A^d\sum_{R=R_0}^{\infty}A^{-Rd}}-1\Big)\\
&\qquad\qquad - \frac{2d\kappa(T-1)}{pT} 
+ \frac{1}{T}\,[1+o_{\kappa}(1)]\,\left(\tfrac12\,\E(T-1) - \frac{1}{p}\,
\log(1/\kappa)\right).
\end{aligned}
\end{equation}
At this point we can copy the rest of the proof of \cite{GdHM11}, Theorem 
1.2(iii), with a few minor adaptations of constants.
\epr


\section{Examples}
\label{S7}

In Section~\ref{S7.1} we prove Corollary~\ref{mixingex}, in Section~\ref{S7.2} we prove
Corollary~\ref{contex}.


\subsection{Proof of Corollary \ref{mixingex}}
\label{S7.1}

In Section~\ref{S7.1.1} we settle Part (1), in Section~\ref{S7.1.2} we settle Part (2). 


\subsubsection{Proof of Corollary \ref{mixingex} {\rm (1)}}
\label{S7.1.1}

{\bf 1.1} 
The first condition in Definition~\ref{Gartnerhypmix} is satisfied by our assumption 
on $\xi$.

\medskip\noindent
{\bf 1.2} 
We show that $\xi$ is type-I G\"artner-mixing. Fix $A>1$, pick $b=c=0$ and $a_1=a_2=2$ 
(see (\ref{equiv})), and define
\begin{equation}
\label{subpedMarkex}
B_{R}^{A,\mathrm{sub}}(x,k)= 
\Bigg(\prod_{j=1}^{d}[(x(j)-1)A^{R},(x(j)+1)A^{R})\cap\Z^d\Bigg)
\times \{kA^{R}\}.
\end{equation}
We start by estimating the probability of the event
\begin{equation}
\label{goodMarkex}
\mathrm{B}(x,k)\stackrel{\mathrm{def}}=
\big\{B_{R+1}^{A}(x,k) \mbox{ is good, but contains a bad $R$-block}\big\}.
\end{equation}
Note that each $R+1$-block contains at most $2^dA^{(1+d)}$ $R$-blocks. For each such $R$-block 
$B_{R}^{A}(y,l)$ there are no more than $A^{Rd}$ blocks
\begin{equation}
\label{QblockMarkex}
Q_{R}^{A}(z_1)=\prod_{j=1}^{d}[(z_1(j),z_1(j)+A^{R})
\end{equation} 
contained in it. For any such block we may estimate, for $C_1 >0$,
\begin{equation}
\label{Markest1ex}
\begin{aligned}
&\P\Bigg(\exists s \in [lA^R,(l+1)A^R)\colon\, 
\sum_{z_2\in Q_{R}^{A}(z_1)} \xi(z_2,s) > C_1 A^{Rd} \Bigg)\\
&\leq \sum_{k=0}^{\lfloor A^R \rfloor}
\P\Bigg( \sum_{z_2\in Q_{R}^{A}(z_1)}\sup_{s\in [k,k+1)} X_{s}(z_2)>C_{1}A^{Rd}\Bigg)\\
&\leq (A^{R}+1)\exp\left\{-C_{1}A^{Rd}\right\} 
\E\left(\exp\left\{\sup_{s\in[0,1)}X_{s}(0)\right\}\right)^{A^{Rd}},
\end{aligned}
\end{equation}
where we use the time stationarity in the first inequality, and the time stationarity 
and the space independence in the second inequality. Thus, for $C_1$ sufficiently large, there is a 
$C_1' >0$ such that
\begin{equation}
\label{Markest2ex}
\P\big(\mathrm{B}(x,k)\big)
\leq e^{-C_1'A^{Rd}}.
\end{equation}
Moreover, for space-time blocks that are disjoint in space, the corresponding events in 
(\ref{goodMarkex}) are independent. Hence we may assume that $B_{R+1}^{A}(x_1,k_1),\ldots,
B_{R+1}^{A}(x_n,k_n)$ are equal in space but disjoint in time. Since 
\begin{equation}
\label{condMarkex}
\P\Bigg(\bigcap_{i=1}^{n} \mathrm{B}(x_i,k_i)\Bigg)
= \P\Bigg(\mathrm{B}(x_n,k_n) ~\Big|~ \bigcap_{i=1}^{n-1} \mathrm{B}(x_i,k_i)\Bigg)
\P\Bigg(\bigcap_{i=1}^{n-1} \mathrm{B}(x_i,k_i)\Bigg),
\end{equation}
it is enough to show that there is a constant $K<\infty$, independent of $R$, such that the conditional 
probability in (\ref{condMarkex}) may be estimated from above by $K\P(\mathrm{B}(x_n,k_n))$.
To do this, we apply the Markov property to obtain
\begin{equation}
\label{Markest3ex}
\begin{aligned}
\P\Bigg(\mathrm{B}(x_n,k_n)~\Big|~ \bigcap_{i=1}^{n-1} \mathrm{B}(x_i,k_i)\Bigg)
&\leq \P\Bigg(\mathrm{B}(x_n,k_n) ~\Big|~ \bigcap_{i=1}^{n-1} \mathrm{B}(x_i,k_i),
B_{R+1}^{A,\mathrm{sub}}(x_n,k_n) \mbox{ is good}\Bigg)\\
&= \P\Bigg(\mathrm{B}(x_n,k_n)~\Big|~ B_{R+1}^{A,\mathrm{sub}}(x_n,k_n) 
\mbox{ is good}\Bigg).
\end{aligned}
\end{equation}
Thus, the left-hand side of (\ref{Markest3ex}) is at most
\begin{equation}
\label{Markupboundex}
\frac{\P(\mathrm{B}(x_n,k_n))}
{\P(B_{R+1}^{A,\mathrm{sub}}(x_n,k_n) \mbox{ is good})}.
\end{equation}
Since $\lim_{R\to \infty} \P(B_{R+1}^{A,\mathrm{sub}}(x_n,k_n) \mbox{ is good}) = 1$,
we obtain that $\xi$ is type-I G\"artner-mixing.

\medskip\noindent
{\bf 1.3}
Condition {\rm (a3)} in Definition~\ref{Gartnerhypmix} follows from the calculations in 
(\ref{Markest1ex}).

\medskip\noindent
{\bf 2.} 
The same strategy as above works to show that $\xi$ is type-II G\"artner-mixing. If $X$ 
has exponential moments of all negative orders, then the same calculations as in the first 
part show that $\xi$ is G\"artner-negative-hyper-mixing. All requirements of 
Theorems~\ref{expgrowthres}--\ref{nonlocLyp} are thus met.


\subsubsection{Proof of Corollary \ref{mixingex}{\rm (2)}}
\label{S7.1.2}

Let $\xi$ be the zero-range process as described in Corollary~\ref{mixingex}{\rm (2)}. We 
will use that each particle, independently of all the other particles, carries an exponential 
clock of parameter one. If there are $k$ particles at a site $x$ and one of these clocks rings,
then the corresponding particle jumps to $y$ with probability $\frac{g(k)}{2dk}$ and it stays 
at $y$ with probability $1-\frac{g(k)}{k}$.

\medskip\noindent
{\bf 1.1}  
By Andjel~\cite{A82}, Theorem 1.9, the product measures in (\ref{ZRmeasure}) are extremal for 
$\xi$. Thus, $\E\left[e^{q \xi(0,0)}\right] <\infty$ for all $q\geq 0$. Consequently, to show 
that $\E\left[e^{q \sup_{s\in[0,1]}\xi(0,s)}\right] <\infty$, it suffices to prove that there 
is a constant $K>0$ such that, for all $k\in\N$ sufficiently large,
\begin{equation}
\label{supest}
\P\Bigg(\sup_{s\in[0,1]} \xi(0,s) \geq k\Bigg) 
\leq K\P\Bigg(\xi(0,1)\geq \frac{e^{-1}}{2} k\Bigg).
\end{equation}
Write $\mathrm{NR}(\frac{e^{-1}}{2}k,\tau)$ for the event that there are at least 
$\frac{e^{-1}}{2}k$ exponential clocks of particles located at zero that do not ring 
in the time interval $[\tau,\tau+1)$. Then we may estimate
\begin{equation}
\label{condprob1ex}
\begin{aligned}
&\P\Bigg(\xi(0,1)\geq \frac{e^{-1}}{2}k ~\Big|~ 
\exists\, \tau\in[0,1]\colon\,\xi(0,\tau)\geq k\Bigg)\\
&\qquad \geq \P\Bigg(\mathrm{NR}\bigg(\frac{e^{-1}}{2}k,\tau\bigg)~\Big|~
\exists\, \tau\in[0,1]\colon\,\xi(0,\tau)\geq k\Bigg).
\end{aligned}
\end{equation}
Since the probability that a clock does not ring within a time interval of length one is
equal to $e^{-1}$, and all clocks are independent, we may estimate the right-hand side of 
(\ref{condprob1ex}) from below by
\begin{equation}
\label{lbBIN}
\P\Big(T\geq \frac{e^{-1}}{2}k\Big),\qquad T = \mathrm{BIN}(k,e^{-1}).
\end{equation}
Finally, note that the probability in (\ref{lbBIN}) is bounded away from zero. Thus,
inserting (\ref{condprob1ex}) and (\ref{lbBIN}) into (\ref{supest}), we get the claim.

\medskip\noindent
{\bf 1.2.} 
We show that $\xi$ is type-I G\"artner-mixing. Fix $A>3$, choose $b=3$, $c=1$, $a_1=13$, 
$a_2=2$ (see (\ref{equiv}), and introduce additional space-time blocks
\begin{equation}
\label{addblockex}
\begin{aligned}
&\tilde{B}_{R}^{A,+}(x,k) 
= \Bigg(\prod_{j=1}^{d}[(x(j)-2)A^{R},(x(j)+2)A^{R})\cap\Z^d\Bigg)
\times [kA^R - A^{R-1}, (k+1)A^R)\\
&\tilde{B}_{R}^{A,\mathrm{sub}}(x,k) 
= \Bigg(\prod_{j=1}^{d}[(x(j)-4)A^{R},(x(j)+4)A^{R})\cap\Z^d\Bigg)
\times\{(k-1)A^R\}.
\end{aligned}
\end{equation}
Given $S \subseteq\Z^d$ and $S'\subseteq\N_0$, write $\partial S$ to denote the inner boundary of $S$
and $\Pi_1(S\times S')$ to denote the projection onto the spacial coordinates. Furthermore, 
$e_j$, $j\in\{1,2,\cdots, d\}$ denotes the $j$-th unit vector (and we agree that $\frac{0}{0}
=0$).

We call a space-time block $B_{R}^{A}(x,k)$ \emph{contaminated} if there is a particle at some 
space-time point $(y,s)\in \tilde{B}_{R}^{A,+}(x,k)$ that has been outside $\Pi_1(\prod_{j=1}^{d}
[(x(j)-4)A^{R}, (x(j)+4)A^{R}))$ at a time $s'$ such that $(k-1)A^R\leq s'< s$. The reason 
for introducing this notion is that events depending on non-contaminated blocks that are equal 
in time but disjoint in space are all independent.

\medskip\noindent
{\bf Contaminated blocks.} 
For $L>0$, define
\begin{equation}
\label{chiZRex}
\chi(x,k)= \one\Big\{B_{R+1}^{A}(x,k) \mbox{ is good, 
but contaminated and intersects } [-L,L]^{d+1}\Big\}
\end{equation}
and fix $(x^*,k^*)\in\Z^d\times\N$.

\bcl{ZRclaim}
There is a $C'>0$ independent of $L$ such that $(\chi(x,k))_{(x,k) \equiv (x^*,k^*)}$ is 
stochastically dominated by independent Bernoulli random variables $(Z(x,k))_{(x,k)=(x^*,k^*)}$ 
with success probability $e^{-C'A^{R+1}}$.
\ecl

\bpr
We use a discretization scheme. More precisely, we construct a discrete-time version of 
the zero-range process where particles are allowed to jump at times $k/n$, 
$k \in\N_0$ only. Here, $n$ is an integer that will later tend to infinity, and we will 
denote by $\xi^n(x,s)$ the number of particles at site $x$ at time $s$. To construct this 
process, we take a family $X^{n}(x,s,q_1,q_2)$ of independent random variables with index set 
$\Z^d\times\frac{1}{n}\N_0\times\N_0\times\N_0$ whose distribution is defined via
\begin{equation}
\label{XnZRex}
\begin{aligned}
&\P\Big(X^n(\cdot,\cdot,\cdot,q_2)=0\Big) = 1-\frac{g(q_2)}{nq_2},\\
&\P\Big(X^n(\cdot,\cdot,\cdot,q_2)=\pm e_j\Big) = \frac{g(q_2)}{2dnq_2},
\qquad j\in\{1,2,\ldots,d\}.
\end{aligned}
\end{equation}
With this family in hand, we proceed as follows. At time zero start with an initial configuration 
that comes from the invariant measure $\pi_{\rho}$. Attach to each particle $\sigma$ a uniform-$[0,1]$ 
random variable $\mathcal{U}(\sigma)$. Take all these random variables independent of each other and 
of $X^{n}(x,s,q_1,q_2)$ for all choices of $(x,s,q_1,q_2)\in\Z^d\times\frac{1}{n}\N_0\times\N_0\times\N_0$. 
For each site $x$, order all particles present at $x$ at time zero so that their uniform random variables 
are increasing. To the $q_1$-th variable attach $X^n(x,0,q_1,\xi^n(x,0))$, i.e., the position of the 
$q_1$-th particle in this ordering at time $\frac{1}{n}$ is $x+X^n(x,0,q_1,\xi^n(x,0))$. In this way 
we obtain the configuration of the system at time $\frac{1}{n}$. To construct the process 
$\frac{1}{n}$ time units further, repeat the first step, but let the particles jump according 
to $X^n(\cdot,\frac{1}{n},\cdot,\xi^n(\cdot,\frac{1}{n}))$. Thus, our construction is such 
that each particle chooses at each step uniformly at random, but dependent on the number of 
particles at the same location, a new jump distribution. In what follows we will use the 
phrase ``at level $n$'' to emphasize that we refer to the discrete-time version of the process. 
For instance, we say that $B_{R}(x,k)$ is good at level $n$ if
\begin{equation}
\label{goodlevelnex}
\sum_{z\in Q_{R}^{A}(y)} \xi^n(z,s)\leq CA^{Rd} \qquad \forall\, y\in\Z^d, 
s\geq 0\, \mbox{ s.t. } Q_{R}(y)\times\{s\} \subseteq \tilde{B}_{R}^{A}(x,k).
\end{equation} 

Next, we introduce 
\begin{equation}
\label{chinZRex}
\begin{aligned}
&\chi^{n}(x,k)\\
&= \one\Big\{B_{R+1}^{A}(x,k) 
\mbox{ is good, but is contaminated at level $n$ and intersects } [-L,L]^{d+1}\Big\}.
\end{aligned}
\end{equation}
It is not hard to show that the joint distribution of $\chi^n$ converges weakly to the joint 
distribution of $\chi$ (use that only finitely many particles can enter a fixed region in 
space-time, so that the above family of random variables may be approximated by a function 
depending on finitely many particles only). Thus, to estimate the joint distribution of $\chi$, 
it is enough to analyze the joint distribution of $\chi^n$, as long as the estimates are uniform 
in $n$. In what follows, $s,s',s''\in\frac{1}{n}\N_0$.

Let $B_{R+1}^{A}(x,k)$ be a good block that is contaminated at level $n$. Then there is a 
particle at a site $y\in\partial \Pi_1(\tilde{B}_{R+1}^{A}(x,k))$ at a time $s\in[(k-1)A^{R+1},
(k+1)A^{R+1})$ (see below \ref{spacetimeblock})) that is at a site $y'\in\partial\Pi_1
(\tilde{B}_{R+1}^{A,+}(x,k))$ at a time $s'\in[(k-1)A^{R+1},(k+1)A^{R+1})$, $s<s'$.
Furthermore, for all $s''$ such that $s<s''<s'$, the particle is inside $\Pi_1
(\tilde{B}_{R+1}^{A}(x,k))$. This implies, when $X^n(y,s,q_1,\xi^{n}(y,s))$ denotes the 
random variable attached to $\sigma$ at time $s$, that $X^n(y,s,q_1,\xi^{n}(y,s))\neq 0$.
Pick any such particle. Since this particle travels over a distance larger than $2A^{R+1}$,
there is at least one coordinate direction along which it makes at least $2A^{R+1}$ steps. 
We call this direction $e_j(\sigma)$, and say that each step in this direction is a success.
Note the uniform estimate 
\begin{equation}
\label{successex}
\P\big(\sigma \mbox{ has a success at time $s''$}\big) \leq \frac{1}{2dn}
\qquad s''\in [(k-1)A^{R+1},(k+1)A^{R+1}).
\end{equation}
Thus, if $\sigma$ contaminates $B_{R+1}^{A}(x,k)$ at level $n$, then from time $s'$ up to 
time $(k+1)A^{R+1}$ it has at least $2A^{R+1}$ successes in $\Pi_1(\tilde{B}^{A}_{R+1}(x,k))$.
We write $\mathrm{S}(y,s,q_1,q_2)$ for the event just described, provided the particle was 
attached to $X^n(y,s,q_1,q_2)$ when it entered $\tilde{B}_{R+1}^{A}(x,k)$. Since 
$B_{R+1}^{A}(x,k)$ is good, at each space-time point $(y'',s'')\in\tilde{B}_{R+1}^{A}(x,k)$
there are at most $CA^{(R+1)d}$ particles that can contaminate $B_{R+1}^{A}(x,k)$. We 
therefore obtain
\begin{equation}
\label{contZRex}
\begin{aligned}
&\Big\{B_{R+1}^{A}(x,k) \mbox{ is good, but contaminated at level $n$}\Big\}\\
&\subseteq \bigcup_{\substack{y\in\partial\Pi_1(\tilde{B}_{R+1}^{A}(x,k))\\
s\in[(k-1)A^{R+1},(k+1)A^{R+1}),\, s\in\frac{1}{n}\N_0\\
q_1\leq CA^{(R+1)d}\\
q_2\leq CA^{(R+1)d}}}
\Big\{X^n(y,s,q_1,q_2)\neq 0, \mathrm{S}(y,s,q_1,q_2)\Big\}\\
&\stackrel{\mathrm{def}}= \mathrm{C}^n(x,k).
\end{aligned}
\end{equation}

Next, note that the event $\mathrm{C}^n(x,k)$ depends on the $X^n(y,s,q_1,q_2)$ with $(y,s)\in
\tilde{B}_{R+1}^{A}(x,k)$ only. Hence, for $x\in\Z^d$, $k\in\N_0$, the family 
$(\mathcal{C}^n(x,k))_{(x,k)\equiv(x^*,k^*)}$ consists of independent events.
We estimate $\P(C^n(x,k))$. By (\ref{successex}), the probability inside the union may 
be bounded from above by 
\begin{equation}
\label{successest1ex}
\frac{1}{n}\P\big(T\geq 2A^{R+1}\big),
\end{equation}
where $T=\mathrm{BIN}(2A^{R+1}n,\frac{1}{2dn})$. Note that the event in (\ref{successest1ex}) 
is a large deviation event, so that Bernstein's inequality guarantees the existence of a 
constant $C'>0$ such that (\ref{successest1ex}) is at most
\begin{equation}
\label{successest2ex}
\frac{1}{n}\exp\left\{-C'A^{R+1}\right\}.
\end{equation}
Recall the definition of $C^{n}(x,k)$ to see that, for a possibly different constant $C''>0$,
\begin{equation}
\label{Cestex}
\P(C^{n}(x,k))\leq \exp\left\{-C''A^{R+1}\right\}.
\end{equation}
Hence there is a family of independent Bernoulli random variables $(Z^n(x,k))_{(x,k)\equiv
(x^*,k^*)}$ that stochastically dominates $(\chi^n(x,k))_{(x,k)\equiv(x^*,k*)}$ and has success 
probability $e^{-C''A^{R+1}}$. Thus, if $f$ is a positive and bounded function that is 
increasing in all of its arguments, then
\begin{equation}
\label{coupl1ex}
\E\Big(f(\chi^n(x_1,k_1),\ldots, \chi^n(x_n,k_n)\Big)
\leq \E\Big(f(Z^n(x_1,k_1),\ldots, Z^n(x_n,k_n)\Big).
\end{equation}
As the right-hand side does not depend on $n$, we obtain, by letting $n\to \infty$,
\begin{equation}
\label{coupl2ex}
\E\Big(f(\chi(x_1,k_1),\ldots, \chi(x_n,k_n)\Big)
\leq \E\Big(f(Z(x_1,k_1),\ldots, Z(x_n,k_n)\Big),
\end{equation}
which proves Claim \ref{ZRclaim}. In particular, since all estimates are independent of 
$L$, we may even set $L=\infty$, to get that the whole field of good but contaminated 
$(R+1)$-blocks may be dominated by an independent family of Bernoulli random variables
with the same success probability as above. This field will be denoted by $Z(x,k)$ as well.
\epr

\medskip\noindent
{\bf Non-contaminated blocks.} 
We begin by estimating the probability of the event
\begin{equation}
\label{ZRgoodbad}
\Big\{B_{R+1}^{A}(x,k) \mbox{ is good, but contains a bad $R$-block}\Big\}.
\end{equation}
Let $B_{R}^{A}(y,l)$ be an $R$-block that is contained in $B_{R+1}^{A}(x,k)$. We bound 
the probability that this block is bad. For that we take $z_1\in\Z^d$ such that
\begin{equation}
\label{ZRQblock}
Q_{R}^{A}(z_1)=\prod_{j=1}^{d}[z_1(j),z_1(j)+A^{R}) \subseteq \tilde{B}_{R}^{A}(y,l).
\end{equation}
Use the time stationarity of $\xi$ and the fact that $[(l-1)A^{R},(l+1)A^{R})$ may be 
divided into at most $2A^{R}+1$ time intervals of length one, to obtain
\begin{equation}
\label{ZRsupestex}
\begin{aligned}
&\P\Bigg(\exists s\in[(l-1)A^{R},(l+1)A^{R})\colon\, 
\sum_{z_2\in Q_{R}^{A}(z_1)} \xi(z_2,s) >CA^{Rd}\Bigg)\\
&\leq (2A^{R}+1) \P\Bigg(\sup_{s\in[0,1]} \sum_{z_2\in Q_{R}^{A}(z_1)}
\xi(z_2,s)>CA^{Rd}\Bigg).
\end{aligned}
\end{equation}
In the same way as in Step 1.1, we may now show that for some constant $K>0$,
\begin{equation}
\label{ZRsupsumest}
\P\Bigg(\sup_{s\in[0,1]}\sum_{z_2\in Q_{R}^{A}(z_1)} \xi(z_2,s) > CA^{Rd}\Bigg)
\leq K\P\Bigg(\sum_{z_2\in Q_{R}^{A}(z_2)} \xi(z_2,1)\geq \frac{e^{-1}}{2}CA^{Rd}\Bigg).
\end{equation}
Next, note that under the invariant measure $\pi_{\rho}$ the sum in the right-hand side 
of (\ref{ZRsupsumest}) is a sum of i.i.d.\ random variables with finite exponential moments. 
Hence (\ref{ZRsupsumest}) is bounded from above by
\begin{equation}
\label{ZRsupsumbound}
K\exp\Big\{-\frac{e^{-1}}{2}CA^{Rd}\Big\} \E\left[e^{\xi(0,0)}\right]^{A^{Rd}}.
\end{equation}
Now choose $C$ large enough so that (\ref{ZRsupestex}) decays superexponentially fast 
in $R$.

In order to estimate the joint distribution of non-contaminated blocks, we let 
$B_{R+1}^{A_1}(x_1,k_1),$ $\dots,B_{R+1}^{A_1}(x_n,k_n)$ be space-time blocks whose 
indices increase in the lexicographic order of $\Z^d\times\N$ and belong to the same 
equivalence class. We abbreviate
\begin{equation}
\label{NoncontZRex}
\begin{aligned}
&N^{\mathrm{sub}}(x_i,k_i) = \Big\{\tilde{B}_{R+1}^{A,\mathrm{sub}}(x_i,k_i) \mbox{ is good, }
B_{R+1}^{A}(x_i,k_i) \mbox{ contains a bad $R$-block,}\\
&\qquad\qquad\qquad\ \ \ \ \mbox{but is not contaminated}\Big\}.
\end{aligned}
\end{equation}
Note that $N^{\mathrm{sub}}(x_i,k_i)$ and $N^{\mathrm{sub}}(x_j,k_j), i\neq j$, are independent 
when they depend on blocks that coincide in time but are disjoint in space. This observation,
together with the Markov property applied, leads to
\begin{equation}
\label{sGartnerncZRex}
\begin{aligned}
&\P\Bigg(N^{\mathrm{sub}}(x_n,k_n)
~\Big|~ \bigcap_{i=1}^{n-1}N^{\mathrm{sub}}(x_i,k_i)\Bigg)\\
&\leq \P\Bigg(N^{\mathrm{sub}}(x_n,k_n)
~\Big| \bigcap_{i=1}^{n-1}N^{\mathrm{sub}}(x_i,k_i), \tilde{B}_{R+1}^{A,\mathrm{sub}}(x_n,k_n)
\mbox{ is good}\Bigg)\\
&= \P\Big(N^{\mathrm{sub}}(x_n,k_n)~\Big|~  \tilde{B}_{R+1}^{A,\mathrm{sub}}(x_n,k_n)
\mbox{ is good}\Big).
\end{aligned}
\end{equation}
Thus, the left-hand side of (\ref{sGartnerncZRex}) is at most
\begin{equation}
\label{condupboundZRex}
\frac{\P(B_{R+1}^{A}(x,k) \mbox{ is good, but contains a bad $R$-block})}
{\P(\tilde{B}_{R+1}^{A_1,\mathrm{sub}}(x_n,k_n) \mbox{ is good})}.
\end{equation}
Note that the denominator tends to one as $R\to\infty$. This comes from the fact that, 
for all $t\geq 0$, $(\xi(x,t))_{x\in\Z^d}$ is an i.i.d.\ field of random variables 
distributed according to $\pi_{\rho}$. Thus, from (\ref{ZRsupestex}) and the lines below, 
we infer that (\ref{condupboundZRex}) decays superexponentially fast in $R$.

Finally, write
\begin{equation}
\label{divconZRex}
\Big\{B_{R+1}^{A}(x_i,k_i) \mbox{ is good, but contains a bad $R$-block }\Big\}
\subseteq C(x_i,k_i) \cup N^{\mathrm{sub}}(x_i,k_i),
\end{equation}
where we denote by $C(x_i,k_i)$ the event that $B_{R+1}^{A}(x_i,k_i)$ is good,
contains a bad $R$-block, and is contaminated. Then
\begin{equation}
\label{intersectZRex1}
\begin{aligned}
&\P\Bigg(\bigcap_{i=1}^{n}\left\{B_{R+1}^{A_1}(x_i,k_i) 
\mbox{ is good, but contains a bad $R$-block}\right\}\Bigg)\\ 
&\qquad \leq \P\Bigg(\bigcap_{i=1}^{n} \big(C(x_i,k_i)\cup 
N^{\mathrm{sub}}(x_i,k_i)\big)\Bigg).
\end{aligned}
\end{equation}
If we denote by $C$ the subset of all $i\in\{1,2,\ldots, n\}$ for which $C(x_i,k_i)$ occurs,
then (\ref{intersectZRex1}) may be rewritten as
\begin{equation}
\label{intersectZRex2}
\sum_{C\subseteq \{1,2,\ldots,n\}} \P\Bigg(\bigcap_{i\in C} C(x_i,k_i) 
\cap \bigcap_{i\notin C} N^{\mathrm{sub}}(x_i,k_i)\Bigg).
\end{equation}
Note that either $|C|\geq n/2$ or $|\{1,2,\ldots,n\}\setminus C|\geq n/2$, so that
by Claim \ref{ZRclaim} and (\ref{ZRsupsumbound}) there is a $C''>0$ such that the 
expression in (\ref{intersectZRex2}) is at most $2^n\exp\{-C''A^{R+1}n/2\}$. A 
comparison with the right-hand side of (\ref{wGartner}) shows that $\xi$ is type-I 
G\"artner-mixing.

\medskip\noindent
{\bf 1.3} 
From the previous calculations we infer that $\xi$ satisfies condition {\rm (a3)} of 
Definition~\ref{Gartnerhypmix}.

\medskip\noindent
{\bf 2.}
Since $\xi$ is bounded from below it is G\"artner-negative-hyper-mixing. Hence it 
remains to show that $\xi$ is type-II G\"artner-mixing. To that end, fix the same 
constants as in the proof of the first part of the corollary. Furthermore, fix 
$\delta >0$ and let $B_{R+1}^{A}(x_1,k_1),\ldots, B_{R+1}^{A}(x_n,k_n)$ be space-time 
blocks whose indices increase in the lexicographic order of $\Z^d\times\N$, and belong 
to the same equivalence class. Take events $\mathcal{A}_i^{R} \in \sigma(B_{R+1}^{A}(x_n,k_n))$, 
$i\in\{1,2,\ldots,n\}$, that are invariant under shifts in space and time, and satisfy
\begin{equation}
\label{AtozeroZRex}
\lim_{R\to\infty}
\P(\mathcal{A}_i^{R})=0.
\end{equation}
As in part 1, we divide space-time blocks into contaminated and 
non-contaminated blocks. With the help of Claim~\ref{ZRclaim}, we may control 
contaminated blocks. To treat non-contaminated blocks, we introduce
\begin{equation}
\label{ZRNCex}
N^{\mathrm{sub}}(x_i,k_i,\mathcal{A}_i^{R}) = 
\Big\{\tilde{B}_{R+1}^{A,\mathrm{sub}}(x_i,k_i) \mbox{ is good, but contaminated, }
\mathcal{A}_i^{R} \mbox{ occurs}\Big\}
\end{equation}
and proceed as in the lines following (\ref{NoncontZRex}), to finish the proof.


\subsection{Proof of Corollary \ref{contex}}
\label{S7.2}

In this section we prove Corollary~\ref{contex}. Suppose that
\begin{equation}
\label{xiMarkovcond}
\begin{aligned}
\blacktriangleright\quad &\xi \mbox{ is Markov with initial distribution } 
\nu \mbox{ and}\\[-.15cm]
&\mbox{generator } L \mbox{ defined on a domain } D(L) \subset L^2(d\nu).
\end{aligned}
\end{equation}
Denote by 
\begin{equation}
\label{symDirichlet}
\bar\varepsilon(f,g) = \tfrac12\big[\langle -L f,g \rangle 
+ \langle -L g,f \rangle\big], 
\qquad f,g \in D(L),
\end{equation}
its symmetrized Dirichlet form, assume that $(\bar\varepsilon,D(L))$ is closable, 
and denote its closure by $(\bar\varepsilon, D(\bar\varepsilon))$. Furthermore, for 
$V\colon\,\Omega \to \R$, define
\begin{equation}
\label{LDPnonreg}
J_V(r) = \inf\left\{\bar\varepsilon(f,f) \colon\, 
f \in D(\bar\varepsilon)\cap L^2(|V|d\nu),\, \int f^2\,d\nu =1,\, 
\int Vf^2\,d\nu = r\right\}, \qquad r\in\R,
\end{equation}
note that $r \mapsto J_V(r)$ is convex, and let $I_V$ be its lower semi-continuous 
regularization. 
For $W\colon\,\Omega\times[0,\infty) \to \R$, define
\begin{equation}
\label{GammacL}
\Gamma_t^L = \sup_{\substack{f \in D(L) \cap L^2(|W(\cdot,t)|d\nu)\\ \|f\|_2=1}}
\left\{\langle L f, f \rangle + \langle W(\cdot,t), f^2 \rangle \right\},
\qquad t\geq 0.
\end{equation}
In the particular case of a static $W_0\colon\,\Omega\to\R$, we denote the corresponding 
variational expression by $\Gamma_0^L$. 

\bl{Cor118Lem} 
Suppose that $W$ is bounded from below, piecewise continuous in the time-coordi\-nate, 
and $\nu$-integrable in the space-coordinate. Suppose further
that:\\
(i) $\Gamma_t^L< \infty$ for all $t>0$.\\
(ii) $\xi$ is reversible in time and c\`adl\`ag.\\
Then 
\begin{equation}
E_\nu\left(\exp\left\{\int_0^t W(\xi(s),s)\,ds\right\}\right) 
\leq \exp\left\{\int_0^t \Gamma_s^L\,ds\right\} \qquad \forall\,t\geq 0.
\end{equation}
\el
\bpr
The proof is based on ideas in Kipnis and Landim~\cite{KL99}, Appendix 1.7, and comes in 
three steps.

\medskip\noindent
{\bf 1.} 
Suppose that $W$ is piecewise continuous, not necessarily bounded from above, 
and define $W^n = W\wedge n$. Then, by an argument similar to that in \cite{KL99}, 
Appendix 1, Lemma 7.2, we have
\begin{equation}
\E_\nu\left(\exp\left\{\int_0^t W^n(\xi(s),s)\,ds\right\}\right) 
\leq \exp\left\{\int_0^t \Gamma^{L,n}_s\,ds\right\},
\end{equation} 
where $\Gamma^{L,n}$ is defined as in (\ref{GammacL}) but with $W$ replaced by $W^n$. It is 
here that we use that $\xi$ is reversible and c\`adl\`ag, since under this condition we have a Feynman-Kac
representation for the parabolic Anderson equation with potential $W^n$ and with $\Delta$ 
replaced by $L$. Because $W^n$ is bounded, we have that $L^2(|W^n(\cdot,t)|d\nu) = L^2(d\nu)$. 
Hence
\begin{equation}
\Gamma^{L,n}_s 
= \sup_{\substack{f \in D(L)\\ \|f\|_2 = 1}}
\left\{\langle Lf,f \rangle + \langle W^n(\cdot,s),f^2\rangle\right\}.
\end{equation}
Since, by monotone convergence,
\begin{equation}
\lim_{n \to \infty} \E_\nu\left(\exp\left\{\int_0^t W^n(\xi(s),s)\,ds\right\}\right) 
= \E_\nu\left(\exp\left\{\int_0^t W(\xi(s),s)\,ds\right\}\right),
\end{equation}
it suffices to show that
\begin{equation}
\label{Gamlim}
\lim_{n \to \infty} \Gamma^{L,n}_s = \Gamma^{L}_s, \qquad s \geq 0.
\end{equation}

\medskip\noindent
{\bf 2.} 
To show that the left-hand side in (\ref{Gamlim}) is an upper bound, fix $\varepsilon>0$ 
and pick $f \in D(L)\cap L^2(|W(\cdot,t)|d\nu)$ such that 
\begin{equation}
\langle L f,f \rangle + \langle W(\cdot,t),f^2 \rangle + \varepsilon 
\geq \Gamma^{L}_s.
\end{equation}
Then, by monotone convergence, 
\begin{equation}
\lim_{n \to \infty} \left\{\langle L f,f \rangle 
+ \langle W^n(\cdot,t), f^2 \rangle\right\} 
= \langle L f,f \rangle + \langle W(\cdot,t), f^2 \rangle.
\end{equation}

\medskip\noindent
{\bf 3.} 
To prove that the left-hand side in (\ref{Gamlim}) is also a lower bound,
we need to assume that $L$ is self-adjoint. Then, by Wu~\cite{Wu94}, Remark on 
p.\ 209, we have
\begin{equation}
\lim_{t\to \infty} \frac{1}{t} \log 
\E_\nu\left(\exp\left\{\int_0^t W^n(\xi(s),0)\,ds\right\}\right) 
= \Gamma_0^{L,n},
\end{equation}
and the same is true when $W^n$ is replaced by $W$. But, obviously,
\begin{equation}
\lim_{t\to\infty} \frac{1}{t} \log 
\E_\nu\left(\exp\left\{\int_0^t W^n(\xi(s),0)\,ds\right\}\right) 
\leq \lim_{t \to \infty} \frac{1}{t} \log 
\E_\nu\left(\exp\left\{\int_0^t W(\xi(s),0)\,ds\right\}\right),  
\end{equation}
which shows that $\Gamma_0^{L,n} \leq \Gamma_0^{L}$. Note that the time point 0 does not 
play any special role. Hence we obtain (\ref{Gamlim}).
\epr

\bp{contDirichlet}
Suppose that $\xi$ satisfies \eqref{xiMarkovcond} and condition (ii) in 
Lemma~{\ref{Cor118Lem}}. Then, for all paths $\Phi$ (recall \eqref{path}), 
\begin{equation}
\label{devinequ}
\P_\nu\left(\frac{1}{n(T-1)} \sum_{j=1}^n \int_{(j-1)T+1}^{jT} 
\big[\xi(\Phi(s),s)\,ds - \rho\big] > \delta\right) 
\leq \exp\big\{-n(T-1)I_{V_0}(\rho+\delta)\big\},
\end{equation}
where $V_0(\eta) = \eta(0)$ and $\rho=\E(\xi(0,0))$.
\ep

\bpr
First note that by Lemma~\ref{Cor118Lem}, 
\begin{equation}
\label{contDirichlet1}
\E_\nu\left(\exp\left\{\sum_{j=1}^n \int_{(j-1)T+1}^{jT} 
\xi(\Phi(s),s)\,ds\right\}\right) 
\leq \exp\big\{n(T-1)\Gamma_0^L\big\},
\end{equation}
where $\Gamma_0^L$ is defined with $W(\cdot,0)=V_{0}$.  
To see that, take a path $\Phi$ which starts in zero and 
define $W\colon\,\R^{\Z^d} \times [0,\infty) \to \R$ by
\begin{equation}
W(\eta,t) =
\begin{cases} 
\eta(\Phi(t)),  
& \mbox{if } t \in [(j-1)T+1,jT) \mbox{ for some } 1\leq j\leq n,\\
0, & \mbox{otherwise}.
\end{cases}
\end{equation}
Let $x \in \Z^d$ be such that $\Phi(t) = x$, and denote by $\tau_x$ the space-shift 
over $x$. Then, using the fact that $\nu$ is shift-ergodic, we get
\begin{equation}
\begin{aligned}
\langle W(\cdot,t),f^2 \rangle 
&= \int_{\R^{\Z^d}}\eta(x)f^2(\eta)\,d\nu(\eta)\\
&= \int_{\R^{\Z^d}}\tau_x\eta(0)f^2(\eta)\,d\nu(\eta)\\
&= \int_{\R^{\Z^d}}\eta(0)(\tau_{-x}f)^2(\eta)\,d\nu(\eta)
= \langle W(\cdot,0),(\tau_{-x}f)^2 \rangle,
\end{aligned}
\end{equation}
which yields $\Gamma_t^L=\Gamma_0^L$. 
Since the space point $0$ does not play
any special role, Lemma \ref{Cor118Lem} leads to (\ref{contDirichlet1})
for any path $\Phi$.
Next apply the Chebyshev inequality to the left-hand side of (\ref{devinequ}). 
After that it remains 
to solve an optimization problem. See Wu~\cite{Wu00} for details.
\epr

\noindent
Note that, for a reversible dynamics, $I_{V_0}$ is the large deviation rate 
function for the occupation time
\begin{equation}
\label{Tt}
T_t = \int_0^t \xi(0,s)\,ds, \qquad t \geq 0.
\end{equation}

We are now ready to give the proof of Corollary~\ref{contex}.

\bpr
All three dynamics in (1)--(3) satisfy condition (a) in Definition~\ref{Gartnerposreg}. 
The proof of (b) below consists of an application of Proposition~\ref{contDirichlet}, 
combined with a suitable analysis of $I_{V_0}$. 

\medskip\noindent
(1) Redig and V\"ollering~\cite{RV11}, Theorem 4.1, shows that for all $\delta_1 >0$ 
there is a $\delta_2=\delta_2(\delta_1)$ such that
\begin{equation}
\label{Redig}
\P_{\nu}\left(\int_0^{nt}\xi(0,s)\,ds \geq \delta_1 nt\right) \leq e^{-\delta_2 nt}.
\end{equation}
A straightforward extension of this result implies that condition (b) in 
Definition~\ref{Gartnerposreg} is satisfied. All requirements in 
Theorem~\ref{continuity} are thus met.

\medskip\noindent
(2) By Landim~\cite{L92}, Theorem 4.2, the rate function of the simple exclusion process 
is non-degenerate (i.e., it has a unique zero at $\rho$). Hence condition (b) in 
Definition~\ref{Gartnerposreg} is satisfied. Thus, all requirements of Theorem~\ref{continuity} are met.

\medskip\noindent
(3) By Cox and Griffeath~\cite{CG84}, Theorem 1, the rate function for independent
simple random walks is non-degenerate. Hence condition (b) in Definition~\ref{Gartnerposreg} 
is satisfied. All requirements in 
Theorem~\ref{continuity} are thus met.
\epr




\begin{thebibliography}{99}

\bibitem{A82}
E. D.\ Andjel,
Invariant measure for the zero range process,
Ann.\ Probab.\ 10 (1982) 525--547. 

\bibitem{CM94}
R.A.\ Carmona and S.A.\ Molchanov, 
\emph{Parabolic Anderson Problem and Intermittency},
AMS Memoirs 518, American Mathematical Society, Providence RI, 1994. 

\bibitem{CG84}
J.T.\ Cox and D.\ Griffeath,
Large deviations for Poisson systems of independent random walks,
Z.\ Wahrsch.\ Verw.\ Gebiete 66 (1984) 543--558.

\bibitem{CMS02}
M.\ Cranston, T.S.\ Mountford and T.\ Shiga,
Lyapunov exponents for the parabolic Anderson model,
Acta Math. Univ. Comenianae\ 71 (2002) 163--188.

\bibitem{DGRS12}
A.\ Drewitz, J.\ G\"artner, A.F.\ Ramirez and R.\ Sun,
Survival probability of a random walk among a Poisson system of moving traps, 
in: \emph{Probability in Complex Physical Systems. In honour of Erwin Bolthausen 
and J\"urgen G\"artner} (eds.\ J.-D.\ Deuschel, B.\ Gentz, W.\ K\"onig, M.-K.\ 
van Renesse, M.\ Scheutzow, U.\ Schmock), Springer Proceedings in Mathematics 11, 
Springer, 2012, Berlin, pp.\ 119--158.

\bibitem{E42}
P.\ Erd\"os,
On an elementary proof of some asymptotic formulas
in the theory of partitions,
Ann.\ Math.\ 43 (1942) 437--450.

\bibitem{GdH06}
J.\ G\"artner and F.\ den Hollander, 
Intermittency in a catalytic random medium,
Ann.\ Probab.\ 34 (2006) 2219--2287.

\bibitem{GdHM11}
J.\ G\"artner, F.\ den Hollander and G.\ Maillard, 
Quenched Lyapunov exponent for the parabolic Anderson model in a dynamic 
random environment, in: \emph{Probability in Complex Physical Systems.
In honour of Erwin Bolthausen and J\"urgen G\"artner} (eds.\ J.-D.\ Deuschel, 
B.\ Gentz, W.\ K\"onig, M.-K.\ van Renesse, M.\ Scheutzow, U.\ Schmock), 
Springer Proceedings in Mathematics 11, Springer, 2012, Berlin, pp.\ 159--193.

\bibitem{GM90}
J.\ G\"artner and S.A.\ Molchanov,
Parabolic problems for the Anderson model,
Comm.\ Math.\ Phys.\ 132 (1990) 613--655.
 
\bibitem{HR18}
G.H.\ Hardy and S.A.\ Ramanujan,
Asymptotic formulae in combinatory analysis,
Proc.\ London\ Math.\ Soc.\ 17 (1918) 75--115.
 
\bibitem{KL99}
C.\ Kipnis and C.\ Landim, 
\emph{Scaling Limits of Interacting Particle Systems},
Grundlehren der Mathematischen Wissenschaften 320, 
Springer, Berlin, 1999.
 
\bibitem{KS03}
H.\ Kesten and V.\ Sidoravicius,
Branching random walks with catalysts,
Electr.\ J.\ Prob.\ 8 (2003) 1--51.

\bibitem{L92}
C.\ Landim,
Occupation time large deviations for the symmetric simple exclusion process,
Ann.\ Probab.\ 20 (1992) 206--231.

\bibitem{L85}
T.M.\ Liggett, 
\emph{Interacting Particle Systems}, 
Grundlehren der Mathematischen Wissenschaften 276, 
Springer, New York, 1985.


\bibitem{RV11}
F.\ Redig and F.\ V\"ollering,
Concentration of additive functionals for Markov processes, 
preprint 2011.
\url{http://arXiv:1003.0006v2} 

\bibitem{Wu94}
L.\ Wu, 
Feynman-Kac semigroups, ground state diffusions, and large deviations,
J.\ Funct.\ Anal.\ 123 (1994) 202--231.

\bibitem{Wu00}
L.\ Wu,
A deviation inequality for non-reversible Markov processes,
Ann.\ Inst.\ H.\ Poincar\'e: Probab.\ Statist.\ 4 (2000) 435--445.

\end{thebibliography}
\end{document}